# Sur le théorème de l'indice des équations différentielles $p$-adiques III

Par G. Christol et Z. Mebkhout

## Table des matières

## 1. Introduction

Cet article fait suite aux articles [C-M$_1$], [C-M$_2$]. Nous poursuivons le programme de Philippe Robba ([R$_1$], [R$_2$], [R$_3$], [R$_4$], [R$_5$], [R$_6$], [R$_7$], [R-C]) sur le théorème de l'indice pour les équations différentielles $p$-adiques.

Dans les articles précédents nous avons défini les exposants de la monodromie locale $p$-adique d'un module différentiel ayant la propriété de Robba dans une couronne, montré le théorème de décomposition en modules de rang un dans le cas modéré, montré le théorème d'existence de type de Riemann et défini la catégorie des coefficients $p$-adiques modérément ramifiés sur les courbes. Dans cet article nous montrons le théorème de décomposition selon les pentes $p$-adiques, pour les modules différentiels solubles, sur l'anneau des fonctions analytiques au bord. Ceci termine la détermination de la structure d'un point singulier d'une équation différentielle $p$-adique et ses conséquences: démonstration du théorème de l'indice, démonstration de la formule de l'indice et définition de la catégorie des coefficients $p$-adiques sur les courbes avec les propriétés de finitude requises.

Soient $p > 0$ un nombre premier, $\mathbb{Q}_p$ le complété du corps des nombres rationnels $\mathbb{Q}$ pour la valeur absolue $p$-adique, $\mathbb{Z}_p$ l'adhérence de $\mathbb{Z}$, $\mathbb{C}_p$ le complété d'une clôture algébrique de $\mathbb{Q}_p$ et $K$ une extension finie de $\mathbb{Q}_p$. Soit $\mathcal{R}_K(1)$ l'anneau des séries de Laurent à coefficients dans le corps $K$ qui convergent dans une couronne $|x| \in ]1 - \varepsilon, 1[$, pour $\varepsilon > 0$ non précisé. Appelons module différentiel sur $\mathcal{R}_K(1)$, un $\mathcal{R}_K(1)$-module libre de type fini muni d'une connexion. Soit $\mathcal{M}$ un module différentiel sur $\mathcal{R}_K(1)$ soluble en 1: c'est-à-dire dont la fonction rayon de convergence $R(\mathcal{M}, \rho)$ tend vers 1 avec $\rho$. Nous montrons le théorème:

THÉORÈME 6.1-14.    *Sous les conditions précédentes $\mathcal{M}$ est extension, dans la catégorie abélienne des modules différentiels sur l'anneau des fonctions analytiques au bord, de sa partie modérée $\mathcal{M}^{\leq 0}$ par sa partie de pentes strictement positives $\mathcal{M}_{>0}$:*

$$0 \to \mathcal{M}_{>0} \to \mathcal{M} \to \mathcal{M}^{\leq 0} \to 0.$$

La partie modérée $\mathcal{M}^{\leq 0}$ a la propriété de Robba au bord: son rayon de convergence est maximum pour tout $\rho \in ]1 - \varepsilon, 1[$ pour un $\varepsilon > 0$ assez petit alors que la partie de pentes strictement positives $\mathcal{M}_{>0}$ est injective au bord: ses solutions locales au point générique $t_\rho$ ont toutes un rayon de convergence strictement plus petit que $\rho$ pour tout $\rho \in ]1 - \varepsilon, 1[$ pour un $\varepsilon > 0$ assez petit. On définit les exposants de la monodromie locale $p$-adique de $\mathcal{M}$ comme ceux de sa partie modérée définis dans [C-M$_1$], [C-M$_2$]. L'obstruction à l'existence de l'indice se définit alors à l'aide des exposants de la monodromie locale $p$-adique de $\mathcal{M}$.

En fait nous montrons que $\mathcal{M}$ admet une filtration décroissante par des modules différentiels $\mathcal{M}_{>\gamma}$ solubles au bord, indexée par les nombres réels $\gamma \geq 0$, dont les sauts sont des nombres rationnels. Nous définissons ainsi le polygone de Newton $p$-adique de $\mathcal{M}$ et nous montrons que les sommets de ce polygone sont à coordonnées entières, c'est l'analogue $p$-adique du théorème de Hasse-Arf. Le théorème de décomposition de la partie modérée [C-M$_2$] joue alors le rôle $p$-adique du théorème de Grothendieck de la quasi-unipotence de la monodromie $\ell$-adique.

Le théorème de décomposition 6.1-14 est indépendant du théorème de décomposition de Dwork-Robba [D-R$_1$] sur le *corps* des éléments analytiques au bord. Il s'apparente plutôt au premier théorème de décomposition de Robba [R$_1$], mais en *famille*. Cependant il y a deux différences dans sa démonstration. D'une part, l'anneau structural n'étant plus un corps valué complet, la théorie des espaces normés est insuffisante et nous utilisons de façon essentielle la théorie des espaces vectoriels topologiques localement convexes sur un corps valué complet telle qu'elle est exposée dans le livre de A. Grothendieck [G$_1$]. D'autre part, les hypothèses du théorème des homomorphismes le plus général ([G$_1$, Chap. IV]) n'étant pas satisfaites, nous faisons appel à la théorie des équations différentielles $p$-adiques proprement dite à travers la structure de Frobenius de Christol-Dwork [C-D$_2$] pour les modules différentiels ayant un grand rayon de convergence, sans être maximum, sur une couronne de diamètre non nul.

Nous utilisons le théorème de l'indice pour montrer enfin une de nos principales motivations. Soient $X$ une variété algébrique affine non singulière sur le corps résiduel de l'anneau des entiers d'un sous corps à valuation discrète $K$ de $\mathbb{C}_p$ et $B_{i,p}(X)$ ses nombres de Betti $p$-adiques définis comme les dimensions de ses $K$-espaces de cohomologie de de Rham $p$-adiques définis par Monsky-Washnitzer [M-W]:

THÉORÈME 7.6-1. *Les nombres de Betti $B_{i,p}(X)$ sont finis pour tout $i$.*

En particulier les polynômes caractéristiques de l'endomorphisme de Frobenius opérant sur les espaces de cohomologies $p$-adiques $P_{i,p}(X)$ qui étaient définis comme des séries entières par P. Monsky [Mo$_2$] sont des éléments de $K[T]$ et peuvent se prêter aux comparaisons avec leurs analogues $l$-adiques. La réduction [Me$_4$] ramène le Théorème 7.6-1 à la finitude de la cohomologie $p$-adique des modules exponentiels $M_{\mathrm{Teich}(\bar{f}),n,m}$, étudiés dans [Me$_4$], qui est conséquence du théorème de l'indice.

La même méthode de réduction ([Me$_4$, 4.2.5]) ramène les propriétés de pureté des racines inverses des polynômes $P_{i,p}(X)$ à la pureté des valeurs propres de l'endomorphisme de Frobenius opérant sur les espaces de cohomologie $p$-adique des modules exponentiels $M_{\mathrm{Teich}(\bar{f}),n,m}$ avec la minoration requise du poids.

L'autre motivation, intimement liée à la précédente [M-N$_2$], est la construction de la catégorie des coefficients $p$-adiques au sens de Grothendieck ([G$_3$], [G$_4$], [G$_5$]) que nous réalisons dans [C-M$_4$] pour les courbes en utilisant tous nos résultats. La catégorie des coefficients $p$-adiques sur les courbes a toutes les propriétés de finitude et tous les invariants de la catégorie des coefficients $\ell$-adiques pour $\ell \neq p$ devrait permettre de démontrer par *voie $p$-adique* les propriétés de pureté des zéros et des pôles de la fonction zêta sur les corps finis ([Me$_4$, 4.2.5, 4.2.6]).

Les modules exponentiels $M_{\mathrm{Teich}(\bar{f}),n,m}$ nous ont servi de guide tout au long de ce travail. Ils fournissent des exemples d'équations provenant de la géométrie où le théorème de transfert pour les singularités irrégulières n'a pas lieu: la décomposition formelle de Turrittin est distincte en général de la décomposition $p$-adique du Théorème 6.1-14. Des contre exemples au principe de transfert pour les singularités irrégulières se trouvent aussi dans [Re]. Ils fournissent de même des exemples de modules rationnels solubles au point générique pour lesquels la décomposition du Théorème 6.1-14 n'a pas lieu sur le corps des éléments analytiques au bord ainsi que des exemples de modules de rang > 1 irréductibles et indécomposables par ramification.

Voici le contenu de ce travail. Nous définissons pour tout nombre réel $\gamma \geq 0$ dans le paragraphe 2 la topologie localement convexe métrisable $\mathcal{T}_\gamma$ sur l'anneau des opérateurs différentiels d'ordre fini à coefficients dans l'anneau des fonctions analytiques dans une couronne et la topologie localement convexe quotient $\mathcal{T}_{\gamma,Q}$ sur les modules de type fini, qui n'est pas en générale séparée. Nous montrons dans le paragraphe 3, que l'adhérence de zéro et le séparé associé d'un module différentiel pour cette topologie sont encore des modules différentiels, du moins quand le corps de base est maximalement complet. Nous montrons dans le paragraphe 4 que la fonction rayon de convergence d'un module différentiel soluble en $r$ est de la forme $\rho(\rho/r)^{pt(\mathcal{M})}$ pour un nombre rationnel $pt(\mathcal{M}) \geq 0$ et pour $\rho \in [r - \varepsilon, r[, \varepsilon > 0$ assez petit et nous définissons ainsi sa plus grande pente $p$-adique. Etant donné un module différentiel soluble muni d'une base, nous construisons, dans le paragraphe 5, une suite d'opérateurs différentiels dont les matrices dans cette base sont minorées et majorées. Nous utilisons dans le paragraphe 6 tous les résultats précédents pour montrer le théorème de décomposition selon la plus grande pente des modules différentiels solubles quand le corps de base est localement compact. Nous définissons la filtration $\mathcal{M}_{>\gamma}$, nous étudions ses propriétés fonctorielles et nous définissons le polygone de Newton Newton$(\mathcal{M}, p)$ $p$-adique. Nous montrons dans le paragraphe 7 que le théorème de décomposition des modules ayant la propriété de Robba [C-M$_2$] a déjà lieu sur le corps de base, le théorème de l'existence de l'indice local sous l'hypothèse (**NL**), pour les différences des exposants et pour les exposants. Nous en déduisons l'existence de l'indice local pour les opérateurs ayant une structure de Frobenius et le théorème de finitude des nombres de Betti $p$-adiques d'une variété algébrique affine non singulière sur un corps fini ou sur sa clôture algébrique. Nous définissons dans le paragraphe 8 l'indice généralisé local et nous étudions ses propriétés. Nous montrons la formule de l'indice local conjecturée par Robba [R$_5$] qui fait apparaître l'irrégularité $p$-adique, quand elle est finie, comme la hauteur du polygone de Newton $p$-adique. Nous en déduisons l'intégralité des sommets du polygone de Newton $p$-adique et la formule de l'indice globale pour les modules algébriques.

En utilisant toute la structure $p$-adique d'un point singulier d'une équation différentielle, nous montrons dans [C-M$_4$] le théorème d'existence de *réseaux* sous certaines conditions portant sur les exposants de la monodromie des $\mathcal{R}_K(1)$-modules différentiels et nous en déduisons le théorème d'algébrisation d'une classe de fibrés analytiques $p$-adiques qui nous permettent de construire la catégorie des coefficients $p$-adiques sur une courbe, contenant la catégorie des coefficients $p$-adiques modérément ramifiés ([C-M$_2$, 6.3]) et qui est stable par immersion *ouverte*. L'article [C-M$_4$] constituait à l'origine l'important paragraphe 9 de cet article sous le titre: La catégorie des coefficients $p$-adiques en dimension un. Suivant la recommandation du referee nous avons préféré le compléter sur quelques points de fondement pour le publier séparément.

Nous souhaitons à l'occasion de cet article rendre hommage à Bernard Dwork. La structure locale $p$-adique d'un point singulier d'une équation différentielle obtenue dans ce travail, qui est au coeur de la définition de la catégorie des coefficients $p$-adiques, nous voudrions insister auprès du lecteur sur ce point là, est l'aboutissement des idées très originales de B. Dwork. Outre l'introduction des méthodes $p$-adiques [D] dans la théorie de la cohomologie des variétés algébriques, il a mis en évidence les points clefs de la théorie des équations différentielles $p$-adiques, rayon de convergence au point générique d'un disque, structures de Frobenius, principe du transfert, propriété (**NL**) des exposants. Du reste il a conjecturé dès le Colloquium de Bombay de 1968 le théorème de l'indice: " Conjecture. *A linear differential operator in one variable with polynomial coefficients operating on functions holomorphic in an "open" disk has finite cokernel...* " ([D$_5$, §4, p. 88]). La démonstration du théorème de l'indice a été la principale motivation dans la théorie des équations différentielles $p$-adiques ces 25 dernières années.

Nous remercions une fois de plus Alberto Arabia pour l'aide considérable à tous les niveaux et l'encouragement qu'il nous a apportés tout au long de ce travail.

## Notations

Nous utilisons les notations des articles ([C-M$_1$], [C-M$_2$]) que nous avons reprises pour la plupart des articles de P. Robba. Soient un nombre premier $p > 0$ et $K$ un corps de caractéristique nulle muni d'une valeur absolue $p$-adique $|.|$ pour laquelle il est complet. Pour tout intervalle $I$ de nombres réels positifs on note $C(I)$ la couronne des nombres $x$ du complété de la clôture algébrique de $K$ tels que $|x| \in I$, $\mathcal{A}_K(I)$ l'espace vectoriel sur $K$ des séries de Laurent $\sum_{k \in \mathbb{Z}} a_k x^k$ à coefficients dans $K$ qui convergent pour tout $x$ dans la couronne $C(I)$. On note $\Delta^k := \frac{1}{k!} \frac{d^k}{dx^k}$ et $\pi$ une solution de l'équation $\pi^{p-1} + p = 0$.

On rappelle le résultat suivant: Pour tout intervalle $I$ la couronne $C(I)$ est une variété analytique rigide de Stein: pour tout faisceau analytique cohérent $\mathcal{F}$ on a le théorème $B$ de Cartan: l'espace de cohomologie cohérente $H^1(C(I), \mathcal{F})$ à valeurs dans le faisceau $\mathcal{F}$ est nul [Ki]. Pour un intervalle fermé c'est le théorème d'acyclicité de Tate.

Nous appelons une fois pour toute module différentiel sur un anneau différentiel un module libre de type fini sur cet anneau muni d'une connexion.

Pour simplifier les démonstrations dans cet article nous faisons, à partir du paragraphe 6, l'hypothèse de locale compacité sur le corps de base. Mais le lecteur pourra vérifier à partir des remarques 6.1-20 et 8.2-10 que tous les résultats de cet article s'étendent au cas d'un corps maximalement complet quelconque, en particulier au cas d'un corps à valuation discrète, ce qui est explicité dans [C-M$_4$].

## 2. Topologies localement convexes $\mathcal{T}_\gamma$ sur les anneaux d'opérateurs différentiels $\mathcal{A}_K(I)[\frac{d}{dx}]$, $\mathcal{R}_K(r)[\frac{d}{dx}]$

2.1. *Topologies localement convexe sur les anneaux* $\mathcal{A}_K(I), \mathcal{R}_K(r), \mathcal{H}_K^\dagger(r)$. Pour un nombre réel $r > 0$, on écrit $\mathcal{A}_K(r)$ pour $\mathcal{A}_K([0, r[$. Pour un nombre $\rho \in I$ et pour une fonction $f(x) = \sum_{k \in \mathbb{Z}} a_k x^k$ de $\mathcal{A}_K(I)$ on note

$$|f|_\rho := \sup_k |a_k| \rho^k$$

la norme de la convergence uniforme sur le cercle de rayon $\rho$. En vertu de la logarithme-convexité des normes $|-|_\rho$, la topologie sur $\mathcal{A}_K(I)$ définie par la famille des normes $|-|_\rho$ est celle de la convergence uniforme sur les sous-couronnes fermées. L'espace $\mathcal{A}_K(I)$ muni de la famille des normes $|-|_\rho$ est un espace de type $\mathcal{F}$ c'est-à-dire un espace localement convexe métrique complet. On note $\mathcal{R}_K(r)$ l'espace des fonctions analytiques au bord, réunion des espaces $\mathcal{A}_K([r-\varepsilon, r[$ pour un $\varepsilon > 0$. On définit l'espace $\mathcal{H}_K^\dagger(r)$ par la suite exacte:

$$0 \to \mathcal{A}_K(r) \to \mathcal{R}_K(r) \to \mathcal{H}_K^\dagger(r) \to 0.$$

Rappelons ([G$_1$, Chap. IV §1]) qu'on dit qu'une topologie localement convexe séparée sur un espace vectoriel est de type $\mathcal{LF}$ si c'est une limite inductive dénombrable de topologies de type $\mathcal{F}$. Si on note $H_K(\infty, r)$ l'espace des fonctions analytiques, à coefficients dans $K$, dans le disque fermé $D(\infty, r^+)$ nulles à l'infini, l'espace $\mathcal{H}_K^\dagger(r)$ est limite inductive des espaces $H_K(\infty, r - \varepsilon)$. De plus l'injection de $\mathcal{H}_K^\dagger(r)$ dans $H_K(\infty, r)$ est continue. La topologie limite inductive sur $\mathcal{H}_K^\dagger(r)$ est séparée, c'est donc une topologie de type $\mathcal{LF}$. C'est aussi un espace $\mathcal{DF}$ ([G$_1$, Chap. IV §3], [M-S]). L'espace $\mathcal{R}_K(r)$ muni de la topologie localement convexe limite inductive des espaces $\mathcal{A}_K([r-\varepsilon, r[$ est somme directe *topologique* de $\mathcal{A}_K(r)$ et de $\mathcal{H}_K^\dagger(r)$:

$$\mathcal{R}_K(r) = \mathcal{A}_K(r) \oplus \mathcal{H}_K^\dagger(r).$$

L'espace $\mathcal{R}_K(r)$ est alors un espace $\mathcal{LF}$.

2.2. *Topologies localement convexes* $\mathcal{T}_\gamma$ *sur les anneaux* $\mathcal{A}_K(I)[\frac{d}{dx}]$, $\mathcal{R}_K(r)[\frac{d}{dx}]$. Pour tout nombre réel $\gamma \geq 0$ nous allons munir les espaces des opérateurs différentiels $\mathcal{A}_K(I)[\frac{d}{dx}]$, $\mathcal{R}_K(r)[\frac{d}{dx}]$ d'une topologie localement convexe *séparée* $\mathcal{T}_\gamma$.

Si $P(x, \frac{d}{dx}) = \sum_k a_k(x)\Delta^k$ est un opérateur différentiel d'ordre fini à coefficients dans l'anneau $\mathcal{A}_K(I)$ on pose pour tout $\rho$ de l'intervalle $I$

$$|P|_{\gamma,\rho} := \sup_k |a_k|_\rho \, \rho^{-(\gamma+1)k}.$$

*Définition* 2.2-1. Pour un nombre réel $\gamma \geq 0$ on définit la topologie $\mathcal{T}_\gamma$ sur l'espace $\mathcal{A}_K(I)[\frac{d}{dx}]$ comme *la topologie définie par la famille des normes* $|-|_{\gamma,\rho} \, \rho \in I$ et la topologie $\mathcal{T}_\gamma$ sur l'espace $\mathcal{R}_K(r)[\frac{d}{dx}]$ comme *la topologie limite inductive des espaces métriques* $\mathcal{A}_K([r-\varepsilon,r[)[\frac{d}{dx}]$.

En particulier l'espace $\mathcal{A}_K(r)[\frac{d}{dx}]$ est muni d'une topologie d'espace métrique donc séparée. L'espace $\mathcal{H}_K^\dagger(r)[\frac{d}{dx}]$ est muni d'une topologie naturelle limite inductive d'espaces métriques $H_K(\infty, r-\varepsilon)[\frac{d}{dx}]$. L'injection de l'espace $\mathcal{H}_K^\dagger(r)[\frac{d}{dx}]$ dans l'espace $H(\infty,r)[\frac{d}{dx}]$ est continue. Comme ce dernier espace est normé donc séparé il en résulte que l'espace $\mathcal{H}_K^\dagger(r)[\frac{d}{dx}]$ est séparé pour la topologie $\mathcal{T}_\gamma$ pour tout réel $\gamma \geq 0$. Finalement l'isomorphisme topologique:

$$\mathcal{R}_K(r)[\frac{d}{dx}] = \mathcal{A}_K(r)[\frac{d}{dx}] \oplus \mathcal{H}_K^\dagger(r)[\frac{d}{dx}]$$

montre que la topologie $\mathcal{T}_\gamma$ sur l'espace $\mathcal{R}_K(r)[\frac{d}{dx}]$ est séparée pour tout réel $\gamma \geq 0$.

*Attention.* L'espace $\mathcal{A}_K(I)[\frac{d}{dx}]$ muni de la famille des normes $|-|_{\gamma,\rho}$ devient un espace métrique qui *n'est pas* complet, ni même limite inductive de complets, tout le problème est là. En effet si on note $\mathcal{A}_K(I)[\frac{d}{dx}]_m$ l'espace des opérateurs différentiels d'ordre au plus $m$ muni de la topologie induite par $\mathcal{T}_0$, la topologie limite inductive sur l'espace $\varinjlim \mathcal{A}_K(I)[\frac{d}{dx}]_m$ est *stricte* et donc ses parties bornées sont formées d'opérateurs différentiels d'ordre *borné* ([G$_1$, Chap. IV, Prop. 3]). Si la topologie $\mathcal{T}_0$ était $\mathcal{LF}$ elle serait équivalente en vertu du théorème des homomorphismes ([G$_1$, Chap. IV, Théorème 2]) à la topologie précédente. Mais les opérateurs $x^m \Delta^m, m \in \mathbb{N}$, forment une partie bornée pour la topologie $\mathcal{T}_0$, ce qui est contradictoire. Donc la topologie $\mathcal{T}_0$ sur $\mathcal{A}_K(I)[\frac{d}{dx}]$ n'est pas $\mathcal{LF}$. Un argument similaire vaut pout toute topologie $\mathcal{T}_\gamma$.

Contrairement au théorème de décomposition de Robba [R$_1$] où l'anneau de base des fonctions est un *corps* valué complet, les hypothèses du théorème des homomorphismes le plus général ([G$_1$, Chap. IV, Théorème 2]) ne sont pas satisfaites dans le cas des couronnes de diamètre non nul. Ceci montre les limites des méthodes de l'analyse fonctionnelle et nous aurons recours aux méthodes de la théorie des équations différentielles $p$-adiques proprement dite.

2.3. *Topologies quotients $\mathcal{T}_{\gamma,Q}$ sur les modules de type fini sur les anneaux* $\mathcal{A}_K(I)[\frac{d}{dx}]$, $\mathcal{R}_K(r)[\frac{d}{dx}]$.

*Définition 2.3-1.* Soit $\mathcal{M}$ un module à gauche de type fini sur l'anneau $\mathcal{A}_K(I)[\frac{d}{dx}]$ resp. $\mathcal{R}_K(r)[\frac{d}{dx}]$ on définit la topologie quotient $\mathcal{T}_{\gamma,Q}$ sur $\mathcal{M}$ comme *la topologie quotient* induite dans une présentation de $\mathcal{M}$

$$(\mathcal{A}_K(I)[\frac{d}{dx}])^m \to \mathcal{M} \to 0, \qquad \text{resp. } (\mathcal{R}_K(r)[\frac{d}{dx}])^m \to \mathcal{M} \to 0.$$

La topologie quotient $\mathcal{T}_{\gamma,Q}$ sur $\mathcal{M}$ ne dépend pas de la présentation. En effet si

$$(\mathcal{A}_K(I)[\frac{d}{dx}])^{m_1} \to \mathcal{M} \to 0$$

et

$$(\mathcal{A}_K(I)[\frac{d}{dx}])^{m_2} \to \mathcal{M} \to 0$$

sont deux présentations de $\mathcal{M}$ il existe un morphisme $\mathcal{A}_K(I)[\frac{d}{dx}]$-linéaire

$$(\mathcal{A}_K(I)[\frac{d}{dx}])^{m_1} \to (\mathcal{A}_K(I)[\frac{d}{dx}])^{m_2}$$

qui relève l'application identique de $\mathcal{M}$. Ce morphisme est représenté par une matrice $(m_1, m_2)$ à coefficients dans l'anneau $\mathcal{A}_K(I)[\frac{d}{dx}]$ et donc est *continu* pour la topologie $\mathcal{T}_\gamma$. Le même argument vaut pour l'anneau $\mathcal{R}_K(r)[\frac{d}{dx}]$.

Pour les mêmes raisons tout morphisme $\mathcal{A}_K(I)[\frac{d}{dx}]$-linéaire: $\mathcal{M}_1 \to \mathcal{M}_2$ entre $\mathcal{A}_K(I)[\frac{d}{dx}]$-modules de type *fini* est automatiquement continu.

La topologie $\mathcal{T}_{\gamma,Q}$ sur un module de type fini $\mathcal{M}$ sur $\mathcal{A}_K(I)[\frac{d}{dx}]$ *n'est en général pas séparée.* Toute notre idée de départ dans cet article est d'analyser la structure de l'adhérence de zéro et du séparé associé dans une situation particulière.

## 3. Propriétés de finitude algébro-topologiques des modules libres à connexion sur l'anneau $\mathcal{A}_K(I)$

3.1. *Propriétés de finitude algébriques.* Les anneaux $\mathcal{A}_K(I)$ et $\mathcal{R}_K(r)$ ne sont pas principaux ni même noethériens. Cependant comme conséquence du théorème de M. Lazard [L] on obtient:

PROPOSITION 3.1-1. *Soit $I$ un intervalle; si la valuation de $K$ est discrète, ou plus généralement si $K$ est maximalement complet, tout idéal de type fini de l'anneau $\mathcal{A}_K(I)$ resp. de l'anneau $\mathcal{R}_K(r)$ est principal.*

*Démonstration.* Soit $\mathcal{I}$ un idéal de type fini de l'anneau $\mathcal{A}_K(I)$. Par récurrence sur le nombre de générateurs on peut supposer qu'il est engendré par deux fonctions $(g_1, g_2)$. Soit $\tilde{h}$ le plus grand diviseur au sens de [L] commun des diviseurs des zéros de $g_1$ et de $g_2$ dans la couronne $C(I)$. En vertu du théorème de M. Lazard [L], si $K$ est maximalement complet, en particulier s'il est de valuation discrète, le diviseur $\tilde{h}$ est celui d'une fonction analytique $h$ de $\mathcal{A}_K(I)$. Nous allons voir que l'idéal $(g_1, g_2)$ est égal à l'idéal $(h)$ dans l'anneau $\mathcal{A}_K(I)$. Si $g_1 = h_1 h$ et $g_2 = h_2 h$ dans $\mathcal{A}_K(I)$, pour tout intervalle fermé $I' \subset I$, l'idéal de l'anneau principal $\mathcal{A}_K(I')$ engendré par $(h_1, h_2)$ contient 1. Par limite projective, l'espace des sections globales sur la couronne $C(I)$ du faisceau d'idéaux engendré par $(h_1, h_2)$ est égal à $\mathcal{A}_K(I)$. Mais comme la couronne $C(I)$ est une variété de Stein en vertu du théorème B de Cartan ceci entraîne que l'idéal $(h_1, h_2)$ est égal à $\mathcal{A}_K(I)$. Cela montre que les idéaux $(h)$ et $(g_1, g_2)$ de l'anneau $\mathcal{A}_K(I)$ sont égaux. Le cas des idéaux de type fini de l'anneau $\mathcal{R}_K(r)$ se ramène au cas des idéaux de type fini d'un anneau $\mathcal{A}_K([r - \varepsilon, r[)$.

En particulier sous les conditions de la proposition 3.1-1 les anneaux $\mathcal{A}_K(I)$ et $\mathcal{R}_K(r)$ sont cohérents: les idéaux de type fini sont de *présentation finie*.

La proposition précédente montre que si le corps $K$ est maximalement complet la dimension homologique plate des anneaux $\mathcal{A}_K(I)$ et $\mathcal{R}_K(r)$ est égale à un et donc tout module sans torsion sur ces anneaux est plat.

COROLLAIRE 3.1-2. *Si le corps de base $K$ est maximalement complet, tout sous-module de type fini d'un module libre de type fini sur les anneaux $\mathcal{A}_K(I)$ et $\mathcal{R}_K(r)$ est libre de type fini.*

*Démonstration.* Soit $\mathcal{M}$ un sous $\mathcal{A}_K(I)$-module de type fini d'un module libre $(\mathcal{A}_K(I))^m$. On raisonne par récurrence sur $m$. Si $m$ est égal à un c'est la proposition 3.1-1. Si $m > 1$ la projection de $\mathcal{M}$ sur $(\mathcal{A}_K(I))^{m-1}$ est libre de type fini par hypothèse de récurrence et son noyau qui est de type fini est un idéal de $\mathcal{A}_K(I)$ et donc libre. Donc $\mathcal{M}$ est libre de type fini. Le même raisonnement vaut pour l'anneau $\mathcal{R}_K(r)$, en fait pour tout anneau où tout idéal de type fini est libre.

3.2. *Propriétés de finitude topologiques.* Nous démontrons dans ce paragraphe le théorème de finitude de nature algébro-topologique qui joue un rôle essentiel dans le théorème de décomposition 6.1-14:

THÉORÈME 3.2-1. *Soit $I$ un intervalle et $K$ un corps maximalement complet, par exemple à valuation discrète. Soit $\mathcal{M}$ un $\mathcal{A}_K(I)$-module libre de rang $m$ à connexion. Alors l'adhérence $\bar{0}_\gamma(\mathcal{M})$ de zéro dans $\mathcal{M}$ pour la topologie quotient $\mathcal{T}_{\gamma,Q}$ et le séparé associé $\mathcal{M}/\bar{0}_\gamma(\mathcal{M})$ sont des $\mathcal{A}_K(I)$-modules libres à connexion de rang fini pour tout nombre réel $\gamma \geq 0$.*

*Démonstration.* Soit $G(x)$ la matrice de la connexion dans une base de $\mathcal{M}$, c'est donc une matrice carrée de rang $m$ à coefficients dans l'anneau $\mathcal{A}_K(I)$. On a alors la présentation

$$0 \to (\mathcal{A}_K(I)[\frac{d}{dx}])^m \xrightarrow{u := \frac{d}{dx} - G(x)} (\mathcal{A}_K(I)[\frac{d}{dx}])^m \to \mathcal{M} \to 0$$

qui est en fait une suite exacte. On note $\mathcal{I}$ l'image du morphisme $u$; son adhérence $\bar{\mathcal{I}}$ dans $(\mathcal{A}_K(I)[\frac{d}{dx}])^m$ pour la topologie $\mathcal{T}_\gamma$ est un sous $\mathcal{A}_K(I)[\frac{d}{dx}]$-module et l'adhérence de zéro $\bar{0}_\gamma(\mathcal{M})$ s'identifie au quotient $\bar{\mathcal{I}}/\mathcal{I}$. La suite exacte

$$0 \to \bar{0}_\gamma(\mathcal{M}) \to \mathcal{M} \to \mathcal{M}/\bar{0}_\gamma(\mathcal{M}) \to 0$$

s'identifie à la suite exacte de $\mathcal{A}_K(I)[\frac{d}{dx}]$-modules

$$0 \to \bar{\mathcal{I}}/\mathcal{I} \to (\mathcal{A}_K(I)[\frac{d}{dx}])^m/\mathcal{I} \to (\mathcal{A}_K(I)[\frac{d}{dx}])^m/\bar{\mathcal{I}} \to 0.$$

Soit $I_\varepsilon$ une famille décroissante d'intervalles fermés contenus dans $I$ tels que $I = \bigcup_{\varepsilon>0} I_\varepsilon$. Pour tout $\varepsilon > 0$ notons $\mathcal{M}_\varepsilon$ la restriction de $\mathcal{M}$ à la couronne fermée $C(I_\varepsilon)$:

$$\mathcal{M}_\varepsilon := \mathcal{A}_K(I_\varepsilon) \otimes_{\mathcal{A}_K(I)} \mathcal{M}.$$

C'est un module libre de rang $m$ sur l'anneau principal $\mathcal{A}_K(I_\varepsilon)$ muni d'une connexion. On peut considérer la topologie quotient $\mathcal{T}_{\gamma,Q}$ sur $\mathcal{M}_\varepsilon$ et l'adhérence de zéro $\bar{0}_\gamma(\mathcal{M}_\varepsilon)$. C'est un sous-module à connexion de $\mathcal{M}_\varepsilon$ et donc libre de rang $m_\varepsilon \leq m$ sur l'anneau principal $\mathcal{A}_K(I_\varepsilon)$; c'est là un résultat général (cf. ([C$_1$, 4.3])). D'autre part pour $\varepsilon$ variable les morphismes de restriction

$$\mathcal{M}_{\varepsilon'} \to \mathcal{M}_\varepsilon, \qquad 0 < \varepsilon' < \varepsilon$$

sont évidemment continus. On a donc des morphismes continus de transition pour les topologies induites

$$\bar{0}_\gamma(\mathcal{M}_{\varepsilon'}) \overset{t_{\varepsilon',\varepsilon}}{\to} \bar{0}_\gamma(\mathcal{M}_\varepsilon).$$

En vertu de la platitude, le module engendré par l'image de $t_{\varepsilon',\varepsilon}$ est égale à

$$\mathcal{A}_K(I_\varepsilon) \otimes_{\mathcal{A}_K(I_{\varepsilon'})} \bar{0}_\gamma(\mathcal{M}_{\varepsilon'}).$$

C'est donc un module libre de rang $m_{\varepsilon'}$ à connexion sur l'anneau $\mathcal{A}_K(I_\varepsilon)$. Cela montre que le rang $m_\varepsilon$ est décroissant avec $\varepsilon$ et donc est stationnaire pour $\varepsilon$ assez petit, disons pour $\varepsilon \leq \varepsilon_0$, de valeur $m_{\varepsilon_0}$.

Considérons le recouvrement de la couronne ouverte $C(I)$ par les couronnes fermées $C(I_\varepsilon)$ pour $0 < \varepsilon \leq \varepsilon_0$. Sur chaque couronne fermée $C(I_\varepsilon)$ soit $\widetilde{\bar{0}_\gamma}(\mathcal{M}_\varepsilon)$ le faisceau associé au module $\bar{0}_\gamma(\mathcal{M}_\varepsilon)$. Puisque le rang $m_\varepsilon$ est constant on a alors les isomorphismes pour $0 < \varepsilon' < \varepsilon \leq \varepsilon_0$:

$$\mathcal{A}_K(I_\varepsilon) \otimes_{\mathcal{A}_K(I_{\varepsilon'})} \bar{0}_\gamma(\mathcal{M}_{\varepsilon'}) \simeq \bar{0}_\gamma(\mathcal{M}_\varepsilon).$$

Cela exprime les conditions de compatibilité entre les faisceaux $\widetilde{\bar{0}_\gamma}(\mathcal{M}_\varepsilon)$. Ces faisceaux se recollent pour donner naissance à un faisceau $\widetilde{\bar{0}_\gamma}(\mathcal{M})$ sur la couronne $C(I)$. C'est un faisceau localement libre sur le faisceau structural $\widetilde{\mathcal{A}}_K(I)$ de rang $m_{\varepsilon_0}$ par construction.

Considérons la suite exacte de faisceaux analytiques cohérents sur la couronne $C(I)$

$$(*) \qquad\qquad 0 \to \widetilde{\bar{0}_\gamma(\mathcal{M})} \to \widetilde{\mathcal{M}} \to \widetilde{\mathcal{M}/\bar{0}_\gamma(\mathcal{M})} \to 0.$$

En vertu du théorème d'acyclicité de Tate l'espace des sections globales sur une couronne fermée $C(I_\varepsilon)$ du faisceau quotient de la suite $(*)$ est isomorphe au module quotient $\mathcal{M}_\varepsilon/\bar{0}_\gamma(\mathcal{M}_\varepsilon)$ *à connexion*. Il est donc libre sur l'anneau $\mathcal{A}_K(I_\varepsilon)$ de rang $m - m_{\varepsilon_0}$. La suite exacte $(*)$ est donc une suite exacte de faisceaux analytiques *localement libres de rang fini*. Le faisceau $\widetilde{\mathcal{M}}$ apparaît comme une extension de $\widetilde{\mathcal{M}/\bar{0}_\gamma(\mathcal{M})}$ par $\widetilde{\bar{0}_\gamma(\mathcal{M})}$ dans la catégorie abélienne des faisceaux analytiques sur la couronne $C(I)$. Il correspond, en vertu du théorème de Yoneda, à un élément du groupe

$$\mathrm{Ext}^1_{\mathcal{A}_K(I)}(C(I); \widetilde{\mathcal{M}/\bar{0}_\gamma(\mathcal{M})}, \widetilde{\bar{0}_\gamma(\mathcal{M})}).$$

Ce groupe est isomorphe, puisque le faisceau quotient $\widetilde{\mathcal{M}/\bar{0}_\gamma(\mathcal{M})}$ est *localement libre*, au premier groupe de cohomologie cohérente

$$H^1\Big(C(I); \mathcal{H}om_{\mathcal{A}_K(I)}(\widetilde{\mathcal{M}/\bar{0}_\gamma(\mathcal{M})}, \widetilde{\bar{0}_\gamma(\mathcal{M})})\Big)$$

qui est nul en vertu du théorème B de Cartan. Le faisceau $\widetilde{\mathcal{M}}$ est donc somme directe:

$$\widetilde{\mathcal{M}} = \widetilde{\bar{0}_\gamma(\mathcal{M})} \oplus \widetilde{\mathcal{M}/\bar{0}_\gamma(\mathcal{M})}$$

et le module de ses sections globales est somme directe:

$$\Gamma(C(I); \widetilde{\mathcal{M}}) = \Gamma(C(I); \widetilde{\bar{0}_\gamma(\mathcal{M})}) \oplus \Gamma(C(I); \widetilde{\mathcal{M}/\bar{0}_\gamma(\mathcal{M})}).$$

Le module de type fini $\Gamma(C(I); \widetilde{\mathcal{M}/\bar{0}_\gamma(\mathcal{M})})$ apparaît comme un sous-module d'un module libre de rang fini. Il est donc libre de rang fini sur l'anneau $\mathcal{A}_K(I)$ en vertu du corollaire 3.1-2. Le sous-module $\Gamma(C(I); \widetilde{\bar{0}_\gamma(\mathcal{M})})$ apparaît comme un module quotient d'un module libre de rang fini. Il est de type fini et donc libre en vertu du corollaire 3.1-2.

Le morphisme de restriction

$$\mathcal{M} \to \mathcal{M}_\varepsilon$$

est continu pour la topologie quotient $\mathcal{T}_{\gamma,Q}$ pour tout $\varepsilon$. On a donc un morphisme pour tout $\varepsilon$:

$$\bar{0}_\gamma(\mathcal{M}) \to \bar{0}_\gamma(\mathcal{M}_\varepsilon)$$

et donc un morphisme

$$(**) \qquad\qquad \bar{0}_\gamma(\mathcal{M}) \to \varprojlim_{\varepsilon \le \varepsilon_0} \bar{0}_\gamma(\mathcal{M}_\varepsilon).$$

LEMME 3.2-2.  *Le morphisme précédent* $(**)$ *est un isomorphisme.*

*Démonstration.* Remarquons que les deux membres de (∗∗) sont des sous-modules de $\mathcal{M} = \varprojlim_{\varepsilon \leq \varepsilon_0} \mathcal{M}_\varepsilon$. Par définition de la topologie quotient $\mathcal{T}_{\gamma, Q}$ sur $\mathcal{M}$, $\bar{0}_\gamma(\mathcal{M})$ est l'intersection de toutes les boules, de l'espace vectoriel topologique localement convexe $\mathcal{M}$, centrées en zéro pour les semi-normes quotients induites par les semi-normes $| - |_{\gamma, \rho}, \rho \in I$. Cette intersection est égale à l'intersection de $\bar{0}_\gamma(\mathcal{M}_\varepsilon) \cap \mathcal{M}$ pour $0 < \varepsilon \leq \varepsilon_0$. On a donc une bijection

$$\bar{0}_\gamma(\mathcal{M}) \simeq \varprojlim_{\varepsilon \leq \varepsilon_0} \bar{0}_\gamma(\mathcal{M}_\varepsilon) \cap \mathcal{M}.$$

D'autre part l'injection

$$\varprojlim_{\varepsilon \leq \varepsilon_0} \bar{0}_\gamma(\mathcal{M}_\varepsilon) \cap \mathcal{M} \to \varprojlim_{\varepsilon \leq \varepsilon_0} \bar{0}_\gamma(\mathcal{M}_\varepsilon)$$

est surjective donc bijective. D'où le lemme.

Remarquons que le membre de droite de (∗∗) n'est rien d'autre que le module des sections globales du faisceau $\widetilde{\bar{0}_\gamma(\mathcal{M})}$ dont nous venons de voir que c'est un module libre de rang fini sur l'anneau $\mathcal{A}_K(I)$. Cela entraîne que le séparé associé du module $\mathcal{M}$ est isomorphe au module des sections globales du faisceau quotient $\widetilde{\mathcal{M}}/\widetilde{\bar{0}_\gamma(\mathcal{M})}$ et donc est libre de rang fini. En fait, on trouve que le rang de $\bar{0}_\gamma(\mathcal{M})$ est égale à $m_{\varepsilon_0}$ alors que celui de $\mathcal{M}/\bar{0}_\gamma(\mathcal{M})$ est égal à $m - m_{\varepsilon_0}$. D'où le théorème 3.2-1.

Si $\mathcal{M}$ est un module à connexion libre de rang fini sur l'anneau $\mathcal{A}_K(I)$. On construit la filtration de $\mathcal{M}$ en posant $\mathcal{M}_0 := \mathcal{M}$ et en définissant $\mathcal{M}_i, i \geq 1$, comme l'adhérence de zéro dans $\mathcal{M}_{i-1}$ pour la topologie quotient $\mathcal{T}_{\gamma, Q}$. C'est un module à connexion libre de type fini sur l'anneau $\mathcal{A}_K(I)$ en vertu du théorème 3.2-1.

COROLLAIRE 3.2-3. *La filtration décroissante $\mathcal{M}_i, i \geq 0$, est stationnaire.*

*Démonstration.* Le rang $m_i$ de $\mathcal{M}_i$ est une fonction décroissante avec $i$ donc stationnaire. Comme les modules $\mathcal{M}_i$ sont à connexion, cela montre que la filtration elle-même $\mathcal{M}_i$ est stationnaire.

*Exemple* 3.2-4. Prenons $\gamma = 0$ et $\mathcal{M}$ défini par l'opérateur $P = (x\frac{d}{dx} - \alpha)^m$ où $\alpha$ est un nombre de $\mathbb{Z}_p$ et $m \geq 1$, on trouve que $\mathcal{M}_i$ est défini par l'opérateur $P_i = (x\frac{d}{dx} - \alpha)^{m-i}$ pour $i = 0, \ldots, m-1$ et $\mathcal{M}_i = 0$ pour $i \geq m$.

Pour un intervalle $I$ resp. un nombre réel $r > 0$ nous notons

$$\text{MLC}(\mathcal{A}_K(I)), \text{ resp. } \text{MLC}(\mathcal{R}_K(r)),$$

la catégorie des $\mathcal{A}_K(I)$-modules, resp. $\mathcal{R}_K(r)$-modules, libres de rang fini à connexion.

THÉORÈME 3.2-5. *Pour tout intervalle $I$ et tout nombre réel $r > 0$ si le corps de base $K$ est maximalement complet les catégories $\text{MLC}(\mathcal{A}_K(I))$ et $\text{MLC}(\mathcal{R}_K(r))$ sont abéliennes.*

*Démonstration.* Soit $v : \mathcal{M}_1 \to \mathcal{M}_2$ un morphisme dans la catégorie $\mathrm{MLC}(\mathcal{A}_K(I))$. Le noyau de $v$ est un $\mathcal{A}_K(I)$-module de type fini, parce que l'anneau $\mathcal{A}_K(I)$ est cohérent, donc libre en vertu du corollaire 3.1-2. En vertu du théorème B de Cartan le conoyau de $v$ est isomorphe à l'espace des sections globales du conoyau du morphisme de faisceaux libres associé à $v$. Le raisonnement de la démonstration du Théorème 3.2-1 montre alors qu'il est isomorphe à un facteur direct de type fini d'un module libre de rang fini. Il est donc libre de rang fini. Le cas de la catégorie $\mathrm{MLC}(\mathcal{R}_K(r))$ se ramène au cas de la catégorie $\mathrm{MLC}(\mathcal{A}_K([r-\varepsilon, r[))$. D'où le théorème 3.2-5.

*Remarque* 3.2-6. Il est facile de trouver des exemples de module quotient (de type fini) à connexion d'un module libre à connexion de rang fini sur l'anneau $\mathcal{A}_K(I)$ *qui ne sont pas libres* contrairement à la situation sur l'anneau $\mathcal{A}_K(I_\varepsilon)$ pour un intervalle fermé $I_\varepsilon$. Cela montre que le théorème 3.2-1 ne peut être atteint par voie purement algébrique en raison de l'existence d'idéaux non de type fini maximaux de l'anneau $\mathcal{A}_K(I)$ formés de fonctions qui n'ont pas de zéro commun. Cela illustre l'importance des méthodes cohomologiques dans ce genre de questions.

*Remarque* 3.2-7. Si on part d'un module libre $\mathcal{M}$ à connexion sur l'anneau $\mathcal{R}_K(r)$ il provient d'un module $\mathcal{M}_{[r-\varepsilon, r[}$ libre à connexion sur l'anneau $\mathcal{A}_K([r-\varepsilon, r[)$ pour $\varepsilon$ assez petit. On a alors une application continue injective

$$\lim_{\substack{\to \\ 0 < \varepsilon}} \bar{0}_\gamma(\mathcal{M}_{[r-\varepsilon, r[}) \to \bar{0}_\gamma(\mathcal{M})$$

qui n'a aucune raison d'être a priori bijective parce qu'une limite inductive d'espaces topologiques séparés n'est pas en général séparée. On ne peut pas conclure à l'analogue du théorème 3.2-1 pour l'anneau $\mathcal{R}_K(r)$ à partir du théorème 3.2-1 pour les anneaux $\mathcal{A}_K([r-\varepsilon, r[)$. Il faut un autre argument.

Nous allons définir des sous-catégories de la catégorie abélienne $\mathrm{MLC}(\mathcal{R}_K(r))$ formées de modules différentiels dont la connexion a des propriétés de plus en plus restrictives pour arriver à la catégorie de base $\mathrm{MLS}(\mathcal{R}_K(r), \mathbf{NL}^{**})$ dans le paragraphe 7 qui nous permettra en particulier de démontrer le théorème de l'indice.

## 4. La plus grande pente $p$-adique d'un $\mathcal{R}_K(r)$-module libre de rang fini à connexion soluble en $r$

4.1. *La fonction rayon de convergence.* Soit $f$ une fonction définie sur une partie de $\mathbb{R}^+$ à valeurs dans $\mathbb{R}^+$. Rappelons (cf. [R-C]) que l'on dit que la fonction $f$ a *logarithmiquement* une propriété (continuité, dérivabilité, convexité,...) si la fonction $\tilde{f} = \mathrm{Log} \circ f \circ \exp$ a cette propriété.

Soit $\mathcal{M}$ un $\mathcal{R}_K(r)$-module libre de rang $m$ à connexion. Sa fonction rayon de convergence $\rho \mapsto R(\mathcal{M}, \rho)$ est définie sur un intervalle $[r - \varepsilon, r[$ par:

$$R(\mathcal{M}, \rho) := \min(\rho, \liminf_{k \to \infty}(|G_k|_\rho)^{-1/k})$$

où $G_k$ est la matrice de l'opérateur $\Delta^k$ dans une base. La fonction rayon de convergence ne dépend pas de la base choisie: c'est le rayon de convergence dans le disque générique $D(t_\rho, \rho^-)$ (cf. [C-D$_2$]). C'est donc une fonction à valeurs réelles positives, logarithmiquement concave (cf. [R-C]). La limite $R(\mathcal{M}, r^-) := \lim_{\rho \to r^-} R(\mathcal{M}, \rho)$ existe et, de plus, la fonction $R(\mathcal{M}, \rho)$ est logarithmiquement dérivable à gauche, de valeur éventuellement infini, quand $\rho \to r^-$. Par construction $R(\mathcal{M}, r^-) \le r$. Posons $\omega := p^{-\frac{1}{p-1}}$ et $\omega' := p^{-\frac{1}{p}}$.

*Définition* 4.1-1. Nous dirons que $\mathcal{M}$ est *soluble* en $r$ si $R(\mathcal{M}, r^-) = r$.

On note la catégorie des $\mathcal{R}_K(r)$-modules libres de type fini à connexion solubles en $r$ par

$$\mathrm{MLS}(\mathcal{R}_K(r)).$$

PROPOSITION 4.1-2. *Si le corps de base $K$ est maximalement complet la catégorie $\mathrm{MLS}(\mathcal{R}_K(r))$ est abélienne.*

*Démonstration.* Pour tout intervalle $I$ et tout corps complet $K$, si

$$0 \to \mathcal{M}_1 \to \mathcal{M} \to \mathcal{M}_2 \to 0$$

est une suite exacte de $\mathcal{A}_K(I)$-modules libres de rang fini à connexion on a:

$$R(\mathcal{M}, \rho) = \min(R(\mathcal{M}_1, \rho), R(\mathcal{M}_2, \rho))$$

pour tout $\rho \in I$ parce que le foncteur qui à un module différentiel associe l'espace de ses solutions dans un disque centré en un point générique $t_\rho$ et de rayon borné par $\rho$ est *exact* ([R$_1$], 4.23). En particulier si $\mathcal{M}$ est soluble en $r$, $\mathcal{M}_1$ et $\mathcal{M}_2$ sont solubles en $r$. La proposition est alors conséquence du théorème 3.2-5.

Soit $\varphi : C(]r, R[) \to C(]r^p, R^p[)$ la ramification de Frobenius $x \to x^p$; on dit qu'un module différentiel $\mathcal{N}$ sur la couronne $C(]r^p, R^p[)$ est un antécédent d'un module différentiel $\mathcal{M}$ sur la couronne $C(]r, R[)$ s'il existe un isomorphisme entre $\mathcal{M}$ et l'image inverse $\varphi^* \mathcal{N}$ de $\mathcal{N}$. Rappelons le théorème de structure de Frobenius de Christol-Dwork ([C-D$_2$], 5.4, 4.4]) pour les modules différentiels ayant un grand rayon de convergence sur une couronne de diamètre non nul:

THÉORÈME 4.1-3. 1) *Soit $\mathcal{M}$ un $\mathcal{A}_K(]r, R[)$-module différentiel sur la couronne $C(]r, R[)$ tel que $R(\mathcal{M}, \rho) > \omega' \rho$ pour tout $\rho$ dans l'intervalle $]r, R[$; alors il existe un $\mathcal{A}_K(]r^p, R^p[)$-module différentiel $\mathcal{N}$ sur la couronne $C(]r^p, R^p[)$* *tel que $(R(\mathcal{M}, \rho))^p = R(\mathcal{N}, \rho^p)$ pour tout $\rho$ dans l'intervalle $]r, R[$, qui soit un antécédent de $\mathcal{M}$ et qui est unique à isomorphisme près.*

2) *Soit $\mathcal{M}$ un $\mathcal{A}_K(]r, R[)$-module différentiel sur la couronne $C(]r, R[)$ tel que $R(\mathcal{M}, \rho) > \omega \rho$ pour tout $\rho$ dans l'intervalle $]r, R[$; alors il existe un $\mathcal{M}_K(]r^p, R^p[)$-module différentiel $\mathcal{N}$ sur la couronne $C(]r^p, R^p[)$, où $\mathcal{M}_K(]r, R[)$ est le corps des fractions de $\mathcal{A}_K(]r, R[)$, tel que $(R(\mathcal{M}, \rho))^p = R(\mathcal{N}, \rho^p)$ pour tout $\rho$ dans l'intervalle $]r, R[$, qui soit un antécédent de $\mathcal{M}$ et qui est unique à isomorphisme près et qui n'admet que des singularités apparentes sur la couronne $C(]r^p, R^p[)$.*

La démonstration du théorème précédent est faite dans [C-D$_2$] pour les éléments analytiques dans une couronne ouverte. Mais le lecteur pourra vérifier qu'elle vaut aussi pour les fonctions analytiques dans une couronne ouverte. D'autre part nous utilisons dans cet article que l'existence de l'antécédent et non son unicité. En itérant on trouve:

THÉORÈME 4.1-4. *Soit $\mathcal{M}$ un $\mathcal{A}_K(]r, R[)$-module différentiel sur la couronne $C(]r, R[)$ tel que $R_\rho(\mathcal{M}) > \rho \omega^{1/p^{h-1}}$ pour tout $\rho$ dans l'intervalle $]r, R[$ pour un entier $h \geq 1$, alors il existe un $\mathcal{M}_K(]r^{p^h}, R^{p^h}[)$-module différentiel $\mathcal{N}_h$ sur la couronne $C(]r^{p^h}, R^{p^h}[)$ tel que $(R(\mathcal{M}, \rho))^{p^h} = R(\mathcal{N}, \rho^{p^h})$ pour tout $\rho$ dans l'intervalle $]r, R[$, qui soit un antécédent d'ordre $h$ de $\mathcal{M}$, qui est unique à isomorphisme près et qui n'admet que des singularité apparentes sur la couronne $C(]r^{p^h}, R^{p^h}[)$.*

*Démonstration.* Si $h = 1$ c'est la deuxième partie du théorème 4.1-3. On peut supposer $h \geq 2$. Comme $\omega'^{1/p^{h-2}} \leq \omega^{1/p^{h-1}}$ la condition 1) du théorème 4.1-3 est réalisée $h - 1$-fois et en itérant on obtient un $\mathcal{A}_K(]r^{p^{h-1}}, R^{p^{h-1}}[)$-module différentiel $\mathcal{N}_{h-1}$ sur la couronne $C(]r^{p^{h-1}}, R^{p^{h-1}}[)$ tel que

$$(R(\mathcal{M}, \rho))^{p^{h-1}} = R(\mathcal{N}_{h-1}, \rho^{p^{h-1}})$$

pour tout $\rho$ dans l'intervalle $]r, R[$, qui soit un antécédent d'ordre $h-1$ de $\mathcal{M}$ et qui est unique à isomorphisme près. Comme $\omega \rho^{p^{h-1}} < R(\mathcal{N}_{h-1}, \rho^{p^{h-1}})$, on est sous les conditions d'application de la deuxième partie du théorème 4.1-3 et l'on obtient l'antécédent $\mathcal{N}_h$ qui a les propriétés du théorème 4.1-4.

4.2. *La plus grande pente $p$-adique.*

THÉORÈME 4.2-1. *Soit $\mathcal{M}$ un $\mathcal{R}_K(r)$-module libre de rang $m \geq 1$ à connexion soluble en $r$, il existe un nombre rationnel $\beta := \beta(\mathcal{M}) \geq 0$ tel que la fonction $R(\mathcal{M}, \rho)$ est égale à $\rho(\rho/r)^\beta$ pour $\rho$ dans l'intervalle $[r - \varepsilon, r[$ pour $\varepsilon > 0$ assez petit.*

*Démonstration.* Nous supposons que $\mathcal{M}$ est défini dans un intervalle $[r - \varepsilon, r[$ pour un $\varepsilon > 0$. Nous pouvons supposer que $\mathcal{M}$ n'a pas la propriété de Robba dans la couronne $C([r - \varepsilon, r[)$, c'est-à-dire la fonction $R(\mathcal{M}, \rho) < \rho$ n'est pas égale à $\rho$ [C-M$_2$]. Par logarithme-concavité nous pouvons supposer que $R(\mathcal{M}, \rho) < \rho$ pour $\rho$ dans l'intervalle $[r - \varepsilon, r[$ pour un $\varepsilon > 0$ assez petit. Considérons la suite de points $r_h$ croissante tendant vers $r$ définis pour $h$ assez grand de l'intervalle $[r - \varepsilon, r[$ tels que

$$\rho\omega^{1/p^{h-1}} < R(\mathcal{M}, \rho) < \rho\omega^{1/p^h}$$

pour $\rho \in ]r_{h-1}, r_h[$. Dans l'intervalle ouvert $]r_{h-1}, r_h[$, il existe, en vertu du théorème 4.1-4, un antécédent de Frobenius de Christol-Dwork $\mathcal{N}_h$ d'ordre $h$ de $\mathcal{M}$ qui est un $\mathcal{M}_K(]r_{h-1}^{p^h}, r_h^{p^h}[)$-module différentiel sur la couronne $C(]r_{h-1}^{p^h}, r_h^{p^h}[)$ tel que $(R(\mathcal{M}, \rho))^{p^h} = R(\mathcal{N}_h, \rho^{p^h})$ pour tout $\rho$ dans l'intervalle $]r_{h-1}, r_h[$. On a donc $R(\mathcal{N}_h, \rho) < \rho\omega$ pour $\rho$ dans l'intervalle $]r_{h-1}^{p^h}, r_h^{p^h}[$. On est dans les conditions d'application du théorème de Robba [R$_6$] et de Young [Y$_1$], [Y$_2$]. En vertu du lemme du vecteur cyclique, $\mathcal{N}_h$ est isomorphe à

$$\mathcal{M}_K(]r_{h-1}^{p^h}, r_h^{p^h}[)[x\frac{d}{dx}]/\mathcal{M}_K(]r_{h-1}^{p^h}, r_h^{p^h}[)((x\frac{d}{dx})^m + a_{m-1}(x)(x\frac{d}{dx})^{m-1} + \cdots + a_0(x)).$$

La fonction rayon de convergence de $\mathcal{N}_h$ est alors donnée explicitement ([C-D$_2$, 1.5]) par

$$R(\mathcal{N}_h, \rho) = \rho\omega \min_{0 \leq i \leq m-1} |a_i|_\rho^{-1/m-i}.$$

Cette expression montre que la fonction $R(\mathcal{N}_h, \rho)$ est, *par morceaux*, de la forme $C_h\rho^{\beta_h}$ où $C_h$ est une constante. Le point est que le dénominateur de $\beta_h$ est *borné par l'ordre* $m$. Par image inverse la fonction $R(\mathcal{M}, \rho)$ a les mêmes propriétés, elle est logarithmiquement linéaire par morceaux dans chaque intervalle $]r_{h-1}, r_h[$. En particulier on peut trouver un point $r_h'$ de l'intervalle $]r_{h-1}, r_h[$ où la dérivée logarithmique de la fonction $R(\mathcal{M}, \rho)$ vaut $\beta_h$, un nombre rationnel dont le dénominateur est *borné par* $m$. La suite des nombres $\beta_h$ est décroissante minorée par 1 quand $h$ tend vers l'infini en vertu de la logarithme-concavité de la fonction $R(\mathcal{M}, \rho)$. Cette suite est donc stationnaire de valeur $\beta + 1$. La fonction $R(\mathcal{M}, \rho)$ est nécessairement égale à $C\rho^{\beta+1}$ dans un intervalle $]r - \varepsilon', r[$ pour un $\varepsilon' > 0$ en vertu de la logarithme-concavité. La constante $C$ vaut $r^{-\beta}$.

*Définition* 4.2-2. Soit $\mathcal{M}$ un $\mathcal{R}_K(r)$-module libre de rang $m \geq 1$ à connexion soluble en $r$; nous appelons plus grande pente $p$-adique de $\mathcal{M}$ le nombre rationnel $pt(\mathcal{M}) := \beta \geq 0$ défini dans le théorème 4.2-1.

*Remarque* 4.2-3. On peut considérer les intervalles $]r_h', r_h[$ où l'on a les inégalités

$$\rho\omega'^{1/p^{h-1}} < R(\mathcal{M}, \rho) < \rho\omega^{1/p^h}$$

pour $\rho \in ]r_h', r_h[$. Si l'intervalle $]r_h', r_h[$ n'est pas vide la partie 1) du théorème 4.1-3 est suffisante pour avoir un antécédent de Frobenius de Christol-Dwork [C-D$_2$] d'ordre $h$ sans singularités dans un intervalle. Mais pour que $]r'_h, r_h[$ soit non vide il est nécessaire que $p < p(p-1)$ ce qui exclu le seul cas $p = 2$. Donc seul le cas de la caractéristique 2 nécessite la partie 2) du théorème 4.1-3 pour la démonstration du théorème 4.2-1. De toute façon on conjecture que la condition $\omega\rho < R(\mathcal{M}, \rho)$ dans un intervalle est suffisante pour l'existence d'un antécédent de Frobenius de Christol-Dwork sans singularité.

## 5. Majorations et minorations explicites

Ce paragraphe est consacré à la démonstration du théorème de majoration et minoration explicites 5.3-3. La majoration nous a été suggérée par le théorème de majorations explicites de Dwork-Robba [D-R$_2$] alors que la minoration nous a été suggérée par le théorème de réduction de Robba [R$_4$] pour les modules différentiels de rang un. Nous fixons le corps $K$ de base.

Si $\mathcal{B}$ est une base d'un module différentiel nous notons quand il n'y a pas de risque de confusion Mat($\mathcal{B}$) la matrice de la connexion dans cette base et Mat($\mathcal{B}, \mathcal{G}$) la matrice de changement de bases. Plus généralement nous notons Mat($P, \mathcal{B}$) la matrice d'un opérateur différentiel $P$ dans la base $\mathcal{B}$.

Si la fonction $f$ est logarithmiquement dérivable à droite (resp. à gauche) en un point $\rho$, nous noterons $d\log^+ f(\rho)$ (resp. $d\log^- f(\rho)$) la dérivée à droite (resp. à gauche) de la fonction $\widetilde{f}$ au point Log($\rho$). En particulier, si la fonction $f$ est dérivable en $\rho$, on trouve:

$$d\log f(\rho) = \rho \, \frac{f'}{f}(\rho).$$

Nous notons $\|A\|(\rho)$ la norme en $\rho$ d'une matrice dont les coefficients sont des fonctions analytiques dans un intervalle.

THÉORÈME 5.0-4. *Soit $c(\rho)$ une fonction logarithmiquement concave sur l'intervalle $]r_1, r_2[$ et soient $\{A_i\}_{1 \leq i \leq n}$ une famille finie de matrices $m \times m$ à coefficients dans $\mathcal{A}_{]r_1, r_2[}$ telles que, pour tout nombre $\rho$ de l'intervalle $]r_1, r_2[$, on ait $\max_{1 \leq i \leq n}(\|A_i\|(\rho)) \geq c(\rho)$. Alors il existe des nombres $\lambda_i$ prenant les valeurs 0 ou 1, tels que, pour tout nombre $\rho$ de l'intervalle $]r_1, r_2[$, on ait $\|\sum_{i=1}^n \lambda_i A_i\|(\rho) \geq c(\rho)$.*

*Démonstration.* On fait une récurrence sur le nombre $n$ d'éléments de la famille. Pour $n = 1$, le résultat est évident. Supposons le résultat démontré pour les familles contenant $n-1$ matrices. La fonction $\|A_n\|(\rho)$ est logarithmiquement convexe (cf. [R-C]) et donc la fonction $\|A_n\|(\rho)/c(\rho)$ est logarithmiquement convexe sur l'intervalle $]r_1, r_2[$. L'ensemble des nombres $\rho$ où elle est strictement inférieure à 1 est donc un intervalle ouvert $]r_3, r_4[$, éventuellement vide, contenu dans $]r_1, r_2[$. On a donc, pour $r_3 < \rho < r_4$:

$$\max_{1 \leq i \leq n-1}(\|A_i\|(\rho)) \geq c(\rho).$$

D'après l'hypothèse de récurrence, il existe des nombres $\lambda_i$ prenant les valeurs 0 ou 1, tels que, pour tout nombre $\rho$ de l'intervalle $]r_3, r_4[$ on ait:

$$\| \sum_{i=1}^{n-1} \lambda_i A_i \|(\rho) \geq c(\rho).$$

Posons $B = \sum_{i=1}^{n-1} \lambda_i A_i$. On a $\|B\|(\rho) \geq c(\rho)$ si $\rho \in ]r_3, r_4[$ et, par construction de $]r_3, r_4[$, $\|A_n\|(\rho) \geq c(\rho)$ pour $\rho \in ]r_1, r_3] \cup [r_4, r_2[$. Donc

$$\max(\|B\|(\rho), \|A_n\|(\rho)) \geq c(\rho)$$

pour tout nombre $\rho$ de $]r_1, r_2[$. Pour conclure, il suffit donc de savoir traiter le cas de deux matrices. Ce sera l'objet du lemme suivant.

LEMME 5.0-5.    *Soit $c(\rho)$ une fonction logarithmiquement concave sur l'intervalle $]r_1, r_2[$ et soient $A$ et $B$ deux matrices $m \times m$ à coefficients dans $\mathcal{A}_{]r_1, r_2[}$ telles que, pour tout nombre $\rho$ de l'intervalle $]r_1, r_2[$, on ait $\max(\|A\|(\rho), \|B\|(\rho)) \geq c(\rho)$. Alors il existe des nombres $\lambda$ et $\mu$ prenant les valeurs 0 ou 1, tels que, pour tout nombre $\rho$ de l'intervalle $]r_1, r_2[$ on ait $\|\lambda A + \mu B\|(\rho) \geq c(\rho)$.*

*Démonstration.* Si l'une des deux fonctions (logarithmiquement convexe) $\|A\|(\rho)/c(\rho)$ ou $\|B\|(\rho)/c(\rho)$ est supérieure à 1 sur tout l'intervalle $]r_1, r_2[$, le résultat est évident. Sinon, quitte à échanger les rôles de $A$ et $B$, il existe des nombres $r_1 \leq r_3 < r_4 \leq r_5 < r_6 \leq r_2$ tels que:

$$\|A\|(\rho)/c(\rho) < 1 \text{ si et seulement si } \rho \in ]r_3, r_4[,$$

$$\|B\|(\rho)/c(\rho) < 1 \text{ si et seulement si } \rho \in ]r_5, r_6[.$$

On a alors:

$$d\log^+\Big(\|A\|(r_4)/c(r_4)\Big) \geq d\log^-\Big(\|A\|(r_4)/c(r_4)\Big) > 0,$$

$$d\log^-\Big(\|B\|(r_5)/c(r_5)\Big) \leq d\log^+\Big(\|B\|(r_5)/c(r_5)\Big) < 0.$$

En particulier, la fonction $\|A\|(\rho)/c(\rho)$ est logarithmiquement strictement croissante sur l'intervalle $]r_4, r_2]$ et la fonction $\|B\|(\rho)/c(\rho)$ logarithmiquement strictement décroissante sur l'intervalle $[r_1, r_5[$. Comme, de plus, on a:

$$\|A\|(r_4)/c(r_4) = \|B\|(r_5)/c(r_5) = 1$$

la fonction $\|A\|(\rho)/\|B\|(\rho)$ est logarithmiquement strictement croissante sur l'intervalle $[r_4, r_5]$ et passe d'une valeur inférieure à 1 à une valeur supérieure

à 1. Elle prend donc la valeur 1 en un unique point $r_7$ de cet intervalle. Compte tenu de la définition des points $r_i$, on en déduit que:

$$\|A\|(\rho) < \|B\|(\rho) \text{ pour } \rho \in [r_3, r_7[,$$
$$\|A\|(\rho) > \|B\|(\rho) \text{ pour } \rho \in ]r_7, r_6].$$

Comme la fonction $\|A + B\|(\rho)$ est continue, on trouve:

$$\|A + B\|(\rho) = \max\left(\|A\|(\rho), \|B\|(\rho)\right) \quad \text{pour } \rho \in [r_3, r_6].$$

En particulier on a:

$$\|A + B\|(\rho)/c(\rho) = \|B\|(\rho)/c(\rho),$$

pour $\rho \in [r_3, r_7]$ et la fonction $\|A + B\|(\rho)/c(\rho)$ est logarithmiquement décroissante sur l'intervalle $[r_3, r_7]$, c'est-à-dire, par logarithmique convexité, sur l'intervalle $]r_1, r_7[$. On démontre de même qu'elle est strictement croissante sur l'intervalle $]r_7, r_2[$. Par ailleurs, comme on a:

$$\|A + B\|(r_7)/c(r_7) = \|A\|(r_7)/c(r_7) \geq \|A\|(r_4)/c(r_4) \geq 1,$$

on a montré que:

$$\|A + B\|(\rho)/c(\rho) \geq 1 \quad \text{pour } \rho \in ]r_1, r_2[.$$

5.1. *Majoration de la matrice de passage.* Soit $\mathcal{M}$ un $\mathcal{A}_{]r-\varepsilon, r[}$-module différentiel libre de rang $m$ soluble en $r$. En vertu du théorème 4.2-1 on peut supposer, pour $\varepsilon > 0$ assez petit, que $R(\mathcal{M}, \rho) = \rho(\rho/r)^\beta$ pour $r - \varepsilon < \rho < r$ où $\beta = \mathrm{pt}(\mathcal{M})$ est la plus grande pente de $\mathcal{M}$ en $r$. Nous supposons dans ce paragraphe la pente $\beta > 0$.

Soit $\rho$ un nombre de $]r - \varepsilon, r[$ fixé. On définit un entier $\lambda$ par:

$$\omega^p < (\rho/r)^{\beta p^\lambda} \leq \omega$$

de sorte que:

$$(1) \qquad \omega^p \, \rho^{p^\lambda} < R(\mathcal{M}, \rho)^{p^\lambda} \leq \omega \, \rho^{p^\lambda},$$

c'est-à-dire, en posant:

$$d(\rho/r) = -(p-1)\beta \, \frac{\mathrm{Log}(\rho/r)}{\mathrm{Log}(p)},$$

$$\frac{1}{d(\rho/r)} \leq p^\lambda < \frac{p}{d(\rho/r)}, \qquad \lambda = \left[-\frac{\mathrm{Log}(d(\rho/r))}{\mathrm{Log}(p)}\right] > 0.$$

On a noté par $[e]$ la partie entière d'un nombre réel $e$.

Notons $\varphi$ le morphisme de Frobenius, la ramification d'ordre $p$, de la couronne $C(]r - \varepsilon, r[)$ dans la couronne $C(](r - \varepsilon)^p, r^p[)$. Le théorème de structure de Frobenius de Christol-Dwork [C-D$_2$] sur une circonférence montre qu'il existe un $E_{\rho^{p^\lambda}}$-module différentiel $\mathcal{N}$ tel que $\varphi^{\lambda*}(\mathcal{N}) \simeq \mathcal{M} \otimes E_\rho$, où $E_\rho$ est le complété pour la norme en $\rho$ du corps $K(x)$. De plus, si $\mathcal{G}$ est une base de $\mathcal{M}$ on construit une base cyclique $\mathcal{B}$ de $\mathcal{N}$ telle que, si on pose:

$$B := \mathrm{Mat}(\mathcal{B}), \quad G := \mathrm{Mat}(\mathcal{G}), \quad F := \mathrm{Mat}(\varphi^{\lambda*}(\mathcal{B})), \quad H := \mathrm{Mat}(\mathcal{G}, \varphi^{\lambda*}(\mathcal{B})),$$

on trouve les majorations:

$$\|H\|(\rho)\,\|H^{-1}\|(\rho) \leq c(m)^{-\lambda}\,|(m-1)!|^{-\lambda-1}\,p^{(m-1)(\lambda+2)}\,\max(1, \|G\|(\rho))^{(m-1)}$$

$$\leq p^{c_1\lambda}c_2(\rho)$$

où $c(m)$ est n'importe quelle constante strictement inférieure à $|\prod_{j=1}^{m-1}\binom{m}{j}|$, où $c_1 > 0$ est une constante ne dépendant que de $m$ et de la base $\mathcal{G}$ et où $c_2(\rho)$ est une fonction logarithmiquement convexe de $\rho$.

5.2. *Passage de $\mathcal{B}$ à $\varphi^{\lambda*}(\mathcal{B})$.*

LEMME 5.2-1.   *Les entiers $\alpha_{k,s}$ définis par*:

$$\left((x+1)^{p^\lambda}-1\right)^s = \sum_{k=s}^{sp^\lambda} \alpha_{k,s}\,x^k$$

*vérifient, pour tout entier $\delta$ tel que $0 \leq \delta \leq \lambda - 1$:*

$$(2) \qquad |\alpha_{k,s}| \leq \omega^{sp+\delta(p-1)s-kp^{\delta+1-\lambda}}$$

*avec égalité si $k = s\,p^{\lambda-\delta-1}$ ou si $k = s\,p^{\lambda-\delta}$. En particulier $|\alpha_{sp^\lambda,s}| = 1$ et $|\alpha_{p^{\lambda-\delta},1}| = |p|^\delta$.*

*Démonstration.* Des égalités:

$$(x+1)^{p^\lambda}-1 = \sum_{k=1}^{p^\lambda}\binom{p^\lambda}{k}x^k, \qquad \left|\binom{p^\lambda}{k}\right| \leq |p|^\nu \text{ pour } k < p^{\lambda-\nu+1},$$

il résulte que, pour:

$$|x| = \omega^{p^{\delta+1-\lambda}}$$

on a:

$$\left|\sum_{k=s}^{sp^\lambda}\alpha_{k,s}\,x^k\right| \leq \left(|p|^\delta\,|x|^{p^{\lambda-\delta}}\right)^s = \omega^{s(p+\delta(p-1))};$$

la majoration (2) s'en déduit immédiatement. Les cas particuliers se vérifient facilement. Ils peuvent aussi très facilement s'obtenir directement.

Posons:

$$B_s := \mathrm{Mat}(\Delta^s, \mathcal{B}), \qquad F_s := \mathrm{Mat}(\Delta^s, \varphi^{\lambda*}(\mathcal{B})).$$

Proposition 5.2-2. *On a, pour $r < \rho < 1$:*

(3)
$$\begin{cases} \|B_s\|(\rho^{p^\lambda}) \leq \omega\, r^{-sp^\lambda}(\rho/r)^{-(\beta+1)sp^\lambda} & \text{pour } s \geq 1 \\ \|B_{p^h}\|(\rho^{p^\lambda}) = \omega\, r^{-p^{h+\lambda}}(\rho/r)^{-(\beta+1)p^{h+\lambda}} & \text{pour } h \geq 0. \end{cases}$$

*Démonstration.* Comme $\mathcal{B}$ est une base cyclique et comme:

$$R(\mathcal{N}, \rho^{p^\lambda}) = R(\mathcal{M}, \rho)^{p^\lambda} \leq \omega\rho^{p^\lambda}$$

on a, d'après ([C-D$_2$, Prop. 3.2 et Lemme 1.4]):

$$\begin{aligned} \|B_s\|(\rho^{p^\lambda}) = \|B\|^s(\rho^{p^\lambda}) &= \left|\frac{1}{s!}\right| \left|\left(\frac{R(\mathcal{N}, \rho^{p^\lambda})}{\omega}\right)^{-s}\right| \\ &= \left|\frac{1}{s!}\right| \left|\omega^s \rho^{-sp^\lambda}(\rho/r)^{-\beta sp^\lambda}\right|. \end{aligned}$$

Il suffit alors de remarquer que:

$$|\frac{1}{s!}|\omega^s \leq \omega \quad \text{pour } s \geq 1$$

avec égalité pour $s = p^h$, pour démontrer la proposition.

Proposition 5.2-3. *On a, pour $r - \varepsilon < \rho < r$:*

$$\begin{cases} \|F_k\|(\rho) \leq \omega\, \rho^{-k}(\rho/r)^{-\beta k} & \text{pour } k \geq 1 \\ \|F_{p^h}\|(\rho) \geq \omega\, |p|^\lambda\, \rho^{-p^h}(\rho/r)^{-\beta p^h} & \text{pour } h \geq 1. \end{cases}$$

*Démonstration.* La solution de l'équation différentielle:

$$X' = B\,X, \qquad X(t_\rho^{p^s}) = I$$

au voisinage du point générique $t_\rho^{p^\lambda}$ (avec $|t_\rho^{p^\lambda}| = \rho^{p^\lambda}$) est :

$$X(x) = \sum_{s=0}^\infty B_s(t_\rho^{p^\lambda})\,\left(x - t_\rho^{p^\lambda}\right)^s.$$

Par définition de la base $\varphi^{\lambda*}(\mathcal{B})$, la solution de l'équation différentielle:

$$Y' = F\,Y, \qquad Y(t_\rho) = I$$

au voisinage du point générique $t_\rho$ $(|t_\rho| = \rho)$ est:

$$\begin{aligned} Y(x) &= \sum_{s=0}^\infty F_s(t_\rho)\,(x - t_\rho)^s \\ &= X(x^{p^\lambda}) = \sum_{s=0}^\infty B_s(t_\rho^{p^\lambda})\,\left(x^{p^\lambda} - t_\rho^{p^\lambda}\right)^s. \end{aligned}$$

Remarquant que:

$$\left(x^{p^\lambda} - t_\rho^{p^\lambda}\right)^s = \sum_{k=s}^{sp^\lambda} t_\rho^{p^\lambda s} \; \alpha_{k,s} \left(\frac{x - t_\rho}{t_\rho}\right)^k,$$

on trouve donc:

$$F_k(t_\rho) = \sum_{p^{-\lambda}k \leq s \leq k} B_s(t_\rho^{p^\lambda}) \, t_\rho^{p^\lambda s - k} \; \alpha_{k,s},$$

c'est-à-dire, d'après les majorations (2) et (3), pour tout entier $\delta$ tel que $0 \leq \delta \leq \lambda - 1$, on trouve

$$
\begin{aligned}
(4) \qquad \|F_k\|(\rho) &\leq \max_{p^{-\lambda}k \leq s \leq k} \|B_s\|(\rho^{p^\lambda}) \, \rho^{p^\lambda s - k} \; |\alpha_{k,s}| \\
&\leq \max_{p^{-\lambda}k \leq s \leq k} (\omega \, (\rho/r)^{-\beta s p^\lambda} \; \rho^{-k} \; \omega^{sp + \delta(p-1)s - kp^{\delta+1-\lambda}}).
\end{aligned}
$$

Par ailleurs la minoration (1) montre:

$$(\rho/r)^{-\beta p^\lambda} \, \omega^{p+\delta(p-1)} < \omega^{-p+p+\delta(p-1)} \leq 1.$$

Les termes apparaissant dans la majoration (4) sont donc, pour chaque valeur de $\delta$, strictement décroissants en fonction de $s$. On en déduit que le maximum est atteint pour la plus petite valeur de $s$ possible, c'est-à-dire pour le plus petit entier supérieur à $p^{-\lambda}k$.

En particulier en prenant $\delta = 0$ et $s = p^{-\lambda}k$ on trouve:

$$\|F_k\|(\rho) \leq \omega \, \rho^{-k} (r/\rho)^{\beta k} \, \omega^{p^{-\lambda}kp - kp^{1-\lambda}} = \omega \, \rho^{-k} r/(\rho)^{\beta k}.$$

Pour obtenir la minoration, nous distinguons deux cas.

5.2.1. *Cas où $h > \lambda$.* On prend $\delta = 0$. Le plus grand terme dans la majoration (4) est obtenu pour $s = p^{h-\lambda}$. Les cas d'égalité dans le lemme 5.2-1 et la proposition 5.2-2 montrent que l'on a:

$$\|F_{p^h}\|(\rho) = \|B_{p^{h-\lambda}}\|(\rho^{p^\lambda}) \, |\alpha_{p^h, p^{h-\lambda}}| = \omega \, \rho^{-p^h} r/(\rho)^{\beta p^h}.$$

5.2.2. *Cas où $1 \leq h \leq \lambda$.* On prend $\delta = \lambda - h$. Comme $0 < p^{h-\lambda} \leq 1$, le plus grand terme dans la majoration (4) est obtenu pour $s = 1$. On trouve donc:

$$\|F_{p^h}\|(\rho) = \|B_1\|(\rho^{p^\lambda}) \, \rho^{p^\lambda - p^h} \; |\alpha_{p^h, 1}| = \omega \, \rho^{-(\beta+1)p^\lambda + p^\lambda - p^h} \, |p|^\delta r^{\beta p^\lambda}.$$

La minoration de la proposition 5.2-3 découle immédiatement de ce que $\beta(p^h - p^\lambda) \leq 0$, $\rho/r < 1$ et $|p|^\delta \geq |p|^\lambda$.

5.3. *Passage de $\varphi^{\lambda *}(\mathcal{B})$ à $\mathcal{G}$.* On pose:

$$G_s := \mathrm{Mat}(\Delta^s, \mathcal{G}).$$

PROPOSITION 5.3-1.    *Il existe une fonction logarithmiquement concave* $c(\rho)$ *et une fonction* $M(\rho)$ *telles que l'on ait, pour* $r - \varepsilon < \rho < r$:

$$\rho^s \|G_s\|(\rho) \quad \leq \quad M(\rho)\,(\rho/r)^{-\beta s} \quad \text{pour } s \geq 0,$$

$$\max_{0 \leq s \leq p^h} \rho^s\,\|G_s\|(\rho) \quad \geq \quad c(\rho)\,(\rho/r)^{(-\beta)p^h} \quad \text{pour } h \geq 1.$$

*Démonstration.* La formule de Leibniz s'écrit:

$$G_s = \sum_{k=0}^{s} \frac{1}{k!}\,H^{(k)}\,F_{s-k}\,H^{-1}, \qquad F_s = \sum_{k=0}^{s} \frac{1}{k!}\,(H^{-1})^{(k)}\,G_{s-k}\,H$$

et donne les majorations:

$$\|G_s\|(\rho) \leq \|H\|(\rho)\,\|H^{-1}\|(\rho)\,\max_{0 \leq k \leq s} \rho^{-k}\,\|F_{s-k}\|(\rho),$$

$$\|F_s\|(\rho) \leq \|H\|(\rho)\,\|H^{-1}\|(\rho)\,\max_{0 \leq k \leq s} \rho^{-k}\,\|G_{s-k}\|(\rho).$$

Compte tenu de la proposition 5.2-3, la première de ces majorations s'écrit (on a $\beta > 0$ et $\rho < r$):

$$\begin{aligned}
\|G_s\|(\rho) \quad &\leq \quad p^{c_1\lambda} c_2(\rho)\omega \max_{0 \leq k \leq s} r^{\beta(s-k)}\rho^{-k-(\beta+1)(s-k)} \\
&\leq \quad \left(\frac{p}{d(\rho/r)}\right)^{c_1} c_2(\rho)\,\rho^{-(\beta+1)s}\omega r^{\beta s}
\end{aligned}$$

qui est bien de la forme annoncée.

Pour $s = p^h$ la deuxième majoration s'écrit:

$$\begin{aligned}
\max_{0 \leq k \leq p^h} \rho^{-p^h+k}\,\|G_k\|(\rho) \quad &\geq \quad p^{-c_1\lambda} c_2(\rho)^{-1}\,\omega\,|p|^{\lambda}\,\rho^{-(\beta+1)p^h}\,r^{\beta p^h} \\
&\geq \quad \omega\left(\frac{d(\rho/r)}{p}\right)^{1+c_1} c_2(\rho)^{-1}\,\rho^{-(\beta+1)p^h}r^{(\beta)p^h}.
\end{aligned}$$

Il suffit alors de remarquer que la fonction:

$$c(\rho) = \omega\left(\frac{d(\rho/r)}{p}\right)^{1+c_1} c_2(\rho)^{-1} = c_3(-\mathrm{Log}(\rho/r))^{1+c_1} c_2(\rho)^{-1}$$

est logarithmiquement concave sur $]0, r[$ parce que $1 + c_1 > 0$.

*Remarque* 5.3-2. La majoration de la proposition précédente montre que si $\mathcal{M}$ est un $\mathcal{A}_{|r-\varepsilon,r[}$-module différentiel tel que $R(\mathcal{M}, \rho) < \rho$, alors pour $a$ tel que $|a| = \rho$ les solutions locales de $\mathcal{M}$, au voisinage de $a$, sont des fonctions bornées dans le disque $D(a, R(\mathcal{M}, \rho)^-)$. On sait que ce résultat n'est pas vrai en général si $R(\mathcal{M}, \rho) = \rho$.

THÉORÈME 5.3-3. *Soit $\mathcal{M}$ un $\mathcal{A}_{]r-\varepsilon,r[}$-module différentiel libre de rang $m$ tel que $R(\mathcal{M},\rho) = r(\rho/r)^{\beta+1}$ pour $r - \varepsilon < \rho < r$, $\beta > 0$, et soit $\mathcal{G}$ une base de $\mathcal{M}$. Il existe des nombres $(\lambda_{s,h})_{0 \le s \le p^h}$ prenant les valeurs 0 ou 1, une fonction $c$ logarithmiquement concave et une fonction $M$, tels que si on pose*:*

$$L_h = \sum_{s=0}^{p^h} \lambda_{s,h}\, x^s \Delta^s,$$

*on ait, pour $r - \varepsilon < \rho < r$ les inégalités*:*

$$c(\rho)\,(\rho/r)^{-\beta p^h} \le \|\operatorname{Mat}(L_h,\mathcal{G})\|(\rho) \le M(\rho)\,(\rho/r)^{-\beta p^h}.$$

*Démonstration.* D'après la proposition 5.3-1, on a, pour $0 \le s \le p^h$, d'une part:

$$\|\operatorname{Mat}(x^s\Delta^s,\mathcal{G})\|(\rho) \le M(\rho)\,(\rho/r)^{-\beta s} \le M(\rho)\,(\rho/r)^{-\beta p^h}$$

ce qui, compte tenu de $|\lambda_{s,h}| \le 1$, donne la majoration du théorème, et, d'autre part:

$$\max_{0 \le s \le p^h} \|\operatorname{Mat}(x^s\Delta^s,\mathcal{G})\|(\rho) \ge c(\rho)\,(\rho/r)^{-\beta p^h}.$$

Si on remarque que la fonction $c(\rho)\,(\rho/r)^{-\beta p^h}$ est logarithmiquement concave, il suffit d'appliquer le théorème 5.0-4 pour pouvoir conclure.

*Exemple* 5.3-4 . Supposons le rang $m$ de $\mathcal{M}$ égal à un. Pour tout entier $s$ on a la majoration explicite de Dwork-Robba [D-R$_2$]

$$\|\operatorname{Mat}(x^s\Delta^s,\mathcal{G})\|(\rho) \le (\rho/r)^{-\beta s}$$

qui fait apparaître la majoration du théorème 5.3-3 comme un cas particulier.

Supposons que dans une base $\mathcal{G}$ la matrice de la connexion s'écrive $a_d\pi/x^{d+1} + g(x)$ où $d \ge 1$ est un entier, $g(x)$ une fonction méromorphe dans le disque $D(0,r^-)$ ayant au plus un pôle en zéro d'ordre $d$ et tel que $|a_d| = r^d$. Par un calcul direct on trouve que

$$\operatorname{Mat}(x^{p^h}\Delta^{p^h},\mathcal{G}) = \frac{(\pi a_d)^{p^h}}{p^h! x^{(d+1)p^h}} + g_h(x)$$

où $g_h(x)$ est une fonction méromorphe dans le disque $D(0,r^-)$ ayant au plus un pôle en zéro d'ordre $(d+1)p^h - 1$. Cette expression donne alors la minoration

$$\omega(r/\rho)^{dp^h} \le \|\operatorname{Mat}(x^{p^h}\Delta^{p^h},\mathcal{G})\|(\rho).$$

Le théorème de réduction de Robba [R$_4$] montre, quitte à passer à une extension maximalement complète assez grande, qu'on peut toujours trouver une telle base pour laquelle $d = pt(\mathcal{M})$, la pente de $\mathcal{M}$ en $r$.

## 6. Le théorème de décomposition selon les pentes $p$-adiques

6.1. *Décomposition par rapport à la plus grande pente.* Soit $K$ un corps de caractéristique nulle valué complet pour une valeur absolue $p$-adique, $r > 0$ un nombre réel et $\mathcal{M}$ un $\mathcal{R}_K(r)$-module libre de rang $m$ muni d'une connexion. Le module différentiel $\mathcal{M}$ est défini dans un intervalle $[r - \varepsilon_0, r[$ pour un $\varepsilon_0 > 0$. Pour tout $\varepsilon, 0 < \varepsilon \le \varepsilon_0$, nous notons $\mathcal{M}_{[r-\varepsilon,r[}$ sa restriction à la couronne $C([r - \varepsilon, r[)$. C'est un $\mathcal{A}_K([r - \varepsilon, r[)$-module libre de rang $m$ muni d'une connexion. Le choix d'une base $\mathcal{G}$ de $\mathcal{M}_{[r-\varepsilon_0,r[}$ définit une suite exacte
$(*)$
$$0 \to (\mathcal{A}_K([r - \varepsilon, r[)[\frac{d}{dx}])^m \xrightarrow{u := \frac{d}{dx} - G(x)} (\mathcal{A}_K([r - \varepsilon, r[)[\frac{d}{dx}])^m \to \mathcal{M}_{[r-\varepsilon,r[} \to 0$$

pour tout $\varepsilon \le \varepsilon_0$ où $G(x)$ est la matrice de la connexion dans la base $\mathcal{G}$.

Faisons l'hypothèse que $\mathcal{M}$ est soluble en $r$ de sorte qu'en vertu du théorème 4.2-1 la fonction rayon de convergence $R(\mathcal{M}, \rho)$ est égale à $r(\rho/r)^{\mathrm{pt}(\mathcal{M})+1}$ pour tout $\rho \in [r - \varepsilon, r[$, pour un $\varepsilon > 0$ assez petit, où la plus grande pente $\mathrm{pt}(\mathcal{M})$ est un nombre rationnel $\ge 0$. Si la plus grande pente $\mathrm{pt}(\mathcal{M})$ est égal à $0$ le module différentiel $\mathcal{M}$ a la propriété de Robba au bord et sa structure a été étudiée en détail dans les articles [C-M$_1$], [C-M$_2$].

Nous supposons que la plus grande pente $\mathrm{pt}(\mathcal{M})$ de $\mathcal{M}$ est *strictement plus grande que* $0$. Soit alors un nombre réel $\gamma$ tel que $0 \le \gamma < \mathrm{pt}(\mathcal{M})$ et considérons la topologie quotient $\mathcal{T}_{\gamma,Q}$ sur $\mathcal{M}_{[r-\varepsilon,r[}$.

THÉORÈME 6.1-1. *Supposons que le corps $K$ est localement compact et que le nombre $r$ appartient au groupe des valeurs absolues du corps de base $K$. Sous les conditions précédentes, pour tout $\varepsilon$ assez petit la topologie quotient $\mathcal{T}_{\gamma,Q}$ sur $\mathcal{M}_{[r-\varepsilon,r[}$ n'est pas séparée. Autrement dit l'adhérence $\bar{0}_\gamma(\mathcal{M}_{[r-\varepsilon,r[})$ de zéro dans $\mathcal{M}_{[r-\varepsilon,r[}$ n'est pas nulle et est un $\mathcal{A}_K([r - \varepsilon, r[)$-module libre de rang fini à connexion.*

LEMME 6.1-2. *Soit $I$ un intervalle ouvert de $\mathbb{R}^+$; si le corps $K$ est localement compact on peut extraire de toute suite bornée pour la tolopogie d'espace de type $\mathcal{F}$ de fonctions analytiques de $\mathcal{A}_K(I)$ une sous-suite convergente.*

*Démonstration.* Soit $I$ un intervalle $]r, R[$; en vertu de la décomposition

$$\mathcal{A}_K(I) = \mathcal{A}_K([0, R[) \oplus \frac{1}{x} \mathcal{A}_K(]r, \infty])$$

on peut supposer que $I$ est de la forme $[0, r[$ pour un réel $r > 0$. Soit $f_n$ une suite bornée de fonctions analytiques dans la couronne $C([0, r[)$; nous allons extraire une sous suite $f_{v(n)}$ de Cauchy pour la topologie d'espace $\mathcal{F}$ de $\mathcal{A}_K([0, r[)$ qui sera donc convergente. Soit une suite strictement croissante de nombre rationnels $r_k$ qui tend vers $r$ par valeurs inférieures. Nous allons d'abord extraire une sous-suite $f_{v_1}$ de Cauchy pour la norme $|,|_{r_1}$. Si $f_n(x) = \sum_{k=0,\infty} a_k^n x^k$ les nombres réels $|a_k^n| r_2^k$ sont bornés par une constante $C_2$ par hypothèse. Quitte à passer à une extension finie $K_2$ de $K$, on peut supposer que $r_2$ est la valeur absolue d'un nombre $b_2$ de $K_2$. La boule $D(0, C_2^+)$ de $K_2$ est *compacte* et en vertu du théorème de Tychonoff le produit $\prod_{k \in \mathbb{N}} D(0, C_2^+)$ est un métrique compact. De la suite $g_n := (a_k^n b_2^k)_{k \in \mathbb{N}}$ de $\prod_{k \in \mathbb{N}} D(0, C_2^+)$ on peut alors extraire une sous-suite convergente $g_{v_1}$. Nous affirmons que la sous-suite $f_{v_1}$ est de Cauchy pour la norme $|,|_{r_1}$. En effet pour un $\varepsilon > 0$ les nombres $|a_k^{v_1(n)} - a^{v_1(n')}| r_1^k = |a_k^{v_1(n)} b_2^k - a^{v_1(n')} b_2^k| \frac{r_1^k}{r_2^k}$ sont plus petits que $\varepsilon$ pour $k \geq k_0$ et les nombres $|a_k^{v_1(n)} - a^{v_1(n')}| r_1^k = |a_k^{v_1(n)} b_2^k - a^{v_1(n')} b_2^k| \frac{r_1^k}{r_2^k}$ pour $k < k_0$ sont plus petits que $\varepsilon$ dès que $n$ et $n'$ sont plus grand que $n(\varepsilon)$ par définition de la topologie produit.

La sous-suite $f_{v_1}$ est bornée pour la norme $|,|_{r_3}$ et par le procédé précédent on peut extraire une sous suite $f_{v_2}$ de Cauchy pour la norme $|,|_{r_2}$. Par récurrence on construit ainsi une suite de sous-suites $f_{v_n}$ telle que la suite $f_{v_{n+1}}$ est une sous-suite de $f_{v_n}$ et de Cauchy pour la norme $|,|_{r_n}$. La sous-suite diagonale $f_{v_n(n)}$ de $f_n$ est de Cauchy pour toutes les normes $|,|_{r_n}$. Elle est donc convergente dans l'espace complet $\mathcal{A}_K([0, r[)$. D'où le lemme.

*Démonstration du théorème* 6.1-1. Nous utilisons les notations précédentes. Si $a$ est un nombre de $K$ tel que $|a| = r$ le changement de variables $x \mapsto x/a$ permet de supposer que $r = 1$, ce que nous ferons pour simplifier les notations. Nous supposons que $R(\mathcal{M}, \rho)$ est égale à $\rho^{\text{pt}(\mathcal{M})+1}$ pour $\rho \in ]1, \infty[$. En vertu du théorème 5.3-3 il existe des nombres $(\lambda_{s,h})_{0 \leq s \leq p^h}$ prenant les valeurs 0 ou 1 tels que, si on pose:

$$L_h = \sum_{s=0}^{p^h} \lambda_{s,h} \; x^s \Delta^s,$$

on ait, pour $1 - \varepsilon < \rho < 1$:

$$(**) \qquad c(\rho) \, \rho^{-\text{pt}(\mathcal{M})p^h} \leq ||\text{Mat}(L_h, \mathcal{G})||(\rho) \leq M(\rho) \, \rho^{-\text{pt}(\mathcal{M})p^h}.$$

Notons $\mathcal{I}$ l'image du morphisme $u$ de la suite exacte $(*)$ de sorte que l'on a la somme directe algébrique d'espaces vectoriels sur $K$:

$$(\mathcal{A}_K([1-\varepsilon, 1[)[\frac{d}{dx}])^m = \mathcal{I} \oplus (\mathcal{A}_K([1-\varepsilon, 1[))^m$$

où $(\mathcal{A}_K([1-\varepsilon, 1[))^m$ est identifié à $\mathcal{M}_{[1-\varepsilon,1[}$ par la base $\mathcal{G}$.

Si on munit les espaces $(\mathcal{A}_K([1-\varepsilon, 1[))^m$ et $\mathcal{I}$ de la topologie induite, dans la décomposition précédente, par la topologie $\mathcal{T}_\gamma$ de $(\mathcal{A}_K([1-\varepsilon, 1[)[\frac{d}{dx}])^m$ cette somme directe *n'est pas* topologique en général.

Notons $p^h = q_h k + s_h, 0 \leq s_h < k$, la division de $p^h$ par le dénominateur $k$ de pt$(\mathcal{M})+1$. Puisque pt$(\mathcal{M})q_h k$ est un entier naturel, nous pouvons considérer l'opérateur différentiel $x^{\mathrm{pt}(\mathcal{M})q_h k} L_h$. Des inégalités $(**)$ on déduit les inégalités pour tout $\rho \in \,]1-\varepsilon,1[$:

$$(***) \qquad c(\rho)\, \rho^{-\mathrm{pt}(\mathcal{M})s_h} \leq ||\mathrm{Mat}(x^{\mathrm{pt}(\mathcal{M})q_h k}L_h, \mathcal{G})||(\rho) \leq M(\rho)\, \rho^{-\mathrm{pt}(\mathcal{M})s_h}.$$

La suite de matrices $\mathrm{Mat}(x^{\mathrm{pt}(\mathcal{M})q_h k}L_h, \mathcal{G})$ est à coefficients dans l'anneau $\mathcal{A}_K([1-\varepsilon,1[)$ en fait dans l'anneau $\mathcal{A}_K([1-\varepsilon_0,1[)$. En vertu de la majoration $(***)$ c'est une suite bornée dans l'espace $(\mathcal{A}_K(]1-\varepsilon,1[))^{m^2}$ pour la topologie naturelle d'espace $\mathcal{F}$ définie par la famille des normes $|,|_\rho, \rho \in \,]1-\varepsilon,1[$. En vertu du lemme précédent, le corps $K$ étant supposé localement compact, cette suite de matrices admet une sous-suite convergente dans l'espace $(\mathcal{A}_K(]1-\varepsilon,1[))^{m^2}$ indexée par $a(h)$, de limite $A$.

En vertu de la minoration $(***)$ la matrice carré limite $A$ d'ordre $m$ à coefficients dans l'anneau $\mathcal{A}_K(]1-\varepsilon,1[)$, n'est pas *nulle*. Elle admet donc une ligne disons, $A^{i_0}$, qui n'est pas nulle.

Considérons alors la suite de $m$-vecteurs ligne d'opérateurs différentiels $x^{\mathrm{pt}(\mathcal{M})q_{a(h)}k}L^{i_0}_{a(h)}$ dont les composantes sont toutes nulles à l'exception de la $i_0$-ème $x^{\mathrm{pt}(\mathcal{M})q_{a(h)}k}L_{a(h)}$ de décomposition dans $(\mathcal{A}_K([1-\varepsilon,1[)[\frac{d}{dx}])^m$:

$$x^{\mathrm{pt}(\mathcal{M})q_{a(h)}k}L^{i_0}_{a(h)} = u(Q_{a(h)}) + A^{i_0}_{a(h)}.$$

Quand $a(h)$ tend vers l'infini, par définition de la topologie $\mathcal{T}_\gamma$ pour $0 \leq \gamma <$ pt$(\mathcal{M})$ ce vecteur tend vers zéro dans l'espace $(\mathcal{A}_K([1-\varepsilon',1[)[\frac{d}{dx}])^m$ pour tout $\varepsilon', 0 < \varepsilon' < \varepsilon$ parce que on a la majoration:

$$(****) \qquad |x^{\mathrm{pt}(\mathcal{M})q_{a(h)}k}L_{a(h)}|_{\gamma,\rho} \leq \rho^{\mathrm{pt}(\mathcal{M})q_{a(h)}k - \gamma p^{a(h)}}.$$

D'autre part la suite $A^{i_0}_{a(h)}$ tend vers $A^{i_0}$ pour la topologie induite par $\mathcal{T}_\gamma$, qui est aussi la topologie naturelle d'espace $\mathcal{F}$, dans l'espace $(\mathcal{A}_K([1-\varepsilon',1[)$. Donc la suite $u(Q_{a(h)})$ tend vers $-A^{i_0}$ dans l'espace $(\mathcal{A}([1-\varepsilon',1[)[\frac{d}{dx}])^m$.

Ceci montre que $A^{i_0}$ qui n'est pas dans l'image de $u$ est adhérent à cette image dans l'espace $(\mathcal{A}([1-\varepsilon',1[)[\frac{d}{dx}])^m$. Sa classe dans $\mathcal{M}_{[1-\varepsilon',1[}$ n'est pas nulle et est adhérente à zéro pour la topologie quotient $\mathcal{T}_{\gamma,Q}$ pour tout $\gamma, 0 \leq \gamma <$ pt$(\mathcal{M})$ et tout $\varepsilon', 0 < \varepsilon' < \varepsilon$. L'espace $\bar{0}_\gamma(\mathcal{M}_{[1-\varepsilon',1[})$, adhérence de zéro dans $\mathcal{M}_{[1-\varepsilon',1[}$ pour la topologie quotient $\mathcal{T}_{\gamma,Q}$ pour tout $\gamma, 0 \leq \gamma <$ pt$(\mathcal{M})$ n'est pas nulle pour tout $\varepsilon' > 0$ assez petit. D'où le théorème 6.1-1 en tenant compte du théorème 3.2-1.

COROLLAIRE 6.1-3. *Supposons le corps $K$ localement compact et soit $\mathcal{M}$ un $\mathcal{A}_K([1-\varepsilon,1[)$-module libre de rang $m$ à connexion tel que $R(\mathcal{M},\rho) = \rho^{\mathrm{pt}(\mathcal{M})+1}$, $\mathrm{pt}(\mathcal{M}) > 0$, pour $\rho \in [1-\varepsilon,1[$ et un nombre réel $\gamma \geq 0$ tel que la topologie quotient $\mathcal{T}_{\gamma,Q}$ sur $\mathcal{M}$ est séparée. Alors $\mathrm{pt}(\mathcal{M}) \leq \gamma$.*

*Démonstration.* En effet si $0 \leq \gamma < \mathrm{pt}(\mathcal{M})$ la topologie $\mathcal{T}_{\gamma,Q}$ sur $\mathcal{M}$ n'est pas séparée en vertu du théorème 6.1-1.

*Remarque* 6.1-4. Sous les conditions du corollaire précédent si $\gamma \neq 0$ et si la topologie quotient $\mathcal{T}_{\gamma,Q}$ sur $\mathcal{M}$ est *séparée*, alors en vertu de la majoration 5.3-1, la topologie quotient est équivalente à la topologie naturelle d'espace $\mathcal{F}$ associée à une base de $\mathcal{M}$. En particulier la topologie quotient $\mathcal{T}_{\gamma,Q}$ sur le $\mathcal{R}_K(1)$-module différentiel associé est *séparée.*

Seul le cas de la topologie $\mathcal{T}_{0,Q}$ dans le cas d'une pente nulle échappe encore à notre connaissance malgré nos efforts:

CONJECTURE 6.1-5. *Soit un $\mathcal{R}_K(1)$-module différentiel $\mathcal{M}$ de pente nulle, alors la topologie $\mathcal{T}_{0,Q}$ sur $\mathcal{M}$ est séparée si et seulement si elle est équivalente à la topologie naturelle d'espace $\mathcal{LF}$ induite par une base.*

Cependant on a un cas particulier important pour les applications où cette conjecture est vérifiée. Rappelons d'abord que nous avons défini dans [C-M$_2$, §5] *l'exposant* $\mathfrak{E}xp_0^r(\mathcal{M})$ de $\mathcal{M}$ comme une classe de $(\mathbb{Z}_p/\mathbb{Z})^m$ modulo la relation d'équivalence $\overset{\mathfrak{E}}{\sim}$ où $m$ désigne le rang de $\mathcal{M}$.

PROPOSITION 6.1-6. *Soit $\mathcal{R}_K(1)$-module différentiel $\mathcal{M}$ de pente nulle dont l'exposant $\mathfrak{E}xp_0(\mathcal{M})$ a des différences qui ont la propriété (**NL**), alors la topologie $\mathcal{T}_{0,Q}$ sur $\mathcal{M}$ est séparée si et seulement si elle est équivalente à la topologie naturelle d'espace $\mathcal{LF}$ induite par une base.*

*Démonstration.* Dans cette situation en vertu du théorème fondamental de décomposition [C-M$_2$], $\mathcal{M}$ admet une forme normale sur $\mathbb{C}_p$ qui en vertu du théorème 7.1-2 se descend au corps de base $K$. On est réduit au cas où $\mathcal{M}$ est défini par un opérateur $(x\frac{d}{dx} - \alpha)^m$ où $\alpha$ est un entier $p$-adique. On voit dans ce cas par un calcul direct que la topologie $\mathcal{T}_{0,Q}$ est séparée si et seulement si $m = 1$.

Soit $\mathcal{M}_{[1-\varepsilon,1[}$ un $\mathcal{A}_K([1-\varepsilon,1[)$ module différentiel ayant la propriété de Robba et dont l'exposant $\mathfrak{E}xp_0(\mathcal{M}_{[1-\varepsilon,1[})$ a des différences qui ont la propriété (**NL**). On definit la filtration de la monodromie par $W_0(\mathcal{M}_{[1-\varepsilon,1[}) := \mathcal{M}_{[1-\varepsilon,1[}$ et $W_{k+1}(\mathcal{M}_{[1-\varepsilon,1[}) := \bar{0}_0(W_k(\mathcal{M}_{[1-\varepsilon,1[}))$ pour $\varepsilon > 0$ assez petit. En vertu du théorème 3.2-1 on a bien une filtration décroissante qui est stationnaire. En vertu de la proposition précédente les quotients sucessives admettent des solutions bornées dans le disque générique $D(t_\rho, \rho^-)$ pour $\rho$ assez proche de 1.

*Définition* 6.1-7. Soit $\mathcal{M}$ un $\mathcal{R}_K(1)$ module différentiel ayant la propriété de Robba et dont l'exposant $\mathfrak{E}xp_0(\mathcal{M}_{[1-\varepsilon,1[})$ a des différences qui ont la propriété (**NL**), on définit sa *filtration de la monodromie* $W_k(\mathcal{M})$ comme la filtration provenant de la filtration de la monodromie de sa restriction à $[1-\varepsilon,1[$ pour $\varepsilon > 0$ assez petit.

Corollaire 6.1-8. *Supposons que le corps $K$ est localement compact et soit $\mathcal{M}$ un $\mathcal{R}_K(1)$-module libre de rang $m$ à connexion soluble en 1 de plus grande pente $\mathrm{pt}(\mathcal{M}) > 0$ et $\gamma$ un nombre réel tel que $0 \leq \gamma < \mathrm{pt}(\mathcal{M})$. Supposons que l'adhérence $\bar{0}_\gamma(\mathcal{M}_{[1-\varepsilon,1[})$ de zéro dans $\mathcal{M}_{[1-\varepsilon,1[}$ pour la topologie quotient $\mathcal{T}_{\gamma,Q}$ n'est pas égale à $\mathcal{M}_{[1-\varepsilon,1[}$ pour tout $\varepsilon > 0$ assez petit, alors il existe une suite exacte de $\mathcal{R}_K(1)$-modules libres à connexion solubles en 1:*

$$0 \to \mathcal{M}_{\mathrm{inj}}^\gamma \to \mathcal{M} \to \mathcal{M}_{\mathrm{sol}}^\gamma \to 0,$$

*telle que $\mathcal{M}_{\mathrm{inj}}^\gamma$ est de plus grande pente égale à $\mathrm{pt}(\mathcal{M})$ et $\mathcal{M}_{\mathrm{sol}}^\gamma$ est non nul et de plus grande pente $\leq \gamma$.*

*Démonstration.* Soit $\mathcal{M}_{[1-\varepsilon,1[}$ la restriction de $\mathcal{M}$ à la couronne $C([1-\varepsilon, 1[)$ pour $\varepsilon$ assez petit. En vertu du théorème 6.1-1, l'adhérence $\bar{0}_\gamma(\mathcal{M}_{[1-\varepsilon,1[})$ de zéro dans $\mathcal{M}_{[1-\varepsilon,1[}$ pour la topologie quotient $\mathcal{T}_{\gamma,Q}$ n'est pas nulle sous les hypothèses du corollaire. En vertu du théorème 3.2-1, $\bar{0}_\gamma(\mathcal{M}_{[1-\varepsilon,1[})$ est un $\mathcal{A}([1-\varepsilon, 1[)$-module libre à connexion de rang $m_\varepsilon \leq m$ et le séparé associé $\mathcal{M}_{[1-\varepsilon,1[}/\bar{0}_\gamma(\mathcal{M}_{[1-\varepsilon,1[})$ est un $\mathcal{A}_K([1-\varepsilon, 1[)$-module libre à connexion de rang $m - m_\varepsilon \leq m$. Ils sont tous les deux automatiquement solubles en 1.

Les morphismes de restriction $\mathcal{M}_{[1-\varepsilon,1[} \to \mathcal{M}_{[1-\varepsilon',1[}, 0 < \varepsilon' < \varepsilon$ sont continus pour la topologie quotient $\mathcal{T}_{\gamma,Q}$ et induisent des morphismes

$$\bar{0}_\gamma(\mathcal{M}_{[1-\varepsilon,1[}) \to \bar{0}_\gamma(\mathcal{M}_{[1-\varepsilon',1[}).$$

Les images des morphismes précédents engendrent des modules à connexion nécessairement libres de rang $m_\varepsilon$ par platitude. Donc le rang $m_\varepsilon$ est une fonction croissante bornée par le rang $m$ et stationnaire de valeur, disons $m_{\varepsilon_0}$, par hypothèse on a $m_{\varepsilon_0} < m$.

Cela veut dire que pour $\varepsilon < \varepsilon_0$ la restriction du séparé associé

$$\mathcal{M}_{[1-\varepsilon,1[}/\bar{0}_\gamma(\mathcal{M}_{[1-\varepsilon,1[})$$

à la couronne $C([1 - \varepsilon', 1[)$ est isomorphe à $\mathcal{M}_{[1-\varepsilon',1[}/\bar{0}_\gamma(\mathcal{M}_{[1-\varepsilon',1[})$ pour tout $\varepsilon', 0 < \varepsilon' \leq \varepsilon$.

Pour $\varepsilon$ assez petit on peut alors supposer que la fonction rayon de convergence de $\mathcal{M}_{[1-\varepsilon,1[}/\bar{0}_\gamma(\mathcal{M}_{[1-\varepsilon,1[})$ est égale à $\rho^\delta$ pour tout $\rho \in [1-\varepsilon, 1[$ où $\delta - 1$ est sa pente en vertu du théorème 4.2-1. Comme il est séparé par construction, en vertu du corollaire 6.1-3, $\delta - 1 \leq \gamma$.

La suite exacte pour $\varepsilon$ assez petit

$$0 \to \bar{0}_\gamma(\mathcal{M}_{[1-\varepsilon,1[}) \to \mathcal{M}_{[1-\varepsilon,1[} \to \mathcal{M}_{[1-\varepsilon,1[}/\bar{0}_\gamma(\mathcal{M}_{[1-\varepsilon,1[}) \to 0$$

fournit par changement de base la suite exacte du corollaire:

$$0 \to \mathcal{M}_{\mathrm{inj}}^\gamma \to \mathcal{M} \to \mathcal{M}_{\mathrm{sol}}^\gamma \to 0.$$

Remarquons pour terminer la démonstration du corollaire que, pour $\varepsilon$ assez petit, parce que $R(\mathcal{M}_{[1-\varepsilon,1[}, \rho) < R(\mathcal{M}_{[1-\varepsilon,1[}/\bar{0}_\gamma(\mathcal{M}_{[1-\varepsilon,1[}, \rho)$ la fonction $R((\overline{0}_\gamma(\mathcal{M}_{[1-\varepsilon,1[}),\rho)$ est égale à la fonction $R(\mathcal{M}_{[1-\varepsilon,1[},\rho)$. En particulier la plus grande pente de $\mathcal{M}_{\mathrm{inj}}^\gamma$ est égale à pt($\mathcal{M}$).

*Définition* 6.1-9. Soit $\mathcal{M}$ un $\mathcal{R}_K(1)$-module libre de rang $m$ à connexion soluble en 1, nous dirons que $\mathcal{M}$ a toutes ses *pentes strictement plus grandes* qu'un nombre réel $\gamma, 0 \leq \gamma$ si pour $\varepsilon > 0$ assez petit toutes les solutions locales au point générique $t_\rho$ ont un rayon de convergence strictement plus petit que $\rho^{\gamma+1}$ pour tout $\rho \in [1-\varepsilon,1[$.

LEMME 6.1-10. *Soit un nombre réel* $\gamma$ *et* $\mathcal{M}$ *un* $\mathcal{A}_K([1-\varepsilon,1[)$-*module différentiel de rang* $m$ *dont le séparé associé pour la topologie* $\mathcal{T}_{\gamma,Q}$ *est nul, alors* $\mathcal{M}$ *n'admet aucune solution non triviale dans le disque générique* $D(t_\rho,\rho^{\gamma+1})$ *pour tout* $\rho \in [1-\varepsilon,1[$ *pour* $\varepsilon > 0$ *assez petit.*

*Démonstration.* Soit $G(x)$ la matrice de la connexion dans une base de $\mathcal{M}$ fournissant une suite exacte:

$$0 \to (\mathcal{A}_K([r-\varepsilon,r[)[\frac{d}{dx}])^m \xrightarrow{u:=\frac{d}{dx}-G(x)} (\mathcal{A}_K([r-\varepsilon,r[)[\frac{d}{dx}])^m \to \mathcal{M} \to 0.$$

L'adhérence de l'image $\mathcal{I}$ de $u$ dans $(\mathcal{A}_K([r-\varepsilon,r[)[\frac{d}{dx}])^m$ pour la topologie $\mathcal{T}_\gamma$ est égale à $(\mathcal{A}_K([r-\varepsilon,r[)[\frac{d}{dx}])^m$ tout entier. Puisque la topologie $\mathcal{T}_\gamma$ est métrisable, cela entraîne qu'il existe une suite $Q_n$ de matrices carrées d'ordre $m$ à coefficients dans $\mathcal{A}_K([r-\varepsilon,r[)[\frac{d}{dx}]$ telle que la suite $Q_n(\frac{d}{dx}-G)$ tend vers la matrice unité $I_m$ d'ordre $m$. Par définition les normes de la topologie $\mathcal{T}_\gamma$ de $\mathcal{A}([r-\varepsilon,r[)[\frac{d}{dx}]$ sont les normes d'opérateurs linéaires de l'espace de Banach des fonctions bornées $\mathcal{W}_{t_\rho}(\rho^{\gamma+1})$ dans le disque générique centré en $t_\rho$ et de rayon $\rho^{\gamma+1}$ pour $\rho \in [r-\varepsilon,r[$. Si $g$ est une solution bornée dans le disque générique centré en $t_\rho$ et de rayon $\rho^{\gamma+1}$ pour tout $\rho \in [r-\varepsilon,r[$ du système différentiel:

$$(\frac{d}{dx}-G(x))(g) = 0$$

alors nécessairement

$$g = \lim_{n\to\infty} Q_n(\frac{d}{dx}-G(x))(g) = 0.$$

Le module différentiel $\mathcal{M}$ n'admet aucune solution $g$ bornée non triviale dans le disque générique centré en $t_\rho$ et de rayon $\rho^{\gamma+1}$ pour tout $\rho \in [r-\varepsilon,r[$. Comme l'injectivité dans l'espace des fonctions bornées $\mathcal{W}_{t_\rho}(\rho^{\gamma+1})$ est équivalente à l'injectivité dans l'espace des fonctions analytiques $\mathcal{A}_{t_\rho}(\rho^{\gamma+1})$ ([R$_1$, 3.5]), les solutions locales au point générique $t_\rho$ du module différentiel $\mathcal{M}$ ont toutes un rayon de convergence strictement plus petit que $\rho^{\gamma+1}$ pour tout $\rho \in [r-\varepsilon,r[$. D'où le lemme.

COROLLAIRE 6.1-11. *Supposons que le corps* $K$ *est localement compact, soit* $\mathcal{M}$ *un* $\mathcal{R}_K(1)$-*module libre de rang* $m > 0$ *à connexion soluble en* 1 *de plus* *grande pente* pt($\mathcal{M}$) > 0 *et un nombre réel* $\gamma, 0 \leq \gamma <$ pt($\mathcal{M}$), *il existe une filtration décroissante finie de* $\mathcal{M}$, $\mathcal{M} = \mathcal{M}_0 \supset \mathcal{M}_1 \supset \cdots \supset \mathcal{M}_{i_0}$ *par des sous-$\mathcal{R}_K(1)$-modules libres de rang* $m_i > 0$ *à connexion soluble en 1 telle que* $\mathcal{M}_{i_0}$ *est de pentes strictement plus grandes que* $\gamma$ *et que les plus grandes pentes des quotients successifs* $\mathcal{M}_i/\mathcal{M}_{i+1}, i = 0 \cdots i_0 - 1$ *sont inférieures ou égales à* $\gamma$.

*Démonstration.* Nous raisonnons par récurrence sur le rang de $\mathcal{M}$. Si le rang de $\mathcal{M}$ est égal à un l'indice $i_0 = 0$ convient. Si le rang de $\mathcal{M}$ est > 1 et que l'adhérence $\bar{0}_\gamma(\mathcal{M}_{[1-\varepsilon,1[})$ de zéro dans $\mathcal{M}_{[1-\varepsilon,1[}$ pour la topologie quotient $\mathcal{T}_{\gamma,Q}$ est égale à $\mathcal{M}_{[1-\varepsilon,1[}$ pour tout $\varepsilon > 0$ assez petit, l'indice $i_0 = 0$ convient en vertu du lemme 6.1-10. Si l'adhérence $\bar{0}_\gamma(\mathcal{M}_{[1-\varepsilon,1[})$ de zéro dans $\mathcal{M}_{[1-\varepsilon,1[}$ pour la topologie quotient $\mathcal{T}_{\gamma,Q}$ n'est pas égale à $\mathcal{M}_{[1-\varepsilon,1[}$ pour tout $\varepsilon > 0$ assez petit, alors en vertu du corollaire 6.1-8 il existe une suite exacte:

$$0 \to \mathcal{M}_{\mathrm{inj}}^\gamma \to \mathcal{M} \to \mathcal{M}_{\mathrm{sol}}^\gamma \to 0,$$

telle que le rang de $\mathcal{M}_{\mathrm{inj}}^\gamma$ est strictement plus petit que le rang de $\mathcal{M}$ et que la plus grande pente de $\mathcal{M}_{\mathrm{sol}}^\gamma$ est inférieure ou égale à $\gamma$. L'hypothèse de récurrence permet de conclure.

COROLLAIRE 6.1-12. *Supposons que le corps* $K$ *est localement compact; soit* $\mathcal{M}$ *un* $\mathcal{R}_K(1)$-*module libre de rang* $m > 0$ *à connexion soluble en 1 de plus grande pente* pt($\mathcal{M}$) > 0 *et un nombre réel* $\gamma, 0 \leq \gamma <$ pt($\mathcal{M}$); *il existe une décomposition*:

$$0 \to \mathcal{M}_{>\gamma} \to \mathcal{M} \to \mathcal{M}^{\leq \gamma} \to 0$$

*par des* $\mathcal{R}_K(1)$-*modules libres à connexion solubles en 1 telle que la plus grande pente de* $\mathcal{M}^{\leq \gamma}$ *est inférieure ou égale à* $\gamma$ *et toutes les pentes de* $\mathcal{M}_{>\gamma}$ *sont strictement supérieures à* $\gamma$.

*Démonstration.* Si on définit $\mathcal{M}_{>\gamma}$ comme le sous-module $\mathcal{M}_{i_0}$ du corollaire 6.1-11 et $\mathcal{M}^{\leq \gamma}$ comme le module quotient dans la suite exacte

$$0 \to \mathcal{M}_{>\gamma} \to \mathcal{M} \to \mathcal{M}^{\leq \gamma} \to 0$$

on a toutes les propriétés du corollaire 6.1-12, car $\mathcal{M}^{\leq \gamma}$ apparaît comme extensions successives des modules quotients $\mathcal{M}_i/\mathcal{M}_{i+1}$.

COROLLAIRE 6.1-13. *Sous les hypothèses précédentes, la dimension*

$$\dim_K \mathrm{Hom}_{\mathcal{A}_K([1-\varepsilon,1[)[\frac{d}{dx}]}(\mathcal{M}_{[1-\varepsilon,1[}, \mathcal{A}_{t_\rho}(\rho^{\gamma+1}))$$

*pour* $\varepsilon$ *assez petit et pour* $\rho \in [1-\varepsilon,1[$, *est égale au rang* $m^{\leq \gamma}$ *du module* $\mathcal{M}^{\leq \gamma}$ *et est indépendante de* $\rho \in [1-\varepsilon,1[$.

*Démonstration.* En effet pour $\varepsilon > 0$ assez petit et tout $\rho \in [1-\varepsilon,1[$, $\mathcal{M}_{>\gamma[1-\varepsilon,1[}$ n'admet aucune solution non triviale dans l'espace $\mathcal{A}_{t_\rho}(\rho^{\gamma+1})$ en vertu du lemme 6.1-10 et donc

$$\dim_K \operatorname{Hom}_{\mathcal{A}_K([1-\varepsilon,1[)[\frac{d}{dx}]}(\mathcal{M}_{[1-\varepsilon,1[}, \mathcal{A}_{t_\rho}(\rho^{\gamma+1}))$$
$$= \dim_K \operatorname{Hom}_{\mathcal{A}_K([1-\varepsilon,1[)[\frac{d}{dx}]}(\mathcal{M}_{[1-\varepsilon,1[}^{\leq\gamma}, \mathcal{A}_{t_\rho}(\rho^{\gamma+1})) = m^{\leq\gamma}.$$

Ce corollaire appliqué au cas $\gamma = 0$ démontre la conjecture 3.1.1 de [C-M₁].

Appliquons le corollaire 6.1-12 pour $\gamma = 0$; on trouve la décomposition fondamentale qui implique déjà, compte tenu des résultats des articles [C-M₁] [C-M₂], le théorème de l'indice 7.4-1:

THÉORÈME 6.1-14. *Supposons que le corps $K$ est localement compact, soit $\mathcal{M}$ un $\mathcal{R}_K(1)$-module libre de rang $m > 0$ à connexion soluble en 1 de plus grande pente $\operatorname{pt}(\mathcal{M}) > 0$, il existe une décomposition*:

$$0 \to \mathcal{M}_{>0} \to \mathcal{M} \to \mathcal{M}^{\leq 0} \to 0$$

*par des $\mathcal{R}_K(1)$-modules libres à connexion solubles en 1 telle que le module $\mathcal{M}^{\leq 0}$ a la propriété de Robba au bord et le module $\mathcal{M}_{>0}$ est injectif dans l'espace $\mathcal{A}_{t_\rho}(\rho)$ pour tout $\rho$ assez proche de 1.*

*Définition* 6.1-15. Supposons que le corps $K$ est localement compact, soit $\mathcal{M}$ un $\mathcal{R}_K(1)$-module libre de rang $m > 0$ à connexion soluble en 1 de plus grande pente $\operatorname{pt}(\mathcal{M}) > 0$. Nous définissons $\mathcal{M}^{\leq 0}$ comme la *partie modérée* de $\mathcal{M}$ et $\mathcal{M}_{>0}$ comme la *partie de pentes strictement positives de $\mathcal{M}$*.

*Exemple* 6.1-16. C'est l'exemple fondamental qui nous a suggéré la forme du théorème 6.1-14. Considérons l'opérateur différentiel d'ordre deux où $\partial = x\frac{d}{dx}$:

$$P(x,\partial) := x^2(\partial + \frac{1}{3})(\partial + \frac{5}{3}) + \frac{8\pi}{27}\partial.$$

Cet opérateur admet zéro et l'infini comme seules singularités. L'infini est une singularité régulière alors que zéro est une singularité irrégulière. En vertu de [Me₄, 4.1.1] pour $p \neq 3$, il est muni d'une structure de Frobenius sur l'anneau des séries de Laurent $(K[x, x^{-1}])^\dagger$ à coefficient dans $K$ qui convergent dans un domaine $1 - \varepsilon \leq |x| \leq 1 + \varepsilon$ pour $\varepsilon > 0$ non précisé. C'est un cas particulier des équations $M_{f,n,m}$, étudiées dans [Me₄] qui permettent de montrer le théorème de finitude des nombres de Betti $p$-adiques. On l'obtient pour $n = 1$, $f(x) = x^2, m = 3$. En particulier il est soluble en 1. Cet opérateur admet la fonction $g(x) = 1 + \sum_{k \geq 1} a_k x^{2k}$ où

$$a_k = -\frac{27(2k + 1/3)(2k + 5/3)}{16\pi k} a_{k-1}, k \geq 1$$

comme solution. Si le corps $K$ contient $\pi$ et si la caractéristique résiduelle est $> 5$, cette fonction est un élément de $\mathcal{A}_K(1)$ qui *n'est pas borné au bord*.

Elle n'appartient pas au corps des éléments analytiques au bord. L'opérateur admet une factorisation sur le corps des fractions de l'anneau $\mathcal{R}_K(1)$ et non sur le corps des éléments analytiques au bord comme l'espérait Robba [R$_5$].

On peut montrer en utilisant le théorème de transfert [C$_3$] que l'opérateur $P$ n'a pas la propriété de Robba dans la classe résiduelle de zéro lorsque la caractéristique résiduelle est différente de 2. Donc, dans ce cas là, le rang de sa partie modérée est égale à un parce que la restriction au disque générique $D(t_\rho, \rho^-)$ de la solution $g$ est une solution non trivial de rayon de convergence minoré par $\rho$ pour tout $\rho$ assez proche de 1. On a alors une factorisation dans $\mathcal{A}_K([1-\varepsilon, 1[)[\frac{d}{dx}]$ pour $\varepsilon > 0$ assez petit:

$$hP = fP_1P_2$$

où $h(x)$ et $f(x)$ sont des fonctions, $P_1$ et $P_2$ sont des opérateurs différentiels d'ordre un tels que $P_2 = g(x)\frac{d}{dx} - g'(x)$. Comme la pente de $P_1$ est strictement plus grande que un l'opérateur identique est adhérent à l'idéal engendré par $fP_1$ pour la topologie $\mathcal{T}_0$. L'opérateur $P_2$ est adhérent à l'idéal engendré par $P$ pour la topologie $\mathcal{T}_0$. L'adhérence de l'idéal $P$ qui contient l'idéal $(P, g(x)\frac{d}{dx} - g'(x))$ est égal à ce dernier parce que le rang de la partie modérée de $P$ est égale à un.

*Remarque* 6.1-17. Comme on vient de le voir on ne peut espérer avoir un théorème de décomposition ayant les propriétés du théorème 6.1-14 sur le corps des éléments analytiques au bord. Ceci montre que le théorème de décomposition de Dwork-Robba [D-R$_1$] est indépendant du théorème de décomposition 6.1-14 aussi bien du point de vue des résultats que du point de vue des démonstrations. C'est bien entendu là un point de structure essentiel. Le théorème de décomposition 6.1-14 s'apparente au premier théorème de décomposition de Robba [R$_1$] en *famille* pour $\rho$ variable dans l'intervalle $[1-\varepsilon, 1[$. Cette confusion entre les deux types de décomposition a pendant longtemps retardé notre compréhension du théorème 6.1-14.

*Remarque* 6.1-18. L'exemple précédent montre que les coefficients des opérateurs différentiels de la décomposition de Robba [R$_1$] d'un polynôme différentiel soluble au point générique $t_1$ ne sont pas superadmissibles en général et sont au plus admissibles dans la terminologie de Dwork-Robba [D-R$_1$].

*Remarque* 6.1-19. Dans le paragraphe 6 nous avons supposé le corps de base localement compact pour pouvoir utiliser le lemme 6.1-2. Cependant la démonstration du théorème 6.1-1 suggère d'élargir la notion de suites extraites et considérer la notion de c-suites extraites dans le cas d'un corps maximalement complet.

*Définition* 6.1-20.  Une c-suite extraite d'une suite $f_n$ d'un espace vectoriel localement convexe sur un corps valué $K$ est une suite $g_n$ pour laquelle il existe des scalaires $\lambda_{i,n}$ tels que

1) les $\lambda_{i,n}$ sont presque tous nuls pour $n$ fixé,

2) $g_n = \sum_{i \geq n} \lambda_{i,n} f_i$,

3) $|\lambda_{i,n}| \leq 1$ et

4) $\sum_{i \geq n} \lambda_{i,n} = 1$.

On peut alors montrer [C-M$_4$] la généralisation du lemme 6.1-2 sur un corps maximalement complet $K$: une suite $f_n$ bornée de $\mathcal{A}_K(I)$ pour la topologie d'espace de type $\mathcal{F}$, pour un intervalle ouvert $I$, et *minorée* en norme par une constante strictement positive en un point de $I$ admet des c-suites extraites convergentes vers une limite *non nulle*. Le lecteur pourra vérifier que ce résultat permet de supposer le corps de base maximalement complet dans les paragraphes 6 et 7.

6.2.  *Décomposition par rapport à toutes les pentes, polygone de Newton p-adique.* Nous allons montrer que les modules différentiels $\mathcal{M}_{>\gamma}$ constituent une filtration décroissante dont les sauts définissent les pentes $p$-adiques de $\mathcal{M}$.

THÉORÈME 6.2-1.  *Supposons que le corps $K$ est localement compact; soit $\mathcal{M}$ un $\mathcal{R}_K(1)$-module libre de rang $m > 0$ à connexion soluble en 1 de plus grande pente* $\mathrm{pt}(\mathcal{M}) > 0$, *pour tous nombres réels* $\gamma', \gamma$ *tels que* $0 \leq \gamma' < \gamma < \mathrm{pt}(\mathcal{M})$, $\mathcal{M}_{>\gamma}$ *est un sous-module à connexion de* $\mathcal{M}_{>\gamma'}$.

*Démonstration.*  Le choix d'une base $\mathcal{G}$ de $\mathcal{M}_{[1-\varepsilon_0,1[}$ définit une suite exacte (∗)

$$0 \to (\mathcal{A}_K([1-\varepsilon,1[)[\frac{d}{dx}])^m \xrightarrow{u := \frac{d}{dx} - G(x)} (\mathcal{A}_K([1-\varepsilon,1[)[\frac{d}{dx}])^m \to \mathcal{M}_{[1-\varepsilon,1[} \to 0$$

pour $\varepsilon > 0$ assez petit. Nous affirmons que pour tout nombre réel $\gamma$ *strictement* plus grand que la plus grande pente $\mathrm{pt}(\mathcal{M})$ la topologie quotient $\mathcal{T}_{\gamma,Q}$ sur $\mathcal{M}_{[1-\varepsilon,1[}$ est *séparée* pour $\varepsilon > 0$ assez petit. En effet si on munit $\mathcal{M}_{[1-\varepsilon,1[}$ de la topologie naturelle d'espace $\mathcal{F}$ par la base $\mathcal{G}$ en vertu de la majoration de la proposition 5.3-1 le dernier morphisme de la suite (∗) est *continu.* Donc son noyau est fermé et la topologie quotient est séparée.

Soit alors deux nombres réels $\gamma', \gamma$ tels que $0 \leq \gamma' < \gamma < \mathrm{pt}(\mathcal{M})$; par construction la plus grande pente $\mathrm{pt}(\mathcal{M}^{\leq \gamma'})$ est inférieure ou égale à $\gamma'$. D'après ce qui précède la topologie quotient $\mathcal{T}_{\gamma,Q}$ est séparée sur $\mathcal{M}^{\leq \gamma'}_{[1-\varepsilon,1[}$ pour $\varepsilon > 0$ assez petit.

Considérons le morphisme canonique

$$\mathcal{M}_{[1-\varepsilon,1[} \to \mathcal{M}^{\leq \gamma'}_{[1-\varepsilon,1[}$$

qui est continu pour la topologie $\mathcal{T}_{\gamma,Q}$. Nous venons de montrer que l'image, par ce morphisme, de l'adhérence de zéro dans $\mathcal{M}_{[1-\varepsilon,1[}$ pour la topologie $\mathcal{T}_{\gamma,Q}$ est *nulle* pour $\varepsilon > 0$ assez petit. En particulier l'image du module $\mathcal{M}_{>\gamma[1-\varepsilon,1[}$ par ce morphisme est nulle. Donc ce dernier est contenu dans le noyau $\mathcal{M}_{>\gamma'[1-\varepsilon,1[}$ pour $\varepsilon > 0$ assez petit. D'où le théorème 6.2-1.

*Remarque* 6.2-2. Ce qui précède est une réciproque de 6.1-3; si pt($\mathcal{M}$) $\neq 0$ la topologie $\mathcal{T}_{\gamma,Q}$ sur $\mathcal{M}_{[1-\varepsilon,1[}$ pour tout $\varepsilon > 0$ assez petit est séparée si et seulement si pt($\mathcal{M}$) $\leq \gamma$.

COROLLAIRE 6.2-3.  *Supposons que le corps $K$ est localement compact; soit $\mathcal{M}$ un $\mathcal{R}_K(1)$-module libre de rang $m > 0$ à connexion soluble en 1 de plus grande pente pt($\mathcal{M}$) > 0, les $\mathcal{R}_K(1)$-modules $\mathcal{M}_{>\gamma}$ libres de rang fini à connexion soluble en 1 constituent une filtration finie décroissante de $\mathcal{M}$ indexée par les nombres réels $\gamma$ tels que $0 \leq \gamma < $ pt($\mathcal{M}$).*

La fonction rang de $\mathcal{M}_{>\gamma}$ bornée supérieurement par le rang $m$ est une fonction décroissante avec $\gamma$, elle n'a donc qu'un nombre fini de sauts.

Nous définissons le $\mathcal{R}_K(1)$-module à connexion gradué $\mathrm{Gr}_\gamma(\mathcal{M})$ par la suite exacte pour un nombre réel $\eta$ assez petit:

$$0 \to \mathcal{M}_{>\gamma} \to \mathcal{M}_{>\gamma-\eta} \to \mathrm{Gr}_\gamma(\mathcal{M}) \to 0$$

ou de façon équivalente par la suite exacte:

$$0 \to \mathrm{Gr}_\gamma(\mathcal{M}) \to \mathcal{M}^{\leq \gamma} \to \mathcal{M}^{\leq \gamma-\eta} \to 0.$$

*Définition* 6.2-4. Nous dirons que $\gamma$ est une *pente* de $\mathcal{M}$ si le module $\mathrm{Gr}_\gamma(\mathcal{M})$ est non nul.

Lorsque le corps de base $K$ est localement compact, le $\mathcal{R}_K(1)$-module $\mathrm{Gr}_\gamma(\mathcal{M})$ est libre de rang fini $m_\gamma$ comme noyau d'un morphisme de $\mathcal{R}_K(1)$-modules libres de rang fini.

*Définition* 6.2-5. Nous dirons qu'un $\mathcal{R}_K(1)$-*module libre de rang $m > 0$ à* connexion soluble en 1 est purement de pente $\gamma$ si toutes ses solutions locales au point générique $t_\rho$ admettent un même rayon de convergence $\rho^{(\gamma+1)}$ pour tout $\rho \in [1-\varepsilon,1[$ pour un $\varepsilon$ assez petit.

COROLLAIRE 6.2-6.  *Supposons que le corps $K$ est localement compact, soit $\mathcal{M}$ un $\mathcal{R}_K(1)$-module libre de rang $m > 0$ à connexion soluble en 1 de plus grande pente pt($\mathcal{M}$) > 0 alors les pentes de $\mathcal{M}$ sont des nombres rationnels > 0 et les modules $\mathrm{Gr}_\gamma(\mathcal{M})$ sont purement de pente $\gamma$.*

*Démonstration.* En vertu de la suite exacte

$$0 \to \mathcal{M}_{>\gamma} \to \mathcal{M}_{>\gamma-\eta} \to \mathrm{Gr}_\gamma(\mathcal{M}) \to 0$$

les solutions locales au point générique $t_\rho$ du gradué $\mathrm{Gr}_\gamma(\mathcal{M})$ admettent toutes un rayon de convergence strictement majoré par $\rho^{\gamma+1-\eta}$ pour tout $\eta > 0$ assez petit et pour tout $\rho \in [1-\varepsilon, 1[$ pour un $\varepsilon$ assez petit. En vertu de la suite exacte

$$0 \to \mathrm{Gr}_\gamma(\mathcal{M}) \to \mathcal{M}^{\leq \gamma} \to \mathcal{M}^{\leq \gamma-\eta} \to 0$$

toutes les solutions locales au point générique $t_\rho$ du gradué $\mathrm{Gr}_\gamma(\mathcal{M})$ admettent un rayon de convergence minoré par $\rho^{\gamma+1}$ pour tout $\rho \in [1-\varepsilon, 1[$ pour un $\varepsilon$ assez petit. Le module $\mathrm{Gr}_\gamma(\mathcal{M})$ est purement de pente de $\gamma$ qui est aussi sa plus grande pente. En vertu du théorème 4.2-1 ce sont des nombres rationnels $> 0$.

Supposons que le corps $K$ est localement compact, soit $\mathcal{M}$ un $\mathcal{R}_K(1)$-module libre de rang $m > 0$ à connexion soluble en 1 de plus grande pente $\mathrm{pt}(\mathcal{M}) > 0$ alors il admet un nombre fini de pentes $0 < \mathrm{pt}(\mathcal{M})_1 \leq \cdots \leq \mathrm{pt}(\mathcal{M})_N = \mathrm{pt}(\mathcal{M})$.

*Définition* 6.2-7. Sous les hypothèses précédentes on appelle *multiplicité d'une pente strictement positive* $\mathrm{pt}(\mathcal{M})_i$ de $\mathcal{M}$ le rang $m_{\mathrm{pt}(\mathcal{M})_i}$ du $\mathcal{R}_K(1)$-module libre $\mathrm{Gr}_{\mathrm{pt}(\mathcal{M})_i}(\mathcal{M})$. On appelle *multiplicité de la pente nulle* le rang $m_0$ de la partie modérée $\mathcal{M}^{\leq 0}$.

*Définition* 6.2-8. Si le corps $K$ est localement compact, on appelle *polygone de Newton p-adique d'un $\mathcal{R}_K(1)$-module $\mathcal{M}$* libre de rang $m > 0$ à connexion soluble en 1 le polygone du plan réel à coordonnées, *a priori* rationnelles, construit à partir des multiplicités et des pentes:

$$\mathrm{Newton}(\mathcal{M}, p) := (m_0, 0; m_1, \mathrm{pt}(\mathcal{M})_1; \cdots; m_N, \mathrm{pt}(\mathcal{M})_N)$$

en partant de l'origine.

Le théorème de l'indice local 8.3-1 montrera que les sommets du polygone de Newton sont des entiers rationnels ce qui constituera l'analogue *p*-adique du théorème de Hasse-Arf.

6.3. *Propriétés fonctorielles de la filtration* $\mathcal{M}_{>\gamma}$. Nous prolongeons la définition de la filtration $\mathcal{M}_{>\gamma}$ en posant $\mathcal{M}_{>\gamma} = 0$ pour tout réel $\gamma \geq \mathrm{pt}(\mathcal{M})$ et $\mathcal{M}_{>\gamma} := \mathcal{M}$ pour tout réel $\gamma < 0$.

6.3.1. *Exactitude.* Nous allons montrer que les modules $\mathcal{M}_{>\gamma}$, $\mathcal{M}^{\leq\gamma}$, et $\mathrm{Gr}_\gamma(\mathcal{M})$ dépendent fonctoriellement de $\mathcal{M}$ pour tout $\gamma$.

PROPOSITION 6.3-1. *Si le corps de base est localement compact, pour tout nombre réel $\gamma \geq 0$ les modules $\mathcal{M}_{>\gamma}$, $\mathcal{M}^{\leq\gamma}$ et $\mathrm{Gr}_\gamma(\mathcal{M})$ dépendent fonctoriellement du $\mathcal{R}_K(1)$-module $\mathcal{M}$ libre de rang fini à connexion soluble en 1* *et définissent des foncteurs exacts de la catégorie abélienne* MLS($\mathcal{R}_K(1)$) *dans elle même.*

*Démonstration.* Il suffit de traiter le cas du module $\mathcal{M}_{>\gamma}$. Si $\mathcal{M} \to \mathcal{N}$ est un morphisme de la catégorie MLS($\mathcal{A}_K([1-\varepsilon,1[)$) pour $\varepsilon > 0$ assez petit il est automatiquement continu pour la topologie quotient $\mathcal{T}_{\gamma,Q}$. Donc il se restreint en une application

$$\bar{O}_\gamma(\mathcal{M}) \to \bar{O}_\gamma(\mathcal{N})$$

qui en vertu du théorème 3.2-1 est un morphisme de la catégorie MLS($\mathcal{A}_K([1-\varepsilon,1[)$) et donc automatiquement continu pour la topologie quotient $\mathcal{T}_{\gamma',Q}$ pour tout nombre réel $\gamma' \geq 0$.

Si $\gamma < \inf(\mathrm{pt}(\mathcal{M}), \mathrm{pt}(\mathcal{N}))$, en vertu des corollaires 6.1-11 et 6.1-12 on obtient ainsi de proche en proche un morphisme

$$\mathcal{M}_{>\gamma} \to \mathcal{N}_{>\gamma}$$

de la catégorie MLS($\mathcal{R}_K(1)$).

Si $\gamma \geq \mathrm{pt}(\mathcal{M})$, $\mathcal{M}_{>\gamma} = 0$ et le morphisme est trivial.

Si $\gamma \geq \mathrm{pt}(\mathcal{N})$, alors l'image de $\mathcal{M}_{>\gamma}$ est automatiquement nulle parce ses solutions au point générique $t_\rho$ ont un rayon de convergence strictement plus petit que $\rho^{\gamma+1}$ alors que les solutions de $\mathcal{N}$ au point générique $t_\rho$ ont un rayon de convergence au moins égal à $\rho^{(\mathrm{pt}(\mathcal{N})+1)}$ pour tout $\rho \in [1-\varepsilon,1[$ pour $\varepsilon > 0$ assez petit.

Dans tous les cas un morphisme $\mathcal{M} \to \mathcal{N}$ de la catégorie MLS($\mathcal{R}_K(1)$) induit un morphisme $\mathcal{M}_{>\gamma} \to \mathcal{N}_{>\gamma}$, pour tout réel $\gamma$, définissant un foncteur.

Soit une suite exacte de la catégorie MLS($\mathcal{R}_K(1)$):

$$0 \to \mathcal{P} \to \mathcal{M} \to \mathcal{N} \to 0,$$

on obtient donc une suite exacte de complexes:

$$
\begin{array}{ccccccccc}
0 & \to & \mathcal{P}_{>\gamma} & \to \mathcal{P} & \to \mathcal{P}^{\leq\gamma} & \to 0 \\
 & & \downarrow & \downarrow & \downarrow \\
0 & \to & \mathcal{M}_{>\gamma} & \to \mathcal{M} & \to \mathcal{M}^{\leq\gamma} & \to 0 \\
 & & \downarrow & \downarrow & \downarrow \\
0 & \to & \mathcal{N}_{>\gamma} & \to \mathcal{N} & \to \mathcal{N}^{\leq\gamma} & \to 0
\end{array}
$$

qui montre que la cohomologie de la première colonne formée de modules différentiels de pentes strictement supérieure à $\gamma$ est isomorphe à la cohomologie de la dernière colonne formée de modules différentiels de pentes inférieure à $\gamma$. Ces cohomologies sont donc nulles, d'où la proposition. Le lecteur remarquera que l'exactitude du foncteur solution au point générique ([R$_1$, 4.23]) est essentielle ici.

COROLLAIRE 6.3-2. *Un morphisme de la catégorie* MLS($\mathcal{R}_K(1)$) *est automatiquement* strict *pour la filtration* $(-)_{>\gamma}$.

COROLLAIRE 6.3-3.  *Pour tout réel $\gamma, 0 \leq \gamma < \mathrm{pt}(\mathcal{M})$, $\mathcal{M}_{>\gamma}$, resp. $\mathcal{M}^{\leq\gamma}$, est le plus grand sous-module, resp. quotient, dont toutes les pentes sont $> \gamma$, resp. $\leq \gamma$.*

*Démonstration.* En effet si $\mathcal{N}$ est un sous $\mathcal{R}_K(1)$-module libre de rang fini à connexion de $\mathcal{M}$ dont toutes les pentes sont strictement plus grandes que $\gamma$ en vertu du corollaire 6.1-12 l'inclusion $\mathcal{N}_{>\gamma} \subset \mathcal{N}$ est une égalité et en vertu de la fonctorialité précédente, $\mathcal{N}_{>\gamma}$ est un sous module de $\mathcal{M}_{>\gamma}$. Un raisonnement similaire vaut pour le module quotient.

6..3.2. *Images inverses.* Soit $\varphi_q$ la ramification d'ordre $q$ de la couronne $C([r,1[)$ dans la couronne $C([r^q,1[)$ induisant un morphisme $\varphi_q^* : \mathcal{R}_K(1) \rightarrow \mathcal{R}_K(1)$. Si $\mathcal{M}$ est un $\mathcal{A}_K([r,1[)$-module libre de rang $m > 0$ à connexion notons $\varphi_q^*(\mathcal{M})$ son image inverse par $\varphi_q^*$.

LEMME 6.3-4.  *On a les inégalités*:

$$R(\mathcal{M}, \rho^q)/\rho^{q-1} \leq R(\varphi_q^*(\mathcal{M}), \rho) \leq \rho$$

*pour tout $\rho \in [r,1[$. En particulier, si $\mathcal{M}$ est soluble en 1 son image inverse $\varphi_q^*(\mathcal{M})$ est aussi soluble en 1. Si l'indice de ramification $q$ est premier avec $p$, la première inégalité est une égalité; en particulier la plus grande pente $\mathrm{pt}(\varphi_q^*(\mathcal{M}))$ est égal à $\mathrm{pt}(\mathcal{M})q$.*

*Démonstration.* Soit $\mathcal{F}$ une base de $\mathcal{M}$ de matrice de connexion $F$ et $\mathcal{G}$ la base de $\varphi_q^*(\mathcal{M})$ image inverse de $\mathcal{F}$ de matrice de connexion $G$. Notons $Y$, resp. $X$, la solution, au voisinage du point générique $t_\rho$ ($r \leq \rho < 1$), resp. $t_\rho^q$ ($r^q \leq |t_\rho^q| = \rho^q < 1$), de l'équation différentielle:

$$Y' = F\,Y, \quad Y(t_\rho) = I, \quad \text{resp.} \quad X' = G\,X, \quad X(t_\rho^q) = I,$$

associée à cette base. Si $q$ est premier avec $p$, on a, pour $|x - t_\rho| < |t_\rho| = \rho$:

$$|x^q - t_\rho^q| = \left| \sum_{i=1}^{q} \binom{q}{i}(x - t_\rho)^i\, t_\rho^{q-i} \right| = |x - t_\rho|\,|t_\rho|^{q-1} = |x - t_\rho|\,\rho^{q-1}.$$

Si bien que les coefficients de la matrice $Y$ sont analytiques dans le disque $D(t_\rho, r)$ si et seulement si ceux de la matrice $X$ convergent dans le disque $D(t_\rho^q, r\,\rho^{q-1})$. Autrement dit, on a:

$$R(\mathcal{M}, \rho^q) = R(\varphi_q^*(\mathcal{M}), \rho)\,\rho^{q-1}.$$

Si $q$ est divisible par $p$, on a encore:

$$|x^q - t_\rho^q| \leq |x - t_\rho|\,\rho^{q-1}$$

d'où l'on déduit que:

$$R(\mathcal{M}, \rho^q) \leq R(\varphi_q^*(\mathcal{M}), \rho)\,\rho^{q-1}.$$

Remarquons que l'application image inverse

$$\varphi_q^{-1}: \qquad \mathcal{A}_K([1-\varepsilon,1[)[\frac{d}{dx}] \to \mathcal{A}_K([(1-\varepsilon)^{1/q},1[)[\frac{d}{dx}]$$

qui à $x$ associe $x^q$ et à $x\frac{d}{dx}$ associe $\frac{1}{q}x\frac{d}{dx}$ induit l'application image inverse:

$$\varphi_q^*: \mathcal{A}_K([(1-\varepsilon)^{1/q},1[)\otimes_{\mathcal{A}_K([1-\varepsilon,1[)}\mathcal{A}_K([1-\varepsilon,1[)[\frac{d}{dx}] \to \mathcal{A}_K([(1-\varepsilon)^{1/q},1[)[\frac{d}{dx}]$$

qui est un isomorphisme impliquant une décomposition algébrique:

$$\bigoplus_{1\le k<q} x^k \otimes \mathcal{A}_K([1-\varepsilon,1[)[\frac{d}{dx}] = \mathcal{A}_K([(1-\varepsilon)^{1/q},1[)[\frac{d}{dx}].$$

C'est là une somme directe topologique pour toute topologie $\mathcal{T}_\gamma$.

PROPOSITION 6.3-5.  *Supposons que le corps $K$ est localement compact, soient $\mathcal{M}$ un $\mathcal{R}_K(1)$-module libre de rang $m > 0$ à connexion soluble en 1 de pente $\mathrm{pt}(\mathcal{M}) > 0$, $q$ un entier et $\gamma$ un réel $0 \le \gamma < \mathrm{pt}(\mathcal{M})$. On a alors des injections de $\mathcal{R}_K(1)$-modules libres à connexion*

$$\varphi_q^*(\mathcal{M})_{>\gamma q} \to \varphi_q^*(\mathcal{M}_{>\gamma}).$$

LEMME 6.3-6.  *Pour tout réel $r, 0 < r \le 1$, et tout opérateur $P$ de $\mathcal{A}_K([r^q,1[)[\frac{d}{dx}]$ on a pour tout $\rho \in [r,1[$*

$$|P|_{\gamma,\rho^q} \le |\varphi_q^{-1}(P)|_{q\gamma,\rho}.$$

*Démonstration.* Notons $\mathcal{W}_{t_\rho}(r)$ l'anneau des fonctions analytiques bornées dans le disque $D(t_\rho, r)$ muni de la norme de la convergence uniforme:

$$\Big|\sum_{n=0,\infty} a_n\,(x-t_\rho)^n\Big|_r = \sup(|a_n|\,r^n).$$

La norme $|\ |_{\gamma,\rho}$ est alors la norme d'opérateur sur l'espace de Banach $\mathcal{W}_{t_{\rho^{\gamma+1}}}(\rho)$. On constate que, pour $|x-t_\rho| = \rho^\beta$, on a, sauf pour un nombre fini de valeurs de $\rho$:

$$\begin{aligned}
|x^q - t_\rho^q|_{\rho^\beta} &= \Big|\sum_{i=1}^q \binom{q}{i}(x-t_\rho)^i\,t_\rho^{q-i}\Big| \\
&= \max_{p^\alpha|q}\Big(|q\,p^{-\alpha}|\,\rho^{\beta p^\alpha}\,\rho^{q-p^\alpha}\Big) \\
&= \rho^{\delta(q,\rho,\beta)}
\end{aligned}$$

où l'on a posé:

$$\delta(q,\rho,\beta) = \min_{p^\alpha|q}\left(p^\alpha(\beta-1) + q + \frac{\mathrm{Log}(|q\,p^{-\alpha}|)}{\mathrm{Log}(\rho)}\right).$$

Sauf pour un nombre fini de valeurs de $\rho$, l'application $\varphi_q$ est donc une surjection du cercle $C_\zeta := \{x \in \mathbb{C}_p; |x - \zeta t_\rho| = \rho^\beta\}$ avec $\zeta^q = 1$, sur le cercle $C := \{x \in \mathbb{C}_p; |x^q - t_\rho^q| = \rho^{\delta(q,\rho,\beta)}}$. En fait si $q$ est premier avec $p$, l'application $\varphi_q$ est une bijection entre chacun des cercles $C_\zeta$ et $C$ et si $q$ est une puissance de $p$ et $\rho$ est assez proche de 1, l'application $\varphi_q$ est un revêtement de degré $q$. Pour une fonction $g = \sum_{n=0,\infty} a_n (x - t_\rho^q)^n$ de $\mathcal{W}_{t_{\rho^{\delta(q,\rho,\beta)}}}(\rho^q)$, on a, par définition:

$$|\varphi_q^{-1}(g)|_{\rho^\beta} = \max_{|x - t_\rho| = \rho^\beta}(|g(x^q)|)$$

$$= \max_{|y - t_\rho^q| = \rho^{\delta(p,\rho,\beta)}}(|g(y)|) = |g|_{\rho^{\delta(p,\rho,\beta)}}$$

car les maximums sont atteints si $x$ n'appartient pas à un nombre fini de classes résiduelles des cercles considérés.

Pour $\gamma$ défini par: $\delta(q,\rho,\beta) = q(\gamma+1)$, c'est-à-dire pour:

$$\beta = \beta(q,\rho,\gamma+1) = \max_{p^\alpha | q}\left(1 + p^{-\alpha}q\gamma - p^{-\alpha}\frac{\mathrm{Log}(|q\,p^{-\alpha}|)}{\mathrm{Log}(\rho)}\right),$$

on trouve

$$|\varphi_q^{-1}(g)|_{\rho^\beta} = |g|_{\rho^{q(\gamma+1)}}.$$

En particulier, $\varphi_q^{-1}$ est une injection de $\mathcal{W}_{t_{\rho^{\gamma+1}}}(\rho^q)$ dans $\mathcal{W}_{t_{\rho^\beta}}(\rho)$. On constate que $\beta(q,\rho,\gamma+1) \leq 1 + q\gamma$ avec égalité si $q$ est premier à $p$. Mais, si $q$ est une puissance de $p$ et $\rho$ suffisamment proche de 1, on a $\beta(q,\rho,\gamma+1) = \gamma+1$. Maintenant, par définition de $\varphi_q^{-1}(P)$, on a:

$$\varphi_q^{-1}(P).\varphi_q^{-1}(g) = \varphi_q^{-1}(P.g).$$

On obtient:

$$|\varphi_q^{-1}(P)|_{\beta,\rho} = \sup_{f \in \mathcal{W}_{t_{\rho^\beta}}(\rho)}\frac{|\varphi_q^{-1}(P).f|_{\rho^\beta}}{|f|_{\rho^\beta}} \geq \sup_{g \in \mathcal{W}_{t_{\rho^{\gamma+1}}}(\rho^q)}\frac{|\varphi_q^{-1}(P.g)|_{\rho^\beta}}{|\varphi_q^{-1}(g)|_{\rho^\beta}}$$

$$= \sup_{g \in \mathcal{W}_{t_{\rho^{\gamma+1}}}(\rho^q)}\frac{|P.g|_{\rho^{q(\gamma+1)}}}{|g|_{\rho^{q(\gamma+1)}}} = |P|_{\gamma,\rho^q}.$$

Il suffit alors de remarquer que la norme $|\;|_{\gamma,\rho}$ est croissante avec $\gamma$ pour conclure.

*Démonstration de la proposition* 6.3-5. Soit $\mathcal{M}$ un $\mathcal{R}_K(1)$-module libre de rang $m > 0$ à connexion soluble en 1; nous allons montrer que, pour $\varepsilon > 0$ assez petit,

$$\bar{0}_{q\gamma}(\varphi_q^*(\mathcal{M}_{[1-\varepsilon,1[}, \mathcal{T}_{q\gamma,Q}))$$

est contenu dans

$$\varphi_q^*(\bar{0}_\gamma(\mathcal{M}_{[1-\varepsilon,1[},\mathcal{T}_{\gamma,Q}))$$

comme sous-modules de $\varphi_q^*(\mathcal{M}_{[1-\varepsilon,1[})$. Le choix d'une base $\mathcal{G}$ de $\mathcal{M}$ définit une suite exacte

$$0 \to (\mathcal{A}_K([1-\varepsilon,1[)[\frac{d}{dx}])^m \overset{u:=\frac{d}{dx}-G(x)}{\longrightarrow} (\mathcal{A}_K([1-\varepsilon,1[)[\frac{d}{dx}])^m \to \mathcal{M}_{[1-\varepsilon,1[} \to 0$$

pour tout $\varepsilon > 0$ assez petit où $G(x)$ est la matrice de la connexion dans la base $\mathcal{G}$. Cette suite exacte induit une suite exacte:

$$\begin{aligned} 0 \to (\mathcal{A}_K([(1-\varepsilon)^{1/q},1[)[\frac{d}{dx}])^m &\overset{u_q}{\longrightarrow} (\mathcal{A}_K([(1-\varepsilon)^{1/q},1[)[\frac{d}{dx}])^m \\ &\longrightarrow \varphi^*(\mathcal{M}_{[1-\varepsilon,1[}) \to 0, \end{aligned}$$

où $u_q := \frac{d}{dx} - qx^{q-1}G(x^q)$. Tout vecteur $P$ de $(\mathcal{A}_K([(1-\varepsilon)^{1/q},1[)[\frac{d}{dx}])^m$ se décompose comme une somme $\sum_{0\le k<q} x^k \otimes P_k$. S'il est adhérent à l'image de $u_q$ pour la topologie $\mathcal{T}_{q\gamma,Q}$ les vecteurs $P_k$ sont adhérents à l'image de $u$ pour la topologie $\mathcal{T}_\gamma$ en vertu de la décomposition topologique précédente et de la majoration du lemme 6.3-6 montrant ainsi que la classe de $P$ est dans $\varphi_q^*(\bar{0}_\gamma(\mathcal{M}_{[1-\varepsilon,1[},\mathcal{T}_{\gamma,Q}))$.

Nous pouvons supposer que $q\gamma$ est strictement plus petit que la plus grande pente de $\varphi_q^*(\mathcal{M})$. Les $\mathcal{A}([(1-\varepsilon)^{1/q},1[)[\frac{d}{dx}]$-modules $\bar{0}_{q\gamma}(\varphi_q^*(\mathcal{M}_{[1-\varepsilon,1[},\mathcal{T}_{q\gamma,Q}))$, $\varphi_q^*(\bar{0}_\gamma(\mathcal{M}_{[1-\varepsilon,1[},\mathcal{T}_{\gamma,Q}))$ sont de type fini en vertu du théorème 3.2-1 et l'inclusion

$$\bar{0}_{q\gamma}(\varphi_q^*(\mathcal{M}_{[1-\varepsilon,1[},\mathcal{T}_{q\gamma,Q})) \to \varphi_q^*(\bar{0}_\gamma(\mathcal{M}_{[1-\varepsilon,1[},\mathcal{T}_{\gamma,Q}))$$

est automatiquement continue pour la topologie $\mathcal{T}_{q\gamma,Q}$ ce qui entraîne par récurrence que pour tout $i$

$$\varphi_q^*(\mathcal{M}_{[1-\varepsilon,1[})_i \subset \varphi_q^*(\mathcal{M}_{[1-\varepsilon,1[i})$$

où $\varphi_q^*(\mathcal{M}_{[1-\varepsilon,1[})_i$ désigne la filtration de $\varphi_q^*(\mathcal{M}_{[1-\varepsilon,1[})$ définie par la topologie $\mathcal{T}_{q\gamma,Q}$ alors que $\mathcal{M}_{[1-\varepsilon,1[i}$ désigne la filtration de $\mathcal{M}_{[1-\varepsilon,1[}$ définie par la topologie $\mathcal{T}_{\gamma,Q}$. On obtient ainsi l'inclusion de la proposition pour $\varepsilon > 0$ assez petit

$$\varphi_q^*(\mathcal{M})_{>q\gamma} \subset \varphi_q^*(\mathcal{M}_{>\gamma}).$$

*Remarque* 6.3-7. On peut montrer, ce n'est pas difficile, que si l'indice de ramification $q$ est premier à $p$, l'inégalité du lemme 6.3-6 est une égalité et que les injections de la proposition 6.3-5 sont des égalités, autrement dit le polygone de Newton de l'image inverse s'obtient par homothétie de rapport l'indice de ramification.

*Remarque* 6.3-8. Si $q$ est une puissance de $p$, on a, pour $|x - t_\rho|$ suffisamment proche de $\rho$:

$$|x^q - t_\rho^q| = |x - t_\rho|^q,$$

c'est-à-dire, pour $R(\varphi_q^*(\mathcal{M}), \rho)$ suffisamment proche de $\rho$:

$$R(\mathcal{M}, \rho^q) = R(\varphi_q^*(\mathcal{M}), \rho)^q$$

relation déjà donnée dans le théorème 4.1-3 et précisée dans le calcul de la proposition 5.2-3.

*Exemple* 6.3-9. D'après ce qui précède la ramification d'ordre une puissance de $p$ ne change pas la plus grande pente. En particulier un module différentiel sur $\mathcal{R}_K(1)$ soluble en 1 dont le dénominateur de la plus grande pente contient des puissances de $p$ *ne se décompose pas* en modules de rang un par ramification. C'est là une différence profonde avec la théorie formelle des équations différentielles sur un corps de caractéristique nulle qui complique considérablement la théorie $p$-adique. Une telle situation se produit même pour les équations provenant de la géométrie. En effet considérons l'opérateur différentiel

$$9x^3 \frac{d^2}{dx^2} + 9x^2 \frac{d}{dx} - x + \frac{\pi^2}{3}.$$

C'est encore un cas particulier des modules exponentiels $M_{f,n,m}$. On l'obtient pour $n = 1$, $f(x) = x$, $m = 3$. En vertu de [Me$_4$, 4.1.1], pour $p \neq 3$, il est muni d'une structure de Frobenius sur l'anneau des séries de Laurent $(K[x, x^{-1}])^\dagger$ à coefficients dans $K$ qui convergent dans un domaine $1 - \varepsilon \leq |x| \leq 1 + \varepsilon$ pour $\varepsilon > 0$ non précisé. On peut montrer qu'il est de pentes strictement positives. Ceci entraîne qu'il est purement de pente $1/2$. Pour $p = 2$, il est irréductible sur l'anneau $\mathcal{R}_K(1)$ et indécomposable par ramification de Frobenius alors qu'il se décompose sur le corps des séries formelles $K((y))$, $x = y^2$. En particulier cet opérateur fournit un contre exemple au principe du transfert pour les singularités irrégulières.

6.3.3. *Compatibilité avec une structure de Frobenius.* Soit **F** le morphisme de Frobenius $\mathcal{R}_K(1) \to \mathcal{R}_K(1)$ induit par la ramification $\varphi_q^*$ d'ordre $q$ pour une puissance $q$ de $p$.

*Définition* 6.3-10. On dit qu'un $\mathcal{R}_K(1)$-module libre de rang $m$ à connexion $\mathcal{M}$ est muni d'une structure de Frobenius s'il existe un isomorphisme

$$\mathbf{F} : \varphi_q^*(\mathcal{M}) \simeq \mathcal{M}$$

de $\mathcal{R}_K(1)$-modules libres de rang $m$ à connexion.

PROPOSITION 6.3-11. *Un $\mathcal{R}_K(1)$-module $\mathcal{M}$ libre de rang $m$ à connexion muni d'une structure de Frobenius est automatiquement soluble en* 1.

*Démonstration.* Soit $X(x, t_\rho) = \sum_{k=0}^{\infty} G_k(t_\rho) \frac{(x-t_\rho)^k}{k!}$ la solution fonda-mentale de $\mathcal{M}$ au voisinage du point générique $t_\rho$ pour $\rho$ assez proche de 1. La solution fondamentale de $\varphi_q^*(\mathcal{M})$ au voisinage de $t_\rho$ est égale à $X(x^q, t_\rho^q)$ qui converge pour $|x^q - t_\rho^q| \leq R(\mathcal{M}, \rho^q)$. Mais $|x^q - t_\rho^q| \leq \max(|x - t_\rho|^q,$ $|x - t_\rho|/q)$. Donc la série $X(x^q, t_\rho^q)$ converge dans le disque centré en $t_\rho$ et de rayon $\min(R(\mathcal{M}, \rho^q)^{1/q}, qR(\mathcal{M}, \rho^q))$. On a alors la majoration pour $\rho$ assez proche de 1:

$$R(\mathcal{M}, \rho) = R(\varphi_q^*\mathcal{M}, \rho) \geq \min(R(\mathcal{M}, \rho^q)^{1/q}, qR(\mathcal{M}, \rho^q))$$

qui donne par passage à la limite quand $\rho$ tend vers un par valeurs inférieures l'inégalité:

$$R(\mathcal{M}, 1^-) \geq \min(R(\mathcal{M}, 1^-)^{1/q}, qR(\mathcal{M}, 1^-)) = R(\mathcal{M}, 1^-)^{1/q}.$$

On obtient l'inégalité $R(\mathcal{M}, 1^-) \geq 1$ et donc l'égalité $R(\mathcal{M}, 1^-) = 1$.

Notons $\mathrm{MLS}(\mathcal{R}_K(1), \mathbf{F})$ la catégorie des $\mathcal{R}_K(1)$-modules libres de rang fini à connexion munis d'une structure de Frobenius; c'est une sous-catégorie pleine de la catégorie $\mathrm{MLS}(\mathcal{R}_K(1))$.

PROPOSITION 6.3-12.  *Supposons que le corps de base $K$ est localement compact; la catégorie $\mathrm{MLS}(\mathcal{R}_K(1), \mathbf{F})$ est stable par les foncteurs exacts $\mathcal{M} \to \mathcal{M}_{>\gamma}$.*

*Démonstration.* Par fonctorialité et en tenant compte de la proposition 6.3-5 on trouve une injection,

$$\mathcal{M}_{>\gamma} \overset{\mathbf{F}^{-1}}{\simeq} \varphi_q^*(\mathcal{M})_{>\gamma} \hookrightarrow \varphi_q^*(\mathcal{M}_{>\gamma/q}).$$

Il suffit de voir que l'image de $\mathcal{M}_{>\gamma}$ est contenue dans $\varphi_q^*(\mathcal{M}_{>\gamma})$ pour qu'elle lui soit isomorphe comme sous-module à connexion de *même rang*. Mais en vertu de 6.3-6 la fonction rayon de convergence ne change pas par ramification d'ordre une puissance de $p$. Le corollaire 6.3-3 montre que $\varphi_q^*((\mathcal{M}_{>\gamma/q})^{\leq \gamma})$ est un quotient de $\varphi_q^*(\mathcal{M}_{>\gamma/q})^{\leq \gamma}$ et donc $\varphi_q^*(\mathcal{M})_{>\gamma}$ est un sous-module de $\varphi_q^*(\mathcal{M}_{>\gamma})$.

*Remarque* 6.3-13.  Nous n'avons considéré que la ramification de Frobenius sur l'anneau $\mathcal{R}_K(1)$ pour simplifier. Mais on peut naturellement considérer des morphismes de Frobenius plus généraux, par exemple obtenus en composant un morphisme de Frobenius du corps de base avec une ramification de Frobenius.

## 7. Le théorème de l'indice

**7.1.** *Descente de la décomposition d'un module différentiel ayant la propriété de Robba au corps de base.* Soit $K$ un sous-corps complet du corps $\mathbb{C}_p$ des complexes $p$-adiques. Si $K$ est à valuation discrète, pour tout intervalle $I$, les extensions $\mathcal{A}_K(I) \to \mathcal{A}_{\mathbb{C}_p}(I)$ et $\mathcal{R}_K(r) \to \mathcal{R}_{\mathbb{C}_p}(r)$ sont plates et donc fidèlement plates.

Soit $\mathcal{M}$ un $\mathcal{R}_{\mathbb{C}_p}(r)$-module libre de rang $m$ à connexion ayant la propriété de Robba: la fonction $R(\mathcal{M}, \rho)$ est égale à $\rho$ pour $\rho$ dans un intervalle $[r - \varepsilon, r[, \varepsilon > 0$. Autrement dit, la plus grande pente pt$(\mathcal{M})$ est nulle. Nous avons défini dans [C-M$_2$, §5] *l'exposant* $\mathfrak{e}xp_0^r(\mathcal{M})$ de $\mathcal{M}$ comme une classe de $(\mathbb{Z}_p/\mathbb{Z})^m$ modulo la relation d'équivalence $\overset{\mathfrak{e}}{\sim}$.

Soit $\mathcal{M}$ un $\mathcal{R}_K(r)$-module libre de rang $m$ à connexion de pente nulle alors son étendu $\mathcal{R}_{\mathbb{C}_p(r)} \otimes_{\mathcal{R}_K(r)} \mathcal{M}$ a la propriété de Robba. Nous définissons son exposant $\mathfrak{e}xp_0^r(\mathcal{M})$ comme $\mathfrak{e}xp_0^r(\mathcal{R}_{\mathbb{C}_p(r)} \otimes_{\mathcal{R}_K(r)} \mathcal{M})$. C'est un élément de $\mathfrak{e}_m := (\mathbb{Z}_p/\mathbb{Z})^m / \overset{\mathfrak{e}}{\sim}$ bien défini. Nous notons Rob$(\mathcal{R}_K(r))$ la catégorie des $\mathcal{R}_K(r)$-modules libres de rang fini à connexion ayant la propriété de Robba. En vertu de la proposition 3.2-5, si $K$ est à valuation discrète, c'est une catégorie abélienne.

Pour un élément de l'ensemble $\mathfrak{e}_m$ nous pouvons parler de la propriété (**NL**) pour ses différences ([C-M$_2$, 4.4]). La classe d'équivalence d'un élément de $\mathfrak{e}_m$ qui a la propriété (**NL**) pour ses différences contient un unique élément représenté par un ensemble $\{\alpha_1, \ldots, \alpha_m\}$ d'éléments de $\mathbb{Z}_p/\mathbb{Z}$ ([C-M$_2$, 4.4-9]). Nous notons $\mathfrak{e}_m^{\textbf{NL}}$ le sous-ensemble des éléments de $\mathfrak{e}^m$ qui ont la propriété (**NL**) pour leur différences.

*Définition 7.1-1.* Nous dirons qu'un *module de la catégorie* Rob$(\mathcal{R}_K(r))$ a la propriété (**NL**$^*$) si son exposant a des différences qui ont la propriété (**NL**).

Nous notons Rob$(\mathcal{R}_K(r), \textbf{NL}^*)$ la catégorie des $\mathcal{R}_K(r)$-modules libres de rang fini à connexion qui ont la propriété (**NL**$^*$). En vertu de [C-M$_2$, 5.4-6] si $K$ est à valuation discrète la catégorie Rob$(\mathcal{R}_K(r), \textbf{NL}^*)$ est une sous-catégorie abélienne de la catégorie Rob$(\mathcal{R}_K(r))$. Si un $\mathcal{R}_K(r)$-module $\mathcal{M}$ libre de rang $m$ à connexion de plus grande pente nulle a la propriété (**NL**$^*$) nous pouvons parler de ses exposants $\{\alpha_1, \ldots, \alpha_m\}$ qui sont des éléments bien définis de $\mathbb{Z}_p/\mathbb{Z}$.

De la même façon on peut parler de la propriété de Robba pour un module différentiel sur l'anneau $\mathcal{A}_K(]r, R[)$. Soit $\mathcal{M}$ un $\mathcal{A}_K(]r, R[)$-module différentiel ayant la propriété de Robba sur une couronne $C(]r, R[)$ dont l'exposant, qui est alors bien défini, a la propriété (**NL**$^*$), en vertu du théorème fondamental ([C-M$_2$]) il se décompose en modules de rang un sur l'anneau $\mathcal{A}_{\mathbb{C}_p}(]r, R[)$.

Nous allons montrer, comme conséquence du théorème Tate-Ax [A], que cette décomposition a lieu déjà sur le corps de base $K$.

Pour un nombre $\alpha$ de $\mathbb{Z}_p$ notons $\Psi_K^\alpha(]r, R[)$ l'espace des sommes finies

$$\sum_{\alpha, k} c_{\alpha, k}(x) x^\alpha \frac{(\mathrm{Log}(x))^k}{k!}$$

où les $c_{\alpha, k}(x)$ sont des fonctions analytiques sur la couronne $C(]r, R[)$ à coefficients dans le corps $K$. C'est un $\mathcal{A}_K(]r, R[)[\partial]$-module libre sur l'anneau $\mathcal{A}_K(]r, R[)$, où $\partial := x \frac{d}{dx}$. Nous choisissons un représentant $\alpha$ pour chaque classe $\bar\alpha$ de $\mathbb{Z}_p/\mathbb{Z}$ et soit $\Psi_K(]r, R[)$ la somme directe des modules $\Psi_K^\alpha(]r, R[)$ pour $\bar\alpha$ variable dans $\mathbb{Z}_p/\mathbb{Z}$. Le module $\Psi_K^\alpha(]r, R[)$ ne dépend, à isomorphisme près, que de la classe de $\alpha$ dans $\mathbb{Z}_p/\mathbb{Z}$.

THÉORÈME 7.1-2. *Soit $\mathcal{M}$ un module différentiel libre de rang $m$ sur l'anneau $\mathcal{A}_K(]r, R[)$ ayant la propriété de Robba sur la couronne $C(]r, R[)$ et dont l'exposant a la propriété $(\mathbf{NL}^*)$. L'espace des solutions multiformes de détermination finie*

$$\mathrm{Hom}_{\mathcal{A}_K(]r, R[)[\partial]}(\mathcal{M}, \Psi_K^\beta(]r, R[))$$

*est nul si $\beta$ n'est pas un exposant de $\mathcal{M}$, l'espace*

$$\mathrm{Ext}^1_{\mathcal{A}_K(]r, R[)[\partial]}(\mathcal{M}, \Psi_K^\beta(]r, R[))$$

*est nul si $\beta$ est un exposant de $\mathcal{M}$ et*

$$\dim_K \mathrm{Hom}_{\mathcal{A}_K(]r, R[)[\partial]}(\mathcal{M}, \Psi_K(]r, R[)) = m.$$

*Démonstration.* Considérons le groupe de Galois $G_K^{\mathrm{cont}}$ des $K$-automorphismes continus de $\mathbb{C}_p$. Un élément $g$ de $G_K^{\mathrm{cont}}$ est une isométrie de $\mathbb{C}_p$; cf. [C₁, 1.6.1]. Ceci permet de définir une action de $G_K^{\mathrm{cont}}$ sur $\mathcal{A}_{\mathbb{C}_p}(]r, R[)$ par $g(\sum_{k \in \mathbb{Z}} a_k x^k) := \sum_{k \in \mathbb{Z}} g(a_k) x^k$ et une action de $G_K^{\mathrm{cont}}$ sur $\Psi_{\mathbb{C}_p}(]r, R[)$ par $g(\sum_{\alpha, k} c_{\alpha, k}(x) x^\alpha \frac{(\mathrm{Log}(x))^k}{k!}) := \sum_{\alpha, k} g(c_{\alpha, k}(x)) x^\alpha \frac{(\mathrm{Log}(x))^k}{k!}$. Le point est qu'en vertu du théorème de Tate-Ax [A] les points fixes de $G_K^{\mathrm{cont}}$ sont les fonctions de $\mathcal{A}_K(]r, R[)$ et de $\Psi_K(]r, R[)$.

Soit $x_0$ un point de la couronne $C(]r, R[)$, $\det(x^\alpha)_{x_0}$ et $\det(\mathrm{Log}(x))_{x_0}$ des déterminations locales en $x_0$ de $x^\alpha$ et $\mathrm{Log}(x)$, c'est-à-dire des solutions locales des équations $(\partial - \alpha)u = 0$ et $\partial u = 1$. On a alors un morphisme de restriction:

$$\mathrm{rest}_{x_0} : \Psi_{\mathbb{C}_p}(]r, R[) \to \mathcal{A}_{x_0}$$

qui à $\sum_{\alpha, k} c_{\alpha, k}(x) x^\alpha \frac{(\mathrm{Log}(x))^k}{k!}$ associe $\sum_{\alpha, k} c_{\alpha, k}(x) \det(x^\alpha)_{x_0} \frac{(\det(\mathrm{Log}(x))_{x_0})^k}{k!}$. Par construction cette restriction commute à l'action des opérateurs différentiels; d'où un morphisme

$$(*) \quad \mathrm{Hom}_{\mathcal{A}_K(]r, R[)[\partial]}(\mathcal{M}, \oplus_{\bar\alpha \in \mathfrak{E}xp_0^{]r, R[}(\mathcal{M})} \Psi_{\mathbb{C}_p}^\alpha(]r, R[)) \to \mathrm{Hom}_{\mathcal{A}_{x_0}[\partial]}(\mathcal{M}_{x_0}, \mathcal{A}_{x_0}),$$

où $\mathcal{A}_{x_0}$ est le $\mathbb{C}_p$-espace des fonctions analytiques au voisinage de $x_0$ et $\mathcal{M}_{x_0} := \mathcal{A}_{x_0} \otimes_{\mathcal{A}_K(]r,R[)} \mathcal{M}$. Si $m$ est égal à un le morphisme précédent est injectif, c'est donc un isomorphisme, parce que la restriction d'une fonction $c(x)x^\alpha$ est nulle si et seulement si la fonction $c(x)$ est nulle. On en déduit, en vertu du théorème de décomposition [C-M$_2$] sur le corps $\mathbb{C}_p$, que le morphisme (*) est un isomorphisme, donnant par la même occasion un sens au prolongement analytique multiforme des solutions locales de $\mathcal{M}$.

Supposons que le point $x_0$ appartient au corps de base $K$, alors l'espace $\mathcal{A}_{x_0}$ est muni d'une action de $G_K^{\text{cont}}$ par

$$g(\sum_{k\in\mathbb{N}} a_k(x-x_0)^k) := \sum_{k\in\mathbb{N}} g(a_k)(x-x_0)^k$$

dont les points fixes sont les séries à coefficients dans le corps $K$ en vertu du théorème de Tate-Ax. D'autre part le morphisme de restriction commute à l'action de $G_K^{\text{cont}}$ et donc le morphisme (*) commute à l'action de $G_K^{\text{cont}}$ puisque $\mathcal{M}$ est défini sur $K$ par hypothèse. D'autre part $\mathcal{M}$ admet une base fondamentale de solutions locales en $x_0$ à coefficients dans le corps $K$, donc il admet une base fondamentale de fonctions multiformes à coefficients dans le corps $K$, ce qui implique la décomposition dans le corps de base $K$ et le théorème 7.1-2 dans ce cas là.

Quitte à faire une extension finie on peut supposer qu'il existe un point $K$-rationnel dans la couronne $C(]r,R[)$. Mais la formation des espaces de solutions commute aux extensions finies. On en déduit la décomposition sur le corps de base $K$ et le théorème 7.1-2.

*Définition* 7.1-3. Nous dirons qu'un *module de la catégorie* $\text{Rob}(\mathcal{R}_K(r))$ a la propriété (**NL$^{**}$**) s'il a la propriété (**NL$^*$**) et si ses exposants ont la propriété (**NL**).

Nous notons $\text{Rob}(\mathcal{R}_K(r), \mathbf{NL}^{**})$ la sous-catégorie de la catégorie $\text{Rob}(\mathcal{R}_K(r), \mathbf{NL}^*)$ des modules dont les exposants ont la propriété (**NL**). Pour la démonstration du théorème de l'indice nous utilisons la conséquence suivante du théorème fondamental:

THÉORÈME 7.1-4. *Pour tout module* $\mathcal{M}$ *de la catégorie* $\text{Rob}(\mathcal{R}_K(r), \mathbf{NL}^{**})$ *l'indice dans l'espace* $\mathcal{R}_K(r)$ *est nul*:

$$\begin{aligned}\chi(\mathcal{M}, \mathcal{R}_K(r)) &:= \dim_K \text{Hom}_{\mathcal{R}_K(r)[\partial]}(\mathcal{M}, \mathcal{R}_K(r)) \\ &- \dim_K \text{Ext}^1_{\mathcal{R}_K(r)[\partial]}(\mathcal{M}, \mathcal{R}_K(r)) = 0.\end{aligned}$$

Remarquons que toute base d'un module différentiel $\mathcal{M}$ sur l'anneau $\mathcal{A}_K(I)$ définit une résolution de ce module de longueur un par des modules libres de type fini sur l'anneau $\mathcal{A}_K(I)[\partial]$. En particulier pour un $\mathcal{R}_K(r)$-module différentiel $\mathcal{M}$ les espaces $\text{Ext}^i_{\mathcal{R}_K(r)[\partial]}(\mathcal{M}, \mathcal{G})$ sont nuls pour $i \geq 2$ et tout $\mathcal{R}_K(r)[\partial]$-module à gauche $\mathcal{G}$.

7.2. *Exposants p-adiques des $\mathcal{R}_K(r)$-modules solubles.* Soit maintenant un module $\mathcal{M}$ de la catégorie MLS($\mathcal{R}_K(r)$); alors en vertu du théorème 6.1-14 si le corps $K$ est localement compact et si $r$ appartient au groupe des valeurs absolues de $K$, il est extension de sa partie modérée $\mathcal{M}^{\leq 0}$ par sa partie de pentes strictement positives $\mathcal{M}_{>0}$:

$$0 \to \mathcal{M}_{>0} \to \mathcal{M} \to \mathcal{M}^{\leq 0} \to 0.$$

La partie modérée $\mathcal{M}^{\leq 0}$ a la propriété de Robba; c'est un objet de la catégorie Rob($\mathcal{R}_K(r)$) et son exposant $\mathfrak{Exp}_0^r(\mathcal{M}^{\leq 0})$ est donc défini comme un élément de $\mathfrak{E}_m$ pour un entier $m$ inférieur ou égal au rang de $\mathcal{M}$.

*Définition* 7.2-1. Supposons le corps $K$ localement compact et soit $\mathcal{M}$ un $\mathcal{R}_K(r)$-module libre de rang fini à connexion soluble en $r$. On définit son exposant $\mathfrak{Exp}_0^r(\mathcal{M})$ comme *l'exposant* $\mathfrak{Exp}_0^r(\mathcal{M}^{\leq 0})$ *de sa partie modérée.*

C'est donc un élément de l'ensemble $\mathfrak{E}_m$ pour un entier $m \leq rg(\mathcal{M})$. On peut considérer alors la sous-catégorie MLS($\mathcal{R}_K(r), \mathbf{NL}^*$) de la catégorie MLS($\mathcal{R}_K(r)$), uniquement définie lorsque le corps de base est localement compact, des modules libres à connexion dont l'exposant a la propriété ($\mathbf{NL}^*$). De même nous pouvons considérer la sous-catégorie pleine

$$\text{MLS}(\mathcal{R}_K(r), \mathbf{NL}^{**})$$

de la catégorie MLS($\mathcal{R}_K(r), \mathbf{NL}^*$) des modules dont les exposants ont la propriété ($\mathbf{NL}$). La catégorie MLS($\mathcal{R}_K(r), \mathbf{NL}^{**}$) est la catégorie de base pour les propriétés de finitude dans la théorie des coefficients $p$-adiques [C-M$_4$].

PROPOSITION 7.2-2. *Si le corps $K$ est localement compact et si $r$ appartient au groupe des valeurs absolues de $K$, les catégories* MLS($\mathcal{R}_K(r), \mathbf{NL}^*$) *et* MLS($\mathcal{R}_K(r), \mathbf{NL}^{**}$) *sont des sous-catégories abéliennes de la catégorie* MLS($\mathcal{R}_K(r)$).

*Démonstration.* En effet le foncteur qui à un module de la catégorie MLS($\mathcal{R}_K(r)$) associe sa partie modérée est un foncteur exact de la catégorie MLS($\mathcal{R}_K(r)$) dans la catégorie Rob($\mathcal{R}_K(r)$) en vertu de 6.3-1. La proposition est conséquence du fait que les catégories Rob($\mathcal{R}_K(r), \mathbf{NL}^*$) et Rob($\mathcal{R}_K(r)$, $\mathbf{NL}^{**}$) sont abéliennes ([C-M$_2$, 5.4-6]).

En vertu du théorème de l'indice de Robba [R$_1$] pour les modules injectifs à coefficients dans le *corps des éléments analytiques* au bord l'indice $\chi(\mathcal{M}_{>0}, \mathcal{R}_K(r - \varepsilon))$ est nul pour $\varepsilon > 0$ assez petit.

Si $\mathcal{M}$ est un objet de la catégorie MLS($\mathcal{R}_K(r), \mathbf{NL}^{**}$) cela entraîne, en vertu du théorème 7.1-4, pour $\varepsilon \geq 0$ assez petit que l'indice $\chi(\mathcal{M}^{\leq 0}, \mathcal{R}_K(r-\varepsilon))$ est nul, donc que pour un $\varepsilon > 0$ assez petit, l'indice $\chi(\mathcal{M}, \mathcal{R}_K(r - \varepsilon))$ est nul. Mais, *attention*, on ne peut rien dire sur les indices $\chi(\mathcal{M}_{>0}, \mathcal{R}_K(r))$ et $\chi(\mathcal{M}, \mathcal{R}_K(r))$ à ce stade.

7.3. *Exposants p-adiques des $E_K^\dagger(r)$-modules.* Notons $E_K^\dagger(r)$ "le corps des
éléments analytiques au bord" à coefficients dans un corps complet $K$, c'est la
limite inductive des anneaux $H_K([r-\varepsilon, r[), \varepsilon > 0$, des éléments analytiques à
coefficients dans le corps $K$ dans la couronne $C([r-\varepsilon, r[)$. Soit $\mathcal{M}$ un espace
vectoriel de dimension finie sur $E_K^\dagger(r)$ à connexion. En vertu du théorème de
décomposition de Dwork-Robba [D-R$_1$], $\mathcal{M}$ est extension d'un $E_K^\dagger(r)$-espace
vectoriel $\mathcal{M}_{sol}$ soluble en $r$ par un espace vectoriel $\mathcal{M}_{inj}$ injectif en $r$:

$$0 \to \mathcal{M}_{inj} \to \mathcal{M} \to \mathcal{M}_{sol} \to 0.$$

Rappelons qu'on dit qu'un module différentiel est soluble en $r$ si la dimension
de l'espace de ses solutions analytiques dans le disque de générique centré en
$t_r$ de rayon $r$ est égale à son rang et qu'on dit qu'un module est injectif si cette
dimension est nulle.

Comme $E_K^\dagger(r)$ est un sous-anneau de $\mathcal{R}_K(r)$ on peut considérer
$\mathcal{R}_K(r) \otimes_{E_K^\dagger(r)} \mathcal{M}_{sol}$ qui est un $\mathcal{R}_K(r)$-module libre de rang fini à connexion
soluble en $r$ en vertu du théorème de continuité de la fonction rayon de con-
vergence ([C-D$_2$, 2.5]).

Si le corps $K$ est localement compact et si $r$ appartient au groupe des
valeurs absolues de $K$, on peut alors définir son exposant:

*Définition* 7.3-1. Soit $\mathcal{M}$ un espace vectoriel de dimension finie sur $E_K^\dagger(r)$
à connexion on définit son exposant $\mathfrak{exp}_0^r(\mathcal{M})$ comme l'exposant $\mathfrak{exp}_0^r(\mathcal{M}_{sol})$
défini auparavant.

On peut alors considérer la catégorie

$$\mathrm{MLC}(E_K^\dagger(r), \mathbf{NL}^*), \ \ \mathrm{resp.} \ \ \mathrm{MLC}(E_K^\dagger(r), \mathbf{NL}^{**}),$$

des espaces vectoriels de dimension finie à connexion dont l'exposant a la pro-
priété ($\mathbf{NL}^*$), resp. ($\mathbf{NL}^{**}$), qui est alors une catégorie abélienne.

7.4. *Le théorème de l'indice.* Soient $K$ un corps complet et $P(x, \frac{d}{dx})$ un
polynôme différentiel à coefficients dans le corps $K$, c'est-à-dire un opérateur
différentiel d'ordre fini à coefficients polynômiaux ou encore un élément de
l'algèbre de Weyl $K[x, \frac{d}{dx}]$. Pour tout nombre réel $r > 0$ le polynôme $P$ est
un élément de l'anneau $E_K^\dagger(r)[\frac{d}{dx}]$ et donc, si $K$ localement compact et si
$r$ appartient au groupe des valeurs absolues de $K$, l'exposant $\mathfrak{exp}_0^r(P)$ de $P$
est défini comme un élément de $\mathfrak{E}_m$ pour un entier $m \le rg(\mathcal{M})$. Si l'exposant
$\mathfrak{exp}_0^r(P)$ est un élément de $\mathfrak{E}_m^{\mathbf{NL}}$ on le note $\{\alpha_1, \ldots, \alpha_m\}$.

On arrive, finalement!, au théorème de l'indice conjecturé d'abord par
Bernard Dwork ([D$_5$, §4]) sous une forme très générale et précisé par la suite
par Philippe Robba (conjecture sur l'indice d'un opérateur différentiel, [R$_5$,
2.2]):

Théorème 7.4-1. *Soient un polynôme différentiel $P(x, \frac{d}{dx})$ à coefficients dans un sous-corps complet et localement compact $K$ de $\mathbb{C}_p$ et un nombre réel $r > 0$ appartenant au groupe des valeurs absolues de $K$. Si l'exposant $\mathfrak{E}xp_0^r(P)$ de $P$ a la propriété $(\mathbf{NL}^{**})$ l'opérateur $P$ est à indice dans les espaces:*

$$\mathcal{A}_{\mathbb{C}_p}(r), \ \mathcal{R}_{\mathbb{C}_p}(r), \ \mathcal{H}_{\mathbb{C}_p}^\dagger(r), \ \mathcal{A}_K(r), \ \mathcal{R}_K(r), \ \mathcal{H}_K^\dagger(r).$$

*De plus on a les égalités*

$$\chi(P, \mathcal{R}_{\mathbb{C}_p}(r)) = \chi(P, \mathcal{A}_{\mathbb{C}_p}(r)) + \chi(P, \mathcal{H}_{\mathbb{C}_p}^\dagger(r)) = 0,$$

$$\chi(P, \mathcal{A}_K(r)) = \chi(P, \mathcal{A}_{\mathbb{C}_p}(r)) = -\chi(P, \mathcal{H}_K^\dagger(r)) = -\chi(P, \mathcal{H}_{\mathbb{C}_p}^\dagger(r)).$$

*Démonstration.* Commençons par le cas des fonctions analytiques à co-efficients dans le corps $\mathbb{C}_p$. La démonstration est la même que [C-M$_1$, 5.1.1], sauf que maintenant nous disposons des propriétés de *finitude* des exposants qui nous faisaient encore défaut dans notre travail précédent ([C-M$_1$, intro.]). Il suffit de construire, en vertu du corollaire ([C-M$_1$, 4.5.3]) du théorème de dualité ([C-M$_1$, 4.1.1]) une suite de nombre réels $r_n$ tendant vers $r$ par valeurs inférieures tels que $\chi(P, \mathcal{R}_{\mathbb{C}_p}(r_n)) = 0$. Notons $\mathcal{M}$ le $E_K^\dagger(r)$-espace défini par $P$ et

$$0 \to \mathcal{M}_{\mathrm{inj}} \to \mathcal{M} \to \mathcal{M}_{\mathrm{sol}} \to 0$$

sa décomposition de Dwork-Robba [D-R$_1$]. Le module $\mathcal{M}_{\mathrm{inj}}$ est défini dans une couronne $C([r - \varepsilon, r[), \varepsilon > 0$, et est injectif en tout $r' \in [r - \varepsilon, r[$ par construction. En vertu du théorème de Robba [R$_1$] $\chi(\mathcal{M}_{\mathrm{inj}}, \mathcal{R}_{\mathbb{C}_p}(r'))$ est nul pour tout $r' \in [r - \varepsilon, r[$. Il suffit de montrer que $\chi(\mathcal{M}_{\mathrm{sol}}, \mathcal{R}_{\mathbb{C}_p}(r'))$ est nul pour tout $r' \in [r - \varepsilon, r[$ pour un $\varepsilon > 0$. Considérons alors la décomposition 6.1-14 de $\mathcal{M}_{\mathrm{sol}}$ en sa partie modérée et sa partie de pentes strictement positives:

$$0 \to \mathcal{M}_{\mathrm{sol}>0} \to \mathcal{R}_{\mathbb{C}_p}(r) \otimes_{E_K^\dagger(r)} \mathcal{M}_{\mathrm{sol}} \to \mathcal{M}_{\mathrm{sol}}^{\leq 0} \to 0.$$

Le module différentiel $\mathcal{M}_{\mathrm{sol}>0}$ est défini dans une couronne $C([r - \varepsilon, r[)$ pour un $\varepsilon > 0$ et est injectif en tout $r' \in [r - \varepsilon, r[$ par construction. En vertu du théorème de Robba [R$_1$] $\chi(\mathcal{M}_{\mathrm{sol}>0}, \mathcal{R}_{\mathbb{C}_p}(r'))$ est nul pour tout $r' \in [r - \varepsilon, r[$. Il suffit de montrer que $\chi(\mathcal{M}_{\mathrm{sol}}^{\leq 0}, \mathcal{R}_{\mathbb{C}_p}(r'))$ est nul pour tout $r' \in [r - \varepsilon, r[$ pour un $\varepsilon > 0$. Mais comme

$$\mathfrak{E}xp_0^r(P) := \mathfrak{E}xp_0^r(\mathcal{M}_{\mathrm{sol}}^{\leq 0}) = \{\alpha_1, \ldots, \alpha_m\}$$

est un élément de $\mathfrak{E}_m^{\mathbf{NL}}$ tels que les exposants $\alpha_i$ ont la propriétés $(\mathbf{NL})$ par hypothèse,

$$\chi(\mathcal{M}_{\mathrm{sol}}^{\leq 0}, \mathcal{R}_{\mathbb{C}_p}(r'))$$

est nul pour tout $r' \in [r - \varepsilon, r[$ en vertu du théorème 7.1-4. D'où le théorème de l'indice de $P$ à valeurs dans les fonctions analytiques à coefficients dans le corps $\mathbb{C}_p$.

Considérons le cas du corps $K$. Le raisonnement précédent s'applique pour tout corps complet contenant le corps de définition de $P$; en particulier l'indice $\chi(P, \mathcal{A}_K(r))$ est fini. Il suffit de montrer que l'indice $\chi(P, \mathbb{C}_p \otimes_K \mathcal{A}_K(r))$ est égal à $\chi(P, \mathcal{A}_{\mathbb{C}_p}(r))$.

Nous allons voir que, pour tout réel $r$, l'espace:

$$\operatorname{Hom}_{K[x, \frac{d}{dx}]}(K[x, \frac{d}{dx}]/P, \mathcal{A}_{\mathbb{C}_p}(r)/\mathbb{C}_p \otimes_K \mathcal{A}_K(r))$$

est *nul*. En effet les séries de $\mathbb{C}_p \otimes_K \mathcal{A}_K(r)$ apparaissent comme des séries formelles à coefficients dans un $K$-sous-espace vectoriel de $\mathbb{C}_p$ de dimension *finie*. Mais $P$ étant un polynôme différentiel à coefficients dans le corps $K$ l'équation

$$P(f) = g$$

entraîne que si la série $g$ est à coefficients dans un $K$-sous-espace vectoriel de $\mathbb{C}_p$ de dimension finie $L$, la série $f$ est aussi à coefficients dans un $K$-sous-espace vectoriel de $\mathbb{C}_p$ de dimension *finie* éventuellement plus grand $L'$.

Soit $\lambda_1, \ldots, \lambda_h$ une $K$-base de $L'$. Alors la série $f = \sum_{k=0}^{\infty} a_k x^k$ est la somme $\sum_{i=1}^{h} \lambda_i f_i$ où les séries $f_i = \sum_{k=0}^{\infty} a_{i,k} x^k$ sont à coefficients dans $K$. Il résulte du théorème de l'équivalence des normes dans un espace vectoriel de dimension finie sur un corps valué complet ([$G_1$, Chap. I 12], [$C_2$, 1.4.2]) qu'il existe une constante $C > 0$ telle que l'on ait les inégalités:

$$|a_{i,k}| \leq C|a_k|$$

pour tout $i$. Cela montre que si la série $f$ converge dans le disque $D(0, r^-)$ il en est de même des séries $f_i, i = 1, \ldots, h$.

On obtient la suite exacte

$$0 \to \operatorname{Ext}^1_{K[x, \frac{d}{dx}]}(K[x, \frac{d}{dx}]/P, \mathbb{C}_p \otimes_K \mathcal{A}_K(r)) \to \operatorname{Ext}^1_{K[x, \frac{d}{dx}]}(K[x, \frac{d}{dx}]/P, \mathcal{A}_{\mathbb{C}_p}(r))$$

$$\to \operatorname{Ext}^1_{K[x, \frac{d}{dx}]}(K[x, \frac{d}{dx}]/P, \mathcal{A}_{\mathbb{C}_p}(r)/\mathbb{C}_p \otimes_K \mathcal{A}_K(r)) \to 0.$$

Supposons que l'espace $\operatorname{Ext}^1_{K[x, \frac{d}{dx}]}(K[x, \frac{d}{dx}]/P, \mathcal{A}_{\mathbb{C}_p}(r))$ est de dimension finie, alors la topologie quotient sur cet espace induite par l'espace $\mathcal{A}_{\mathbb{C}_p}(r)$ est *séparée*. En effet soit $H$ un supplémentaire algébrique de l'image de

$$P : \mathcal{A}_{\mathbb{C}_p}(r) \to \mathcal{A}_{\mathbb{C}_p}(r)$$

qui est de codimension finie. Cette image apparaît comme le complémentaire de l'image d'un ouvert par l'application naturelle continue $\mathcal{A}_{\mathbb{C}_p}(r) \times H \to \mathcal{A}_{\mathbb{C}_p}(r)$ qui à $(f, g)$ associe $P(f) + g$. L'image est donc *fermée* en vertu du théorème des homomorphismes de Banach pour les espaces de type $\mathcal{F}$ ([$G_1$, Chap. I 14]). Cela entraîne que l'espace

$$\mathrm{Ext}^1_{K[x,\frac{d}{dx}]}(K[x,\frac{d}{dx}]/P, \mathcal{A}_{\mathbb{C}_p}(r)/\mathbb{C}_p \otimes_K \mathcal{A}_K(r))$$

est nul. En effet l'espace $\mathbb{C}_p \otimes_K \mathcal{A}_K(r)$ qui contient l'espace dense $\mathbb{C}_p[x]$ dans $\mathcal{A}_{\mathbb{C}_p}(r)$ est lui même dense dans $\mathcal{A}_{\mathbb{C}_p}(r)$. Ceci entraîne que l'image du morphisme:

$$\mathrm{Ext}^1_{K[x,\frac{d}{dx}]}(K[x,\frac{d}{dx}]/P, \mathbb{C}_p \otimes_K \mathcal{A}_K(r)) \to \mathrm{Ext}^1_{K[x,\frac{d}{dx}]}(K[x,\frac{d}{dx}]/P, \mathcal{A}_{\mathbb{C}_p}(r))$$

est partout dense, donc *surjectif* puis bijectif. L'indice $\chi(P, \mathcal{A}_{\mathbb{C}_p}(r))$ est égal à l'indice $\chi(P, \mathcal{A}_K(r))$. D'où le théorème 7.4-1.

Sous les hypothèses du théorème, l'indice $\chi(P, \mathcal{A}_K(r))$ est limite de la suite stationnaire $\chi(P, \mathcal{A}_K(r_n))$ pour une suite de nombre réels $r_n$ tendant vers $r$ par valeurs inférieures ([C-M$_1$, 4.5.1]). Remarquons dans cette situation *rationnelle* que les indices suivants sont nuls: $\chi(\mathcal{M}_{\mathrm{sol}>0}, \mathcal{R}_K(r)) = \chi(\mathcal{M}_{\mathrm{sol}>0}, \mathcal{R}_{\mathbb{C}_p}(r)) = 0$.

Considérons plus généralement un module holonome $\mathcal{M}$ sur l'algèbre de Weyl $K[x, \frac{d}{dx}]$, c'est-à-dire un $K[x, \frac{d}{dx}]$-module de type fini dont la restriction à un ouvert non vide de la droite affine est libre de type fini sur l'anneau des fonctions régulières de cet ouvert. Pour tout réel $r > 0$ on a une extension $K[x, \frac{d}{dx}] \to E_K^\dagger(r)[\frac{d}{dx}]$ qui permet de définir l'exposant $\mathfrak{E}xp_0^r(\mathcal{M})$ quand le corps $K$ est localement compact et que $r$ appartient au groupe des valeurs absolues de $K$. On obtient à partir du théorème 7.4-1, à l'aide du lemme du vecteur cyclique par exemple:

COROLLAIRE 7.4-2. *Soit un module holonome $\mathcal{M}$ sur l'algèbre de Weyl $K[x, \frac{d}{dx}]$ à coefficients dans un sous-corps complet localement compact $K$ de $\mathbb{C}_p$ et un nombre réel $r > 0$ appartenant au groupe des valeurs absolues de $K$. Si l'exposant $\mathfrak{E}xp_0^r(\mathcal{M})$ de $\mathcal{M}$ a la propriété ($\mathbf{NL^{**}}$), les indices*

$$\chi(\mathcal{M}, \mathcal{A}_{\mathbb{C}_p}(r)), \chi(\mathcal{M}, \mathcal{R}_{\mathbb{C}_p}(r)), \chi(\mathcal{M}, \mathcal{H}_{\mathbb{C}_p}^\dagger(r)),$$

$$\chi(\mathcal{M}, \mathcal{A}_K(r)), \chi(\mathcal{M}, \mathcal{R}_K(r)), \chi(\mathcal{M}, \mathcal{H}_K^\dagger(r))$$

*sont finis. De plus on a les égalités*:

$$\chi(\mathcal{M}, \mathcal{R}_{\mathbb{C}_p}(r)) = \chi(\mathcal{M}, \mathcal{A}_{\mathbb{C}_p}(r)) + \chi(\mathcal{M}, \mathcal{H}_{\mathbb{C}_p}^\dagger(r)) = 0,$$

$$\chi(\mathcal{M}, \mathcal{A}_{\mathbb{C}_p}(r)) = \chi(\mathcal{M}, \mathcal{A}_K(r)) = -\chi(\mathcal{M}, \mathcal{H}_{\mathbb{C}_p}^\dagger(r)) = -\chi(\mathcal{M}, \mathcal{H}_K^\dagger(r)).$$

### 7.5. *Le théoréme de l'indice dans le cas d'une structure de Frobenius.*

*Définition 7.5-1.* On dit qu'un $E_K^\dagger(1)$-*espace de dimension finie à connexion $\mathcal{M}$ a une structure de Frobenius dans la classe résiduelle de zéro s'il existe une structure de Frobenius d'ordre $q$, une puissance de $p$, sur l'étendu $\mathcal{R}_K(1) \otimes_{E_K^\dagger(1)} \mathcal{M}$.

THÉORÈME 7.5-2.    *Soit un module holonome $\mathcal{M}$ sur l'algèbre de Weyl $K[x, \frac{d}{dx}]$ à coefficients dans un sous-corps complet localement compact $K$ de $\mathbb{C}_p$ muni d'une structure de Frobenius; alors les indices*

$$\chi(\mathcal{M}, \mathcal{A}_{\mathbb{C}_p}(1)), \ \ \chi(\mathcal{M}, \mathcal{R}_{\mathbb{C}_p}(1)), \ \ \chi(\mathcal{M}, \mathcal{H}^{\dagger}_{\mathbb{C}_p}(1)),$$

$$\chi(\mathcal{M}, \mathcal{A}_K(1)), \ \ \chi(\mathcal{M}, \mathcal{R}_K(1)), \ \ \chi(\mathcal{M}, \mathcal{H}^{\dagger}_K(1))$$

*sont finis. De plus on a les égalités*

$$\chi(\mathcal{M}, \mathcal{R}_{\mathbb{C}_p}(1)) = \chi(\mathcal{M}, \mathcal{A}_{\mathbb{C}_p}(1)) + \chi(\mathcal{M}, \mathcal{H}^{\dagger}_{\mathbb{C}_p}(1)) = 0,$$

$$\chi(\mathcal{M}, \mathcal{A}_K(1)) = \chi(\mathcal{M}, \mathcal{A}_{\mathbb{C}_p}(1)) = -\chi(\mathcal{M}, \mathcal{H}^{\dagger}_K(1)) = -\chi(\mathcal{M}, \mathcal{H}^{\dagger}_{\mathbb{C}_p}(1)).$$

*Démonstration.* Le module $\mathcal{M}$ est soluble en 1: le module $\mathcal{M}_{\mathrm{inj}}$ est nul. En effet en vertu de 6.3-11 la limite $R(\mathcal{M}, 1^-)$ est égale à un et donc $\mathcal{M}$ est complétement soluble dans le disque générique de rayon un en vertu du théorème de continuité de la fonction rayon de convergence ([C-D$_2$, 2.5]). Le théorème de décomposition de Dwork-Robba [D-R$_1$] n'est pas nécessaire dans cette situation. Soit la décomposition 6.1-14 de $\mathcal{R}_K(1)[\frac{d}{dx}] \otimes_{K[x, \frac{d}{dx}]} \mathcal{M}$ en sa partie modérée $\mathcal{M}^{\leq 0}$ et sa partie de pentes strictement positives $\mathcal{M}^{>0}$:

$$0 \to \mathcal{M}^{>0} \to \mathcal{R}_K(1)[\frac{d}{dx}] \otimes_{K[x, \frac{d}{dx}]} \mathcal{M} \to \mathcal{M}^{\leq 0} \to 0.$$

En vertu de 6.3-12 c'est là une suite exacte de $\mathcal{R}_K(1)$-modules libres de rang fini à connexion munis de structure de Frobenius. En vertu de [C-M$_2$, 5.5-3] l'exposant $\mathfrak{e}xp^r_0(\mathcal{M}^{\leq 0})$ de la partie modérée se relève en $\{\alpha_1, \ldots, \alpha_m\}$ où les $\alpha_i$ sont des classes d'éléments de $\mathbb{Z}_p \cap \mathbb{Q}$. Donc l'exposant $\mathfrak{e}xp^r_0(\mathcal{M}^{\leq 0})$ a la propriété (**NL**\*\*). On est dans les conditions du théorème 7.4-1 et de son corollaire 7.4-2.

*Remarque* 7.5-3. Nous attirons l'attention du lecteur sur le fait que nous *ne savons pas* démontrer le théorème 7.5-2 sans passer par le théorème 7.4-1 et donc sans la structure $p$-adique générale d'un point singulier d'une équation différentielle. C'est sans doute là un point essentiel de structure dans la théorie $p$-adique qui explique son retard sur les autres théories cohomologiques des coefficients.

7.6.  *Le Théorème de finitude des nombres de Betti $p$-adiques d'une variété non singulière.* Nous allons utiliser le théorème 7.5-2 et le théorème de réduction ([Me$_4$, 3.3.6]) pour montrer le théorème de finitude 7.6-1 des nombres de Betti $p$-adiques d'une variété ouverte non singulière et leur invariance par changement de base infini.

7.6.1.  *La cohomologie de Monsky-Washnitzer.* Soit $V$ l'anneau des entiers d'un sous corps $K$ à valuation discrète du corps $\mathbb{C}_p$ et $X$ une variété algébrique affine non singulière sur le corps résiduel $k$ d'algèbre $\bar{A}$. En vertu du théorème d'Elkik [E] il existe une $V$-algèbre $A$ de type fini non singulière dont la réduction modulo l'idéal maximum $\mathfrak{m}$ est isomorphe à $\bar{A}$. Si $V[x_1,\ldots,x_n]/I \simeq A$ est une présentation de $A$ on définit l'algèbre $A^\dagger$ comme $(V[x_1,\ldots,x_n])^\dagger/I$ où $(V[x_1,\ldots,x_n])^\dagger$ désigne l'algèbre des séries à coefficients dans $V$ qui admettent un rayon de convergence strictement plus grand que 1. En vertu du théorème de M. Artin le couple $(A^\dagger, \hat{A})$ a la propriété d'approximation où $\hat{A}$ désigne le complété $\mathfrak{m}$-adique de $A$. On peut montrer alors que l'algèbre $A^\dagger$ est formellement très lisse [M-W]. En fait il y a *équivalence* entre être formellement très lisse et formellement lisse et l'hypothèse imposée par Monsky-Washnitzer n'est pas restrictive comme prévu ([M-W, 3, p. 189]).

On définit ([M-W]) la cohomologie de de Rham $p$-adique $H_{\mathrm{DR}}^\bullet(A^\dagger/V)$ de $A^\dagger$, comme la cohomologie du complexe de de Rham $D(A^\dagger)$ des formes différentielles $\mathfrak{m}$-séparées de l'algèbre $A^\dagger$. Si $V$ est non ramifié la cohomologie $p$-adique $H^\bullet(X;V) := H_{\mathrm{DR}}^\bullet(A^\dagger/V)$ ne dépend pas du relèvement $A^\dagger$ et est fonctorielle de façon contravariante en $X$ [M-W]. Si $V$ est éventuellement ramifié la cohomologie $p$-adique $H^\bullet(X;K) := H_{\mathrm{DR}}^\bullet(A^\dagger/V) \otimes_{\mathbb{Z}} \mathbb{Q}$ ne dépend pas du relèvement $A^\dagger$ et est fonctorielle de façon contravariante en $X$ [M-W].

7.6.2 *Le théorème de finitude.* Nous notons $B_{i,p}(X)$ la dimension des $K$-espaces $H^i(X;K)$. Les nombres de Betti $B_{i,p}(X)$ sont nuls pour $i > \dim X$ et en vertu du théorème de Monsky [Mo₃] les nombres $B_{0,p}(X)$ et $B_{1,p}(X)$ sont finis. On obtient finalement le résultat qui était le principal obstacle dans la théorie $p$-adique de la fonction zêta :

THÉORÈME 7.6-1. *Pour toute variété algébrique $X$ affine non singulière sur le corps résiduel $k$ les nombres de Betti $B_{i,p}(X)$ sont finis pour tout $i$ et sont invariants par changement de base d'anneaux de valuation discrète.*

*Démonstration.* Soit $f(x) = f(x_1,\ldots,x_n)$ un polynôme à coefficients dans l'anneau des entiers $\mathcal{O}_{\mathbb{C}_p}$ du corps $\mathbb{C}_p$ et $m$ un entier positif. On définit le module exponentiel $M_{f,n,m}$ comme le quotient de $\mathcal{O}_{\mathbb{C}_p}[x_1,\ldots,x_n,\Gamma,\Gamma^{-1}]$ par les images des opérateurs différentiels

$$\partial_{x_i} + \pi(\partial_{x_i}(f) + m\Gamma x_i^{m-1})$$

pour $i = 1,\ldots,n$. Si $m$ est strictement plus grand que le degré total de $f$ le module $M_{f,n,m}$ est libre de rang $(m-1)^n$ sur l'anneau $\mathcal{O}_{\mathbb{C}_p}[\Gamma,\Gamma^{-1},m^{-1}]$ et est muni d'une connexion. L'infini est une singularité régulière et zéro est une singularité irrégulière.

Si de plus $f = \mathrm{Teich}(\bar{f})$ est le relèvement de Teichmüller d'un polynôme $\bar{f}$ à coefficients dans le corps résiduel et que $m$ est premier avec $p$, le module différentiel $M_{\mathrm{Teich}(\bar{f}),n,m}$ qui est défini sur une extension finie de $\mathbb{Q}_p$ est muni d'une structure de Frobenius ([Me₄, 4.1.1]) sur l'algèbre $(K[\Gamma,\Gamma^{-1}])^\dagger$.

La structure de la singularité régulière à l'infini du module différentiel $M_{\mathrm{Teich}(\bar{f}),n,m}$ est conséquence du théorème de transfert [C$_3$], alors que la structure de la singularité irrégulière en zéro est conséquence de [C-M$_1$], [C-M$_2$] et du présent travail. Les équations $M_{\mathrm{Teich}(\bar{f}),n,m}$ illustrent d'autre part tous les résultats de l'article [C-M$_4$]. Nous ferons remarquer au lecteur que la classe des modules différentiels $M_{f,n,m}$ qu'on construit explicitement à partir de polynômes, a été le principal catalyseur dans la structure $p$-adique générale d'une équation différentielle. Inversement nous ne savons pas démontrer l'existence de l'indice pour la classe des modules exponentiels $M_{\mathrm{Teich}(\bar{f}),n,m}$ sans passer par la théorie générale.

En vertu du théorème de réduction ([Me$_4$, 3.3.6]) qui est parallèle au théorème de réduction en caractéristique nulle [Mo$_4$], la finitude de l'indice local de $M_{\mathrm{Teich}(\bar{f}),n,m}$ dans la classe résiduelle de zéro pour tout triplet $(\bar{f},n,m)$ tel que $m > \deg(\bar{f})$ et $(p,m) = 1$ entraîne la finitude des nombres Betti $B_{i,p}(X)$ pour tout $i$ et pour toute variété $X$ algébrique affine non singulière sur le corps résiduel et leur invariance par changement de base. Le théorème 7.5-2 entraîne alors le théorème 7.6-1.

Dans le cas complexe une réduction similaire, indépendante de l'article de Monsky [Mo$_4$], est faite dans [M-N$_3$]. Nous renvoyons le lecteur à l'introduction de l'article [Me$_4$] pour les différentes étapes ([M-N$_1$], [M-N$_2$], [M-N$_3$]) qui nous a conduit à ce point de vue. On trouvera d'autre part dans cet article les exemples d'équations de la classe $M_{\mathrm{Teich}(\bar{f}),n,m}$ montrant les diverses situations d'un point singulier d'une équation différentielle $p$-adique que l'on rencontre.

Cette méthode ramène aussi le théorème de pureté des valeurs propres de l'endomorphisme de Frobenius [Mo$_2$] opérant sur la cohomologie $p$-adique d'une variété affine non singulière sur un corps fini au résultat analogue pour la cohomologie $p$-adique des équations $M_{\mathrm{Teich}(\bar{f}),n,m}$ ([Me$_4$, 4.2.2]). Il est donc important de transposer à la situation $p$-adique les arguments $\ell$-adiques de Deligne-Laumon ([De$_2$], [De$_3$], [La]). On dispose pour cela de la catégorie des coefficients $p$-adiques sur les courbes ([C-M$_2$], [C-M$_4$]) qui a toutes les propriétés de finitude de la catégorie des coefficients $\ell$-adiques ce qui est essentiel. En particulier la cohomologie $p$-adique intermédiaire d'un fibré $p$-adique à connexion ayant une structure de Frobenius est définie et est de dimension *finie* [C-M$_4$].

Les arguments de Grothendieck ([G$_2$], [G$_3$]) permettent de passer du théorème de finitude de la cohomologie $p$-adique des variétés affines non-singulières sur $k$ au théorème de finitude de la cohomologie $p$-adique des variétés non-singulières sur $k$ telles qu'elle est définie dans [G$_3$, 2.2] comme l'hypercohomologie du site naturel des relévements locaux †-adiques à valeurs dans le complexe de de Rham.

Comme résultat préalable au théorème 7.6-1 nous obtenons la finitude de la cohomologie $p$-adique de tous les modules exponentiels $\exp(\pi\mathrm{Teich}(\bar{f}))$ pour tous les polynômes $\bar{f}$ et de leur cohomologie locale le long de toute hypersurface de l'espace affine, en fait la finitude de la cohomologie $p$-adique de tous les fibrés algébriques à connexion intégrable sur l'espace affine munis d'une structure de Frobenius.

**7.7.** *La conjecture de la propriété* (**NL**\*\*) *des exposants.* Soit $P(x, \frac{d}{dx})$ un polynôme différentiel à coefficients dans le corps des nombres algébriques $\bar{\mathbb{Q}}$, il provient d'un corps $K$ localement compact. Pour tout nombre rationnel $r > 0$ son exposant $\mathfrak{E}xp_0^r(P)$ est défini par voie purement $p$-adique. Nous rappelons la conjecture [C-M$_1$, 3.3.5] :

CONJECTURE 7.7-1. *L'exposant* $\mathfrak{E}xp_0^r(P)$ *a la propriété* (**NL**\*\*). *De plus pour* $r$ *variable il n'y a qu'un nombre fini d'exposants* $\mathfrak{E}xp_0^r(P)$.

La conjecture est vraie en rang un en vertu du théorème de Robba [R$_4$] et pour les équations hypergéométriques un résultat partiel a été obtenu dans cette direction ([C-M$_1$, 3.3.7]). Il est naturel de penser dans cette situation que les exposants sont algébriques. Mais comme nous l'a fait remarquer Y. André les exposants de la monodromie complexe de $P(x, \frac{d}{dx})$ ne sont pas en général algébriques. Cependant nous ne connaissons pas d'exemple de ce phénomène dans le cas $p$-adique qui est tout de même assez différent du cas complexe.

**7.8.** *Cas d'une classe résiduelle d'une courbe de genre supérieur.* Les résultats précédents, qui sont de nature locale, ne dépendent pas de la coordonnée $x$ et se transpose sur toute courbe $X_K$ non singulière sur un sous corps complet $K$ de $\mathbb{C}_p$ qui est extension finie de $\mathbb{Q}_p$, quitte à choisir un modèle entier. Pour tout réel $r, 0 < r \leq 1$, et tout point $x$ de $X_K^{\mathrm{an}}$ de la variété analytique associé à $X_K$ on peut parler du disque $D(x, r^-)$ et des espaces

$$\mathcal{A}_{Kx}(r), \mathcal{A}_{Kx}([r-\varepsilon, r[), \mathcal{R}_{Kx}(r), \mathcal{H}_{Kx}^\dagger(r), E_{Kx}^\dagger(r), \mathcal{A}_{Kx}([r-\varepsilon][\frac{d}{dx}].$$

Si $\gamma$ est un nombre réel la norme $|P|_{\gamma,\rho}$ d'un opérateur différentiel est la norme de Banach de $P$ opérant l'espace des fonctions bornées dans le disque générique $D(t_\rho, \rho^{\gamma-})$. On définit alors les topologies $\mathcal{T}_\gamma$ sur les espaces $\mathcal{A}_{Kx}([r-\varepsilon, r[)[\frac{d}{dx}]$ et les topologies quotients $\mathcal{T}_{\beta,Q}$ sur les modules de type fini sur les anneaux $\mathcal{A}_{Kx}([r-\varepsilon, r[)[\frac{d}{dx}]$. Pour un module libre de rang fini à connexion $\mathcal{M}$ sur l'anneau $\mathcal{R}_{Kx}(r)$ on a la fonction rayon de convergence $R(\mathcal{M}, \rho)$ ([C-M$_2$, 3.1.3]). On peut définir les catégories $\mathrm{MLS}(\mathcal{R}_{Kx}(r))$. Pour un module de la catégorie $\mathrm{MLS}(\mathcal{R}_{Kx}(r))$ on définit sa filtration $\mathcal{M}_{>\gamma}$ et sa cofiltration $\mathcal{M}^{\leq\gamma}$ du moins si $r$ appartient au groupe des valeurs absolues de $K$ et en particulier sa décomposition en partie modérée et partie de pentes strictement positives. D'où la notion d'exposant, qui est intrinsèque ([C-M$_2$, 5.5.4]) et donc des catégories $\mathrm{MLS}(\mathcal{R}_{Kx}(r), \mathbf{NL}^*)$, $\mathrm{MLS}(\mathcal{R}_{Kx}(r), \mathbf{NL}^{**})$, $\mathrm{MLS}(\mathcal{R}_{Kx}(1), \mathbf{F})$. On aura les théorèmes 7.4-1 et 7.5-2 pour les modules holonomes sur $X/K$. Pour les démonstrations le seul point qu'il faut savoir en plus est que l'espace $\mathrm{Hom}_{\mathcal{D}_{X/K}}(\mathcal{M}, \mathcal{H}_{Kx}^{\dagger}(r))$ est de dimension finie sur $K$. La question est locale pour la topologie de Zariski, si $A$ est l'algèbre affine d'un ouvert affine de $X/K$ voisinage de $x$ qui est une extension finie étale d'une algèbre $K[x, 1/g]$ provenant d'une situation similaire sur le corps résiduel, alors l'anneau des opérateurs différentiels $D_{A/K}$ est une extension de $K[x, 1/g][\frac{d}{dx}]$. Un module holonome $\mathcal{M}$ sur l'anneau $D_{A/K}$ reste holonome sur l'anneau $K[x, 1/g][\frac{d}{dx}]$, en particulier il est de type fini. C'est là la propriété essentielle de stabilité de la catégorie des modules holonomes par image directe. Ceci montre que l'espace $\mathrm{Hom}_{D_{A/K}}(\mathcal{M}, \mathcal{H}_{Kx}^{\dagger}(r))$ est un sous-espace de l'espace $\mathrm{Hom}_{K[x,1/g][\frac{d}{dx}]}(\mathcal{M}, \mathcal{H}_{Kx}^{\dagger}(r))$ qui est de dimension finie.

## 8. La formule de l'indice

Nous allons établir la formule de l'indice local conjecturée par Robba [R$_5$] et montrer que les sommets du polygone de Newton d'un $\mathcal{R}_K(r)$-module libre de rang fini muni d'une connexion et soluble sont à coordonnées *entières* lorsque le corps de base est localement compact. Nous définissons dans la situation locale les nombres $\mathrm{Irr}_x(\mathcal{M}, p)$. Nous avons besoin d'une part du théorème de l'invariance de l'indice par perturbation compacte pour les espaces de type $\mathcal{LF}$ ([G$_1$, Chap. V]) et d'autre part de la notion d'indice généralisé introduite par Robba [R$_5$].

8.1. *Le théorème de perturbation compacte.* Soit $K$ un corps valué complet. Rappelons que l'on dit qu'une application linéaire $u$ entre deux $K$-espaces vectoriels topologiques localement convexes $E$ et $F$ séparés, est compacte si elle transforme un voisinage convenable de zéro en une partie relativement compacte; cf. [G$_1$, Chap. V]. On a alors le théorème de l'invariance de l'indice par perturbation compacte ([G$_1$, Chap. V]):

THÉORÈME 8.1-1. *Supposons que le corps $K$ est localement compact et soit $v$ une application linéaire continue à indice entre deux $K$-espaces vectoriels topologiques de type $\mathcal{LF}$ et $u$ une application linéaire compacte, alors $v + u$ est à indice et l'on a l'égalité des indices*

$$\chi(u + v) = \chi(v).$$

*Exemple* 8.1-2. Une application complètement continue $u$ entre deux espaces de Banach sur un corps $p$-adique localement compact $K$ est compacte. En effet, par définition, $u$ est limite pour la norme de Banach d'applications de rang fini. Elle transforme la boule unité en partie précompacte, donc relativement compacte.

8.2. *La notion d'indice généralisé.* Soit $K$ un corps complet muni d'une valeur absolue $p$-adique et un réel $r > 0$. Nous ne considérons dans ce paragraphe que des espaces vectoriels sur $K$.

8.2.1. *Cas d'une couronne.* Soit $I$ un intervalle tel que $I = I_+ \cap I_-$ où $I_+$ est un intervalle d'extrémité zéro et $I_-$ est intervalle d'extrémité $\infty$. On a la décomposition

$$\mathcal{A}_K(I) = \mathcal{A}_K(I_+) \bigoplus \frac{1}{x} \mathcal{A}_K(I_-).$$

On note $\gamma^+$, resp. $\gamma^-$, la projection de $\mathcal{A}_K(I)$ sur $\mathcal{A}_K(I_+)$, resp. sur $\frac{1}{x}\mathcal{A}_K(I_-)$. De même on note $\gamma_+$, resp. $\gamma_-$ l'injection de $\mathcal{A}_K(I_+)$, resp. de $\frac{1}{x}\mathcal{A}_K(I_-)$, dans $\mathcal{A}_K(I)$.

*Définition* 8.2-1. Soit $m$ un entier et $u(I)$ un endomorphisme de l'espace $(\mathcal{A}_K(I))^m$. On définit *les indices généralisés* de $u(I)$ comme les indices

$$\widetilde{\chi}(u(I), \mathcal{A}(I_+)) \quad := \quad \chi(\gamma^+ \circ u(I) \circ \gamma_+, (\mathcal{A}_K(I_+))^m)$$

et

$$\widetilde{\chi}(u(I), \frac{1}{x}\mathcal{A}_K(I_-)) \quad := \quad \chi(\gamma^- \circ u(I) \circ \gamma_-, (\frac{1}{x}\mathcal{A}_K(I_-))^m).$$

Si $u$ est un endomorphisme de $(\mathcal{A}_K(I))^m$ nous noterons pour simplifier $\gamma^+ u \gamma^-$ pour $\gamma_+ \circ \gamma^+ \circ u \circ \gamma_- \circ \gamma^-$ et $\gamma^- u \gamma^+$ pour $\gamma_- \circ \gamma^- \circ u \circ \gamma_+ \circ \gamma^+$.

PROPOSITION 8.2-2. *Supposons le corps de base localement compact et soit $u$ une matrice carré d'ordre $m$ à coefficients dans l'anneau $\mathcal{A}_K(I)[\frac{d}{dx}]$; alors les endomorphismes $\gamma^- u \gamma^+$ et $\gamma^+ u \gamma^-$ de $(\mathcal{A}_K(I))^m$ muni de sa topologie d'espace de type $\mathcal{F}$ sont compacts.*

*Démonstration.* Notons $((\mathcal{A}_K(I))^m, \mathcal{T}_0)$ l'espace $(\mathcal{A}_K(I))^m$ muni de la topologie canonique de type $\mathcal{F}$ et $((\mathcal{A}_K(I))^m, |-|_r)$ l'espace $(\mathcal{A}_K(I))^m$ muni de la norme $|-|_r$ pour un nombre réel $r$ de l'intervalle $I$. Soit

$$L_b(((\mathcal{A}_K(I))^m, |-|_r), ((\mathcal{A}_K(I))^m, \mathcal{T}_0))$$

l'ensemble des applications linéaires continues $K$-linéaires, muni de la topologie de la convergence uniforme sur les parties bornées ([G$_1$, Chap. III]), de $((\mathcal{A}_K(I))^m, |-|_r)$ dans l'espace topologique $((\mathcal{A}_K(I))^m, \mathcal{T}_0)$.

Si $u$ est une matrice à coefficients dans l'anneau $K[x, 1/x][\frac{d}{dx}]$ les applications $\gamma^+ u \gamma^-$ et $\gamma^- u \gamma^+$ sont de rang fini, en particulier transforment la boule unité de $((\mathcal{A}_K(I))^m, |-|_r)$ en partie relativement compacte (le corps $K$ étant localement compact, elles sont donc complètement continues).

Si $u$ est une matrice à coefficients dans l'anneau $\mathcal{A}_K(I)[\frac{d}{dx}]$ elle est limite de matrices $u_n$ à coefficients dans l'anneau $K[x, 1/x][\frac{d}{dx}]$.

Si $s$ et $v$ sont des entiers de $\mathbb{Z}$ et $\ell$ un entier naturel de $\mathbb{N}$ on a $\gamma^+ \circ x^s (\frac{d}{dx})^\ell \circ \gamma^- (x^v) = 0$ si $v > 0$ ou si $s - \ell + v < 0$. Pour $a, b \in K$, $\rho, r \in I$, $R = \max(r, \rho)$, on en déduit:

$$|\gamma^+ \circ ax^s (\frac{d}{dx})^\ell \circ \gamma^- (bx^v)|_\rho \leq |ab| \, \rho^{s-\ell+v}$$

$$\leq |ax^s|_R |bx^v|_r (\rho/R)^{s-\ell+v} (r/R)^{-v} R^{-\ell},$$

c'est-à-dire, pour $f \in \mathcal{A}_K(I)$:

$$|\gamma^+ \circ (u - u_n) \circ \gamma^- (f))|_\rho \leq |f|_r ||u - u_n||_R \max(1, I/R)^{\text{ordre}(u)}.$$

On voit que, les suites $\gamma^- u_n \gamma^+$ et $\gamma^+ u_n \gamma^-$ convergent uniformément sur la boule unité de $((\mathcal{A}_K(I))^m, |-|_r)$ vers $\gamma^- u \gamma^+$ et $\gamma^+ u \gamma^-$ .

Mais (cf. $[G_1$, Chap. 0, 4.1.6']) l'ensemble des endomorphismes linéaires continus qui transforment la boule unité en partie précompacte est fermé dans $L_b(((\mathcal{A}_K(I))^m, |-|_r), ((\mathcal{A}_K(I))^m, \mathcal{T}_0))$. Comme l'espace $((\mathcal{A}_K(I))^m, \mathcal{T}_0)$ est complet, il y a identité entre parties relativement compactes et partie précompactes. Donc les applications $\gamma^- u \gamma^+$ et $\gamma^+ u \gamma^-$ transforment cette boule unité en partie relativement compacte. Mais la boule unité pour la norme $|-|_r$ est un voisinage de l'espace métrique $((\mathcal{A}_K(I))^m, \mathcal{T}_0)$. Les endomorphismes $\gamma^- u \gamma^+$ et $\gamma^+ u \gamma^-$ de $(\mathcal{A}_K(I))^m$ sont donc *compacts*.

PROPOSITION 8.2-3. *Supposons le corps de base $K$ localement compact, soient $u$ et $v$ deux matrices carrées d'ordre $m$ à coefficients dans l'anneau $\mathcal{A}_K(I)[\frac{d}{dx}]$ pour un intervalle $I$. Si deux des matrices $u, v, uv$ ont un indice généralisé, la troisième a un indice généralisé et l'on a l'égalité:*

$$\widetilde{\chi}(uv, \mathcal{A}_K(I)) = \widetilde{\chi}(u, \mathcal{A}_K(I)) + \widetilde{\chi}(v, \mathcal{A}_K(I)).$$

*Démonstration.* En effet on a l'égalité

$$\gamma^+ u \gamma^+ v = \gamma^+ u (v - \gamma^- v) = \gamma^+ uv - \gamma^+ u \gamma^- v.$$

Mais l'endomorphisme $\gamma^+ u \gamma^- v$ est compact en vertu de la démonstration précédente. La proposition 8.2-3 est conséquence du théorème 8.1-1.

PROPOSITION 8.2-4. *Supposons le corps de base $K$ localement compact, soit $u$ une matrice carrée d'ordre $m$ à coefficients dans l'anneau $\mathcal{A}_K(I)[\frac{d}{dx}]$ pour un intervalle $I$. Alors $u$ est à indice dans l'espace $\mathcal{A}_K(I)$ si et seulement si elle admet des indices généralisés et l'on a l'égalité:*

$$\chi(u, \mathcal{A}_K(I)) = \widetilde{\chi}(u, \mathcal{A}_K(I_+)) + \widetilde{\chi}(u, \tfrac{1}{x}\mathcal{A}_K(I_-)).$$

*Démonstration.* Considérons l'endomorphisme:

$$\tilde{u} := u - \gamma^+ u \gamma^- - \gamma^- u \gamma^+$$

qui respecte la décomposition

$$\mathcal{A}_K(I) = \mathcal{A}_K(I_+) \bigoplus \frac{1}{x} \mathcal{A}_K(I_-).$$

Mais l'endomorphisme $\gamma^+ u \gamma^- + \gamma^- u \gamma^+$ de $\mathcal{A}_K(I)$ est compact. La restriction de $\tilde{u}$ à $\mathcal{A}_K(I_+)$ est égal à $\gamma^+ u$ et la restriction de $\tilde{u}$ à $\frac{1}{x} \mathcal{A}_K(I_-)$ est égal à $\gamma^- u$. La proposition 8.2-4 est conséquence de la proposition 8.2-2.

Pour une fonction $f$ de $\mathcal{A}_K(I)$ où $I$ est un intervalle, on note $\mathrm{ord}_I(f)$ le degré, éventuellement infini, du diviseur des zéros de $f$ contenus dans $C(I)$.

*Définition* 8.2-5. Si l'intervalle $I$ est non réduit à un point, on définit l'ordre généralisé

$$\widetilde{\mathrm{ord}}_{I_+}(f) := d \, \mathrm{Log}^- |f|(\rho) + \mathrm{ord}_{[\rho, \infty[ \cap I}(f)$$

pour un nombre $\rho$ de l'intérieur de $I$.

Cette définition est indépendante du nombre $\rho$. On définit de la même manière $\widetilde{\mathrm{ord}}_{I_-}(f)$.

PROPOSITION 8.2-6. *Supposons le corps de base $K$ localement compact, soit $u$ une matrice carrée d'ordre $m$ à coefficients dans l'anneau $\mathcal{A}_K(I)$ pour un intervalle $I$ non réduit à un point. On a alors les égalités*:

$$\widetilde{\chi}(u, \mathcal{A}_K(I_+)) = -\widetilde{\mathrm{ord}}_{I_+}(\det(u)),$$
$$\widetilde{\chi}(u, \frac{1}{x} \mathcal{A}_K(I_-)) = \widetilde{\mathrm{ord}}_{I_-}(\det(u)).$$

*Démonstration.* Supposons d'abord que la matrice $u$ est inversible. En vertu de la décomposition de Birkhoff en facteurs singuliers [C_4] la matrice $u$ est produit $x^a u_1 u_2$ où $u_1$ est matrice inversible à coefficients dans $\mathcal{A}_K(I_+)$, $u_2$ une matrice inversible à coefficients dans $\mathcal{A}_K(I_-)$ et $x^a = (x^{a_1}, \ldots, x^{a_m})$ une matrice diagonale. En vertu de la proposition 8.2-4

$$\chi(u_1, \mathcal{A}_K(I_+)) = \widetilde{\chi}(u_2, \mathcal{A}_K(I_+)) = 0.$$

En vertu de la proposition 8.2-3

$$\widetilde{\chi}(u, \mathcal{A}_K(I_+)) = \widetilde{\chi}(x^a, \mathcal{A}_K(I_+)) = -\sum_i a_i = -\widetilde{\mathrm{ord}}_{I_+}(\det(u)).$$

Pour un intervalle fermé $I$ l'anneau $\mathcal{A}_K(I)$ est principal. La théorie des facteurs invariants montre que la matrice $u$ est égale au produit de matrices carrées $u_1 d u_2$ où $u_1$ et $u_2$ sont inversibles et $d$ diagonale. La proposition 8.2-6 est conséquence de la proposition 8.2-3. Le cas d'un intervalle quelconque s'obtient par passage à la limite, en vertu de la propriété (**ML**) de Mittag-Leffler pour les systèmes projectifs d'espaces métriques complets à images denses ([EGA III], Chap. 0, 13.2.4), à partir du cas des intervalles fermés.

Soit $\mathcal{M}$ un module différentiel sur l'anneau $\mathcal{A}_K(I)$ de rang $m$ et $G(x)$ la matrice représentant l'opérateur $x\frac{d}{dx}$ dans une base.

*Définition* 8.2-7. On définit les indices généralisés

$$\widetilde{\chi}(\mathcal{M}, \mathcal{A}_K(I_+)), \quad \widetilde{\chi}(\mathcal{M}, \frac{1}{x}\mathcal{A}_K(I_-))$$

comme *les indices généralisés de* $x\frac{d}{dx} - G(x)$.

En vertu de la proposition 8.2-3 et de la proposition 8.2-6 les indices généralisés ne dépendent pas de la base choisie.

Proposition 8.2-8. *Si le corps de base $K$ est localement compact, dans une suite exacte de $\mathcal{A}_K(I)$-modules différentiels pour un intervalle*

$$0 \to \mathcal{M}_1 \to \mathcal{M} \to \mathcal{M}_2 \to 0,$$

*les indices généralisés sont additifs.*

*Démonstration.* En effet si on choisit pour base de $\mathcal{M}$ la somme d'une base de $\mathcal{M}_1$ et d'une base de $\mathcal{M}_2$, la matrice de la connexion est triangulaire et induit une suite exacte de complexes:

$$
\begin{array}{ccccccc}
0 & \to & (\mathcal{A}_K(I_+))^{m_1} & \to & (\mathcal{A}_K(I_+))^{m} & \to & (\mathcal{A}_K(I_+))^{m_2} & \to 0 \\
 & & \downarrow \gamma^+(x\frac{d}{dx} - G_1(x)) & & \downarrow \gamma^+(x\frac{d}{dx} - G(x)) & & \downarrow \gamma^+(x\frac{d}{dx} - G_2(x)) & \\
0 & \to & (\mathcal{A}_K(I_+))^{m_1} & \to & (\mathcal{A}_K(I_+))^{m} & \to & (\mathcal{A}_K(I_+))^{m_2} & \to 0
\end{array}
$$

qui montre l'additivité des indices généralisés dans ce cas là. La proposition 8.2-8 est conséquence de l'indépendance de l'indice généralisé de la base choisie.

8.2.2. *Cas du bord d'une couronne.* Soit $r$ un nombre réel, on a la décomposition:

$$\mathcal{R}_K(r) = \mathcal{A}_K(r) \bigoplus \mathcal{H}_K^{\dagger}(r).$$

On note $\gamma^+$, resp. $\gamma^-$, la projection de $\mathcal{R}_K(r)$ sur $\mathcal{A}_K(r)$, resp. sur $\mathcal{H}_K^{\dagger}(r)$. De même on note $\gamma_+$, resp. $\gamma_-$, l'injection de $\mathcal{R}_K(r)$ sur $\mathcal{A}_K(r)$, resp. sur $\mathcal{H}_K^{\dagger}(r)$.

*Définition* 8.2-9. Soit $m$ un entier et $u(r)$ un endomorphisme de l'espace $(\mathcal{R}_K(r))^m$. On définit *les indices généralisés de* $u(r)$ comme

$$\widetilde{\chi}(u(r), \mathcal{A}_K(r)) \quad := \quad \chi(\gamma^+ \circ u(r) \circ \gamma_+, (\mathcal{A}_K(r))^m)$$

et

$$\widetilde{\chi}(u(r), \mathcal{H}_K^{\dagger}(r)) \quad := \quad \chi(\gamma^- \circ u(r)\gamma_-, (\mathcal{H}_K^{\dagger}(r))^m).$$

Si l'endomorphisme $u$ est défini dans un intervalle $]r - \varepsilon, r[$, on a alors

$$\widetilde{\chi}(u(r), \mathcal{A}_K(r)) := \widetilde{\chi}(u(]r - \varepsilon, r[), \mathcal{A}_K(r))$$

et

$$\widetilde{\chi}(u(r), \mathcal{H}_K^{\dagger}(r)) = \lim_{\varepsilon \to 0} \widetilde{\chi}(u(]r - \varepsilon, r[), \frac{1}{x}\mathcal{A}_K(]r - \varepsilon, \infty])).$$

Ceci permet d'appliquer les résultats des couronnes.

*Remarque* 8.9-10. Sur un corps maximalement complet il faut remplacer la notion de partie compacte, par la notion de partie c-compacte, c'est-à-dire dans laquelle toute intersection dénombrable de fermés convexes non vides emboités est non vide. On a alors la notion d'application c-compacte et l'invariance de l'indice par pertubation par une application c-compacte du théorème 8.1-1 avec la même démonstration. Une application complètement continue entre espaces de Banach sur un corps maximalement complet est c-compacte. Le lecteur pourra vérifier à partir de là que tous les résultats du paragraphe 8 s'étendent au cas d'un corps maximalement complet quelconque.

### 8.3. *La formule de l'indice local.*

THÉORÈME 8.3-1. *Soit $K$ un corps de caractéristique nulle complet localement compact pour une valuation $p$-adique et $\mathcal{M}$ un $\mathcal{R}_K(r)$-module libre de rang $m$ à connexion soluble en $r$ et purement de pente $\mathrm{pt}(\mathcal{M}) > 0$ alors $\mathcal{M}$ admet des indices généralisés et l'on a les égalités*

$$\widetilde{\chi}(\mathcal{M}, \mathcal{A}_K(r)) = -\widetilde{\chi}(\mathcal{M}, \mathcal{H}_K^{\dagger}(r)) = m\,\mathrm{pt}(\mathcal{M})$$

*et pour $\varepsilon > 0$ assez petit les égalités*

$$\chi(\mathcal{M}, \mathcal{A}_K(]r - \varepsilon, r[)) = \chi(\mathcal{M}, \mathcal{R}_K(r)) = 0.$$

*Démonstration.* Il suffit de montrer que, pour tout $\varepsilon > 0$, assez petit les indices généralisés $\widetilde{\chi}(\mathcal{M}, \mathcal{A}_K([0, r - \varepsilon[)$, resp. $\widetilde{\chi}(\mathcal{M}, \frac{1}{x}\mathcal{A}_K(]r - \varepsilon, \infty]))$, existent et sont égaux à $m\beta$, resp. $-m\beta$ où $\beta = \mathrm{pt}(\mathcal{M})$. Cela entraînera le Théorème 8.3-1. En effet en vertu de la propriété (ML) pour les systèmes projectifs d'espaces métriques complets à images denses ([EGA III, Chap. 0, 13.2.4])

$$\widetilde{\chi}(\mathcal{M}, \mathcal{A}_K(r)) = m\beta.$$

Par limite inductive

$$\widetilde{\chi}(\mathcal{M}, \mathcal{H}_K^{\dagger}(r)) = -m\beta.$$

En vertu de la proposition 8.2-4

$$\chi(\mathcal{M}, \mathcal{A}_K(]r-\varepsilon, r[) = \widetilde{\chi}(\mathcal{M}, \mathcal{A}_K([0, r[) + \widetilde{\chi}(\mathcal{M}, \frac{1}{x}\mathcal{A}_K(]r-\varepsilon, \infty])) = 0$$

pour $\varepsilon > 0$ assez petit. Donc par limite inductive

$$\chi(\mathcal{M}, \mathcal{R}_K(r)) = 0.$$

En vertu du théorème 4.2-1 on peut supposer que la fonction rayon de convergence $R(\mathcal{M}, \rho)$ est égale à $\rho(\rho/r)^\beta$ pour $\rho \in [r-\varepsilon, r[, \varepsilon > 0$. Considérons les intervalles $]r_{h-1}, r_h[$, définis pour $h$ assez grand, où l'on a les inégalités

$$\rho\omega^{1/p^{h-1}} < R(\mathcal{M}, \rho) < \rho\omega^{1/p^h}.$$

Soit alors $\mathcal{N}_h$ un antécédent de Frobenius de Christol-Dwork d'ordre $h$ de $\mathcal{M}$. C'est un $\mathcal{M}_K(]r_{h-1}^{p^h}, r_h^{p^h}[)$-module libre de rang $m$ à connexion qui n'a que des singularités apparentes et tel que $R(\mathcal{N}_h, \rho) = r^{-\beta p^h}\rho^{\beta+1} < \rho\omega$ pour tout $\rho \in ]r_{h-1}^{p^h}, r_h^{p^h}[$. Remarquons que, par hypothèse, toutes les solutions locales de $\mathcal{M}$ au point générique $t_\rho$ pour $\rho \in [r-\varepsilon, r[$ ont même rayon de convergence égal à $R(\mathcal{M}, \rho)$. Ceci entraîne que toutes les solutions de $\mathcal{N}_h$ locales au point générique $t_\rho$ pour $\rho \in ]r_{h-1}^{p^h}, r_h^{p^h}[$ ont même rayon de convergence égal à $R(\mathcal{N}_h, \rho)$.

LEMME 8.3-2.  *Soit $\mathcal{N}$ un $\mathcal{A}_K(]r_1, r_2[)$-module différentiel de rang $m$ dont toutes les solutions locales au point générique $t_\rho$ pour $\rho \in ]r_1, r_2[$ ont même rayon de convergence égal à $\lambda\rho^{\beta+1} < \omega\rho$. Alors pour $\rho \in ]r_1, r_2[$ on a les égalités*:

$$\widetilde{\chi}(\mathcal{N}, \mathcal{A}_K([0, \rho])) = -\widetilde{\chi}(\mathcal{N}, \frac{1}{x}\mathcal{A}_K([\rho, \infty])) = m\beta.$$

*Démonstration.* Soit

$$G(x) := \begin{pmatrix} 0 & 1 & & \\ & \ddots & \ddots & \\ & & 0 & 1 \\ a_0(x)/a_m(x) & \cdots & \cdots & a_{m-1}(x)/a_m(x) \end{pmatrix}$$

la matrice de la connexion dans une base cyclique de $\mathcal{N}$ où les coefficients $a_i$ sont des fonctions analytiques sur l'intervalle $]r_1, r_2[$. En vertu de la proposition 8.2-3 et de la définition 8.2-6 en tout point $\rho$ l'indice $\widetilde{\chi}(\mathcal{N}, \mathcal{A}_K([0, \rho]))$ + $\chi(a_m, \mathcal{A}_K([0, \rho]))$ est égale à l'indice généralisé de l'opérateur

$$P(x, x\frac{d}{dx}) = a_m(x)(x\frac{d}{dx})^m - a_{m-1}(x)(x\frac{d}{dx})^{m-1}\cdots - a_0(x).$$

Le fait que toutes les solutions locales de $\mathcal{N}$ au point générique $t_\rho$ ont *même* rayon de convergence égal à $\lambda\rho^{\beta+1} < \omega\rho$ implique, en vertu du théorème de

Robba [$R_6$] et de Young [$Y_2$], que la norme

$$|a_m(x)(x\frac{d}{dx})^m - a_{m-1}(x)(x\frac{d}{dx})^{m-1} \cdots a_1(x)\frac{d}{dx}|_{0,\rho}$$

est strictement plus petite que la norme $|a_0|_\rho$ pour tout $\rho \in ]r_1, r_2[$.

De plus, en vertu du théorème de Young [$Y_2$] on a l'égalité

$$R(\mathcal{N}, \rho) = \omega\rho(|a_0/a_m|_\rho)^{-\frac{1}{m}} = \lambda\rho^{\beta+1},$$

en particulier

$$m - d\operatorname{Log}|a_0/a_m|(\rho) = m(\beta + 1).$$

Le lemme suivant termine la démonstration du lemme précédent.

LEMME 8.3-3.   *Soit $P(x, x\frac{d}{dx})$ un opérateur différentiel à coefficients dans $\mathcal{A}_K(]r_1, r_2[)$, une fonction a analytique dans la couronne $C(]r_1, r_2[)$ et un nombre $\rho \in ]r_1, r_2[$ tels que $|P|_{0,\rho} < |a|_\rho$ alors l'indice généralisé*

$$\widetilde{\chi}(P + a, \mathcal{A}_K([0, \rho])), \ \text{resp.} \ \widetilde{\chi}(P + a, \frac{1}{x}\mathcal{A}_K([\rho, \infty))),$$

*existe et est égal à*

$$\widetilde{\chi}(a, \mathcal{A}_K([0, \rho])) = -d\operatorname{Log}^+|a|(\rho), \ \text{resp.} \ \widetilde{\chi}(a, \frac{1}{x}\mathcal{A}_K([\rho, \infty)) = d\operatorname{Log}^-|a|(\rho).$$

*Démonstration.* Nous faisons la démonstration pour l'indice sur $\mathcal{A}_K([0, \rho])$; le cas de $\frac{1}{x}\mathcal{A}([\rho, \infty])$ se traitant de manière analogue.

Si $n$ est un entier, on a d'après la proposition 8.2-3:

$$\widetilde{\chi}(x^n(P + a), \mathcal{A}_K([0, \rho])) = -n + \widetilde{\chi}(P + a, \mathcal{A}_K([0, \rho])).$$

Posons $a^+ = \gamma^+(a)$, $a^- = \gamma^-(a)$.

Pour démontrer le lemme on peut donc remplacer $a$ par $x^n a$ et $P$ par $x^n P$. En particulier, on peut supposer que:

$$a = \sum_{n \in \mathbb{Z}} a_n x^n, \qquad |a|_\rho = \max_{n \in \mathbb{Z}}|a_n|\rho^n = |a_0|, \qquad |a^-|_\rho = \max_{n < 0}|a_n|\rho^n < |a_0|.$$

Comme $a^+ \neq 0$, l'opérateur "multiplication par $a^+$" est injectif dans $\mathcal{A}_K([0, \rho])$. Si on note pr la projection de $H_K([0, \rho[)$ dans $\mathcal{A}_K([0, \rho])$ de la décomposition de Mittag-Leffler, l'opérateur pr $\circ \frac{1}{a^+}$ en est un inverse à gauche ($1/a^+$ appartient à $\mathcal{A}_K([0, \rho[)$ mais, en général, pas à $\mathcal{A}_K([0, \rho])$).

On trouve, pour $b$ dans $\mathcal{A}_K([0, \rho])$:

$$|\operatorname{pr}(\frac{1}{a^+}\,b)|_\rho \leq |\frac{1}{a^+}\,b|_\rho = |a^+|_\rho^{-1}|b|_\rho = |a|_\rho^{-1}|b|_\rho,$$

$$|\gamma^+(P\,b)|_\rho \leq |P|_{0,\rho}|b|_\rho, \qquad\qquad |\gamma^+(a^-\,b)|_\rho \leq |a^-|_\rho|b|_\rho.$$

En notant $\| \; \|$ la norme d'opérateur sur l'espace de Banach $\mathcal{A}_K([0, \rho])$, on a donc:

$$\|\gamma^+ \circ (P + a^-)\| \, \|\mathrm{pr} \circ \frac{1}{a^+}\| \leq \max(|a^-|_\rho, \; |P|_{0,\rho}) |a|_\rho^{-1} < 1.$$

Maintenant, l'opérateur $a^+$ a un indice:

$$\chi(a^+, \mathcal{A}_K([0, \rho])) = -\mathrm{ord}_{[0,\rho]}(a^+) = -d \, \mathrm{Log}^+ |a|(\rho).$$

Le lemme est alors une conséquence immédiate de la relation suivante entre opérateurs sur $\mathcal{A}_K([0, \rho])$:

$$\gamma^+ \circ (P + a) = (1 + \gamma^+ \circ (P + a^-) \circ \mathrm{pr} \circ \frac{1}{a^+}) \circ a^+.$$

LEMME 8.3-4. *Sous les conditions précédentes on a les égalités des indices généralisés pour tout $\rho \in \, ]r_{h-1}, r_h[$:*

$$\begin{array}{rcl}
\widetilde{\chi}(\mathcal{M}, \mathcal{A}_K([0, \rho])) & = & \widetilde{\chi}(\mathcal{N}_h, \mathcal{A}_K([0, \rho^{p^h}])), \\
\widetilde{\chi}(\mathcal{M}, \frac{1}{x}\mathcal{A}_K([\rho, \infty])) & = & \widetilde{\chi}(\mathcal{N}_h, \frac{1}{x}\mathcal{A}_K([\rho^{p^h}, \infty])).
\end{array}$$

*Démonstration.* On a la décomposition en somme directe

$$\mathcal{A}_K([0, \rho]) = \oplus_{0 \leq k < p^h} x^k \mathcal{A}_K([0, \rho^{p^h}]).$$

En vertu des propositions 8.2-8 et 8.2-6 l'indice généralisé $\widetilde{\chi}(\mathcal{M}, \mathcal{A}_K([0, \rho]))$ est égale à l'indice de l'opérateur $\gamma^+(\frac{1}{p^h}x\frac{d}{dx} - G_h(x^{p^h}))$ où $G_h(x)$ est la matrice de la connexion dans une base de $\mathcal{N}_h$. Cet opérateur respecte la décomposition précédente et on trouve que

$$\widetilde{\chi}(\mathcal{M}, \mathcal{A}_K([0, \rho])) = \sum_{0 \leq k < p^h} \chi(\gamma^+(\frac{1}{p^h}x\frac{d}{dx} - G_h(x^{p^h})), x^k \mathcal{A}_K([0, \rho^{p^h}])).$$

Mais alors

$$\chi(\gamma^+(\frac{1}{p^h}x\frac{d}{dx} - G_h(x^{p^h})), x^k \mathcal{A}_K([0, \rho^{p^h}]))$$

$$= \widetilde{\chi}(\mathcal{N}_h \otimes_{\mathcal{A}_K(]r_{h-1}^{p^h}, r_h^{p^h}[)} \mathcal{A}_K(]r_{h-1}^{p^h}, r_h^{p^h}[)x^{-k/p^h}, \mathcal{A}_K([0, \rho^{p^h}])).$$

La fonction rayon de convergence du module $\mathcal{A}_K(]r_{h-1}^{p^h}, r_h^{p^h}[)x^{-k/p^h}$ est égal à $|p^h/k|\omega\rho$. Pour tout $0 < k < p^h$, elle est strictement en dessous de la fonction $R(\mathcal{N}_h, \rho)$ dans l'intervalle $]r_{h-1}^{p^h}, r_h^{p^h}[$. Ceci montre que toutes les solutions locales au point générique $t_\rho$ du module tordu $\mathcal{N}_h \otimes_{\mathcal{A}_K(]r_{h-1}^{p^h}, r_h^{p^h}[)}$ $\mathcal{A}_K(]r_{h-1}^{p^h}, r_h^{p^h}[)x^{-k/p^h}$ ont même rayon de convergence égal à $|p^h/k|\omega\rho < \omega\rho$

pour $0 < k < p^h$ et $\rho \in ]r_{h-1}^{p^h}, r_h^{p^h}[$. Mais en vertu du lemme 8.3-2 l'indice généralisé est nul:

$$\widetilde{\chi}(\mathcal{N}_h \otimes_{\mathcal{A}_K(]r_{h-1}^{p^h}, r_h^{p^h}[)} \mathcal{A}_K(]r_{h-1}^{p^h}, r_h^{p^h}[)x^{-k/p^h}, \mathcal{A}_K([0, \rho^{p^h}])) = m(1-1) = 0.$$

Un raisonnement analogue montre l'autre égalité. D'où le lemme.

Ceci entraîne par limite projective d'espaces métriques complets à images denses que pour tout $\varepsilon > 0$ assez petit les indices

$$\widetilde{\chi}(\mathcal{M}, \mathcal{A}_K([0, r - \varepsilon[)), \text{ resp. } \widetilde{\chi}(\mathcal{M}, \frac{1}{x}\mathcal{A}_K(]r - \varepsilon, \infty])),$$

existent et sont égaux à $m\beta$, resp. $-m\beta$. D'où le théorème 8.3-1.

*Remarque* 8.3-5. Nous n'avons pas montré que les indices

$$\widetilde{\chi}(\mathcal{M}, \mathcal{A}_K([0, r_h])), \text{ resp. } \widetilde{\chi}(\mathcal{M}, \frac{1}{x}\mathcal{A}_K([r_h, \infty])),$$

existent et sont égaux à $m\beta$, resp. $-m\beta$ pour les extrémités $r_h$. Ceci nécessite une précision sur les normes des changements de base qui fait passer de l'image inverse d'une base de $\mathcal{N}_h$ à une base de $\mathcal{M}$.

COROLLAIRE 8.3-6. *Soit $K$ un corps de caractéristique nulle complet localement compact pour une valuation discrète $p$-adique et $\mathcal{M}$ un $\mathcal{R}_K(r)$-module libre de rang $m$ à connexion soluble en $r$ appartenant au groupe des valeurs absolues; alors les sommets du polygone de Newton de $\mathcal{M}$*

$$\text{Newton}(\mathcal{M}, p) := (m_0, 0; m_1, m_1\beta_1; \cdots; m_N, m_N\beta_N)$$

*sont à coordonnées entières.*

*Démonstration.* En effet les nombres $m_i\beta_i$ sont des indices en vertu du théorème 8.3-1 et sont donc des entiers.

THÉORÈME 8.3-7. *Soit $K$ un corps de caractéristique nulle complet localement compact pour une valuation $p$-adique et $\mathcal{M}$ un $\mathcal{R}_K(r)$-module libre de rang $m$ à connexion soluble en $r$ appartenant au groupe des valeurs absolues dont l'exposant a la propriété ($\mathbf{NL^{**}}$), alors $\mathcal{M}$ admet des indices généralisés:*

$$\widetilde{\chi}(\mathcal{M}, \mathcal{A}_K(r)) = -\widetilde{\chi}(\mathcal{M}, \mathcal{H}_K^{\dagger}(r)) = \sum_{i>0} m_i\beta_i, \ \ \chi(\mathcal{M}, \mathcal{R}_K(r)) = 0$$

*où $\sum_{i>0} m_i\beta_i$ est la hauteur du polygone de Newton de $\mathcal{M}$.*

*Démonstration.* En vertu du théorème 6.1-14, $\mathcal{M}$ est extension de sa partie modérée $\mathcal{M}^{\leq 0}$ par sa partie de pentes strictement positives $\mathcal{M}_{>0}$. En vertu du théorème fondamental 7.1-2 la partie modérée se décompose en modules de rang un qui admettent des formes normales de la forme $x\frac{d}{dx} - \alpha$ pour des exposants $\alpha$ ayant la propriété (**NL**) et qui ont donc des indices généralisés nuls dans les espaces $\mathcal{A}_K(r)$ et $\mathcal{H}_K^\dagger(r)$. En vertu de la proposition 8.2-8 l'indice généralisé est additif dans les suites exactes, donc la partie modérée $\mathcal{M}^{\leq 0}$ admet des indices généralisés nuls. En vertu du théorème 6.2-1 la partie de pentes strictement positives admet une filtration finie dont les quotients successifs sont les modules $\mathrm{Gr}_{\beta_i}(\mathcal{M})$ purement de pentes $\beta_i > 0$. On est alors réduit au théorème d'indice local 8.3-1.

Nous pouvons étendre maintenant au cas local la définiton globale de P. Robba ([R$_4$, 10.1]) du nombre de Fuchs-Malgrange:

*Définition* 8.3-8. Soit $K$ un corps de caractéristique nulle complet localement compact pour une valuation $p$-adique et $\mathcal{M}$ un $\mathcal{R}_K(r)$-module libre de rang $m$ à connexion soluble en $r$, appartenant au groupe des valeurs absolues et dont l'exposant a la propriété (**NL**\*\*). On définit

$$\mathrm{irr}_0^r(\mathcal{M}, p) := \widetilde{\chi}(\mathcal{M}, \mathcal{A}_K(r)) = -\widetilde{\chi}(\mathcal{M}, \mathcal{H}_K^\dagger(r)) = \sum_{i>0} m_i \beta_i.$$

Le nombre $\mathrm{irr}_0^r(\mathcal{M}, p)$ apparaît comme la hauteur du polygone de Newton comme dans la théorie formelle des équations différentielles sur un corps de caractéristique nulle. Remarquons que la hauteur du polygone de Newton est toujours définie même si l'indice généralisé n'est pas fini. On trouve alors le théorème de positivité de l'irrégularité et la caractérisation de la nullité de l'irrégularité promises dans [C-M$_1$, 6.2.3]:

COROLLAIRE 8.3-9. *Soit $K$ un corps de caractéristique nulle complet localement compact pour une valuation $p$-adique et $\mathcal{M}$ un $\mathcal{R}_K(r)$-module libre de rang $m$ à connexion soluble en $r$, appartenant au groupe des valeurs absolues, dont l'exposant a la propriété (**NL**\*\*); le nombre $\mathrm{irr}_0^r(\mathcal{M}, p)$ est un entier positif ou nul et il est nul si et seulement si le module $\mathcal{M}$ est de plus grande pente nulle en $r$.*

*Démonstration.* En effet puisque les pentes $\beta_i$ sont positives, le nombre $\mathrm{irr}_0^r(\mathcal{M}, p)$ est positif on nul. Si $\mathcal{M}$ est de plus grande pente nulle en $r$, il n'a pas de pentes strictement plus grandes que 0 et donc $\mathrm{irr}_0^r(\mathcal{M}, p) = 0$. Si $\mathrm{irr}_0^r(\mathcal{M}, p) = 0$ cela entraîne que $\mathcal{M}$ n'a pas de pentes strictement plus grandes que 0, il est donc de plus grande pente nulle en $r$.

Si $\mathcal{M}$ est un module holonome de rang $m$ sur une courbe $X_K$ non singulière sur le corps $K$ et $r \leq 1$ un nombre appartenant au groupe des valeurs absolues, en vertu de [C-M$_1$, 6.2.3] quand l'indice existe, le nombre $\mathrm{irr}_{x_0}^r(\mathcal{M}, p)$ est égal à

$$\chi(\mathcal{M}, \mathcal{A}_{Kx_0}(r)) - (m - \mathrm{ord}_{x_0}^{r^-}(\mathcal{M}))$$

où $\operatorname{ord}_{x_0}^{r^-}(\mathcal{M})$ est la somme des multiplicités des composantes verticales de la variété caractéristique de $\mathcal{M}$ en chaque point singulier contenu dans le disque ouvert $D(x_0, r^-)$ de la courbe analytique. C'est le cas quand $\mathcal{M}$ a la propriété (**NL**$^{**}$) en particulier si $r = 1$ et si $\mathcal{M}$ est muni d'une structure de Frobenius dans les classes singulières.

COROLLAIRE 8.3-10. *Soient $K$ un corps de caractéristique nulle complet localement compact pour une valuation $p$-adique et $\mathcal{M}$ un $\mathcal{R}_K(r)$-module libre de rang $m$ à connexion soluble en $r$ appartenant au groupe des valeurs absolues; alors*:

1) *Si $\mathcal{M}$ a la propriété (**NL**$^*$) toute extension de sa partie modérée $\mathcal{M}^{\leq 0}$ par sa partie de pentes strictement positives $\mathcal{M}_{>0}$ est triviale. En particulier on a l'isomorphisme*:

$$\mathcal{M} \simeq \mathcal{M}^{\leq 0} \oplus \mathcal{M}_{>0}.$$

2) *Si le module des endomorphismes $\operatorname{End}_{\mathcal{R}_K(r)}(\mathcal{M}_{>0})$ de la partie de pentes strictement positives a la propriété (**NL**$^{**}$), on a la décomposition en somme directe*:

$$\mathcal{M}_{>0} \simeq \oplus_{\gamma > 0} \operatorname{Gr}_\gamma(\mathcal{M}).$$

*Démonstration.* 1) Il suffit de montrer que l'espace

$$\operatorname{Ext}^1_{\mathcal{R}_K(r)[\frac{d}{dx}]}(\mathcal{M}^{\leq 0}, \mathcal{M}_{>0})$$

est nul. Soit $\{\alpha_1, \ldots, \alpha_m\}$ l'ensemble des exposants de $\mathcal{M}$. En vertu du théorème fondamental 7.1-2 [C-M$_2$], $\mathcal{M}^{\leq 0}$ s'obtient par extensions successives de modules de rang un isomorphes aux modules $\mathcal{R}_K(r)x^{\alpha_i}$. En vertu du théorème 6.2-1, $\mathcal{M}_{>0}$ admet une filtration dont les gradués $\operatorname{Gr}_\gamma(\mathcal{M})$ sont purement de pentes $\gamma > 0$. On est réduit à montrer que $\operatorname{Ext}^1_{\mathcal{R}_K(r)[\frac{d}{dx}]}(\mathcal{R}_K(r)x^{\alpha_i},$ $\operatorname{Gr}_\gamma(\mathcal{M}))$ est nul. Mais $\operatorname{Gr}_\gamma(\mathcal{M}) \otimes_{\mathcal{R}_K(r)} \mathcal{R}_K(r)x^{-\alpha_i}$ est purement de pente $\gamma > 0$. En vertu du théorème de l'indice local 8.3-1,

$$\begin{aligned}
\chi(\mathcal{R}_K(r), &\operatorname{Gr}_\gamma(\mathcal{M}) \otimes_{\mathcal{R}_K(r)} \mathcal{R}_K(r)x^{-\alpha_i}) \\
&= \dim_K \operatorname{Ext}^1_{\mathcal{R}_K(r)[\frac{d}{dx}]}(\mathcal{R}_K(r)x^{\alpha_i}, \operatorname{Gr}_\gamma(\mathcal{M})) \\
&= \dim_K \operatorname{Hom}_{\mathcal{R}_K(r)[\frac{d}{dx}]}(\mathcal{R}_K(r)x^{\alpha_i}, \operatorname{Gr}_\gamma(\mathcal{M})) = 0,
\end{aligned}$$

parce que le morphismes horizontaux entre le module $\mathcal{R}_K(r)x^{\alpha_i}$ et le module $\operatorname{Gr}_\gamma(\mathcal{M})$ sont automatiquement nuls.

2) Soit $\beta_1 > 0$ la plus petite pente de $\mathcal{M}$ non nulle. Le module

$$\operatorname{Hom}_{\mathcal{R}_K(r)}(\mathcal{M}_{>\beta_1}, \operatorname{Gr}_{\beta_1}(\mathcal{M}_{>0}))$$

est un sous-quotient de $\mathrm{End}_{\mathcal{R}_K(r)}(\mathcal{M}_{>0})$ qui, en vertu de ([C-$M_2$, 5.4.6]) a la propriété ($\mathbf{NL^{**}}$). En vertu de la nullité de l'indice local 8.3-7 sur l'anneau $\mathcal{R}_K$ on a les égalités:

$$\dim_K \mathrm{Ext}^1_{\mathcal{R}_K(r)[\frac{d}{dx}]}(\mathcal{M}_{>\beta_1}, \mathrm{Gr}_{\beta_1}(\mathcal{M}_{>0}))$$
$$= \dim_K \mathrm{Hom}_{\mathcal{R}_K(r)[\frac{d}{dx}]}(\mathcal{M}_{>\beta_1}, \mathrm{Gr}_{\beta_1}(\mathcal{M}_{>0})) = 0$$

parce les morphismes horizontaux entre le module $\mathcal{M}_{>\beta_1}$ et le module $\mathrm{Gr}_{\beta_1}(\mathcal{M}_{>0})$ sont automatiquement nuls. En particulier $\mathcal{M}_{>0}$ est somme directe de $\mathrm{Gr}_{\beta_1}(\mathcal{M}_{>0})$ et de $\mathcal{M}_{>\beta_1}$. Mais $\mathrm{End}_{\mathcal{R}_K(r)}(\mathcal{M}_{>\beta_1})$ a la propriété ($\mathbf{NL^{**}}$) comme sous-quotient de $\mathrm{End}_{\mathcal{R}_K(r)}(\mathcal{M}_{>0})$. Une récurence sur le nombre de pentes permet de conclure.

*Remarque* 8.3-11. La structure de la partie modérée $\mathcal{M}^{\leq 0}$ est déterminée par le théorème fondamental 7.1-2 [C-$M_2$]. Le corollaire précédent ramène la structure générale d'un $\mathcal{R}_K(1)$-module différentiel soluble en 1 aux modules irréductibles purement de pentes $\gamma > 0$ (cf. [C-$M_4$]).

**8.4. *La formule de l'indice global.*** Soit $X/V$ une courbe propre connexe et lisse sur l'anneau des entiers d'une extension finie $K$ de $\mathbb{Q}_p$ de corps résiduel $k$. On note $\mathcal{X}^\dagger = (X_k, \mathcal{O}_{\mathcal{X}^\dagger/V})$, $X_K$ et $X_K^{\mathrm{an}}$ le schéma †-adique de Meredith [Mr], la courbe algébrique et la courbe analytique sur $K$ associés à $X/V$.

Sur $X/V$ on a le faisceau des opérateurs différentiels d'ordre fini $\mathcal{D}_{X/V}$ ([EGA IV, §16]) et sur $\mathcal{X}^\dagger$ on a le faisceau des opérateurs différentiels d'ordre fini $\mathcal{D}_{\mathcal{X}^\dagger/V}$ ([M-$N_1$, Appendice]).

Nous posons

$$\mathcal{O}_{X/K} := \mathcal{O}_{X/V} \otimes_{\mathbb{Z}} \mathbb{Q}, \quad \mathcal{O}_{\mathcal{X}^\dagger/K} := \mathcal{O}_{\mathcal{X}^\dagger/V} \otimes_{\mathbb{Z}} \mathbb{Q},$$
$$\mathcal{D}_{X/K} := \mathcal{D}_{X/V} \otimes_{\mathbb{Z}} \mathbb{Q}, \quad \text{et} \quad \mathcal{D}_{\mathcal{X}^\dagger/K} := \mathcal{D}_{\mathcal{X}^\dagger/V} \otimes_{\mathbb{Z}} \mathbb{Q}.$$

Soient $j : U \hookrightarrow X$ un ouvert affine de $X/V$ et $\mathcal{M}$ un $\mathcal{D}_{U/K}$-module à gauche. On lui associe le $\mathcal{D}_{\mathcal{U}^\dagger/K}$-module à gauche $\mathcal{M}^\dagger := \mathcal{O}_{\mathcal{U}^\dagger/K} \otimes_{\varepsilon^{-1}\mathcal{O}_{U_K}} \varepsilon^{-1}\mathcal{M}$ où $\varepsilon$ est le morphisme d'espace annelés $(U_k, \mathcal{O}_{\mathcal{U}^\dagger/K}) \to (U, \mathcal{O}_{U/K})$ naturel [Mr].

Rappelons que si $\mathcal{M}$ est un $\mathcal{O}_{U/K}$-module localement libre de rang $m$ on a la formule d'Euler-Poincaré de Deligne ([De$_1$, 6.21.1]):

$$\chi(U_K, \mathrm{DR}(\mathcal{M})) := \dim_K \mathrm{Hom}_{\mathcal{D}_{U/K}}(U; \mathcal{O}_{U/K}, \mathcal{M})$$
$$- \dim_K \mathrm{Ext}^1_{\mathcal{D}_{U/K}}(U; \mathcal{O}_{U/K}, \mathcal{M})$$
$$= m\chi(U_K) - \sum_{x \in \mathrm{sing}(\mathcal{M})} \mathrm{irr}_x(\mathcal{M}, \infty)$$

où $\mathrm{irr}_x(\mathcal{M}, \infty)$ est le nombre de Fuchs de $\mathcal{M}$ au point singulier $x$.

De même on définit

$$\chi(U_k, \mathrm{DR}(\mathcal{M}^\dagger)) \quad := \quad \dim_K \mathrm{Hom}_{\mathcal{D}_{\mathcal{U}^\dagger/K}}(U_k; \mathcal{O}_{\mathcal{U}^\dagger/K}, \mathcal{M}^\dagger)$$
$$- \dim_K \mathrm{Ext}^1_{\mathcal{D}_{\mathcal{U}^\dagger/K}}(U_k; \mathcal{O}_{\mathcal{U}^\dagger/K}, \mathcal{M}^\dagger).$$

On définit le saut $\phi(0)$ comme la différence

$$\chi(U_k, \mathrm{DR}(\mathcal{M}^\dagger)) - \chi(U_K, \mathrm{DR}(\mathcal{M})).$$

Posons $\mathrm{Irr}_x(\mathcal{M}^\dagger, p) := \mathrm{irr}^1_x(\mathcal{M}^\dagger, p)$. Sous l'hypothèse de l'existence de l'indice dans les classes résiduelles singulières de $\mathcal{M}$ et en vertu du théorème de semi-continuité ([C-M$_1$, 6.5.1]) le saut $\phi(0)$ est égal à la somme des dimensions des conoyaux de $\mathcal{M}$ opérant dans l'espace des fonctions analytiques dans les classes singulières. En vertu du théorème de l'indice 8.3-7 et en particulier de la finitude de l'indice on obtient:

THÉORÈME 8.4-1. *Soient $K$ un corps localement compact et $\mathcal{M}$ un $\mathcal{O}_{U/K}$-module localement libre de rang fini muni d'une connexion. Si $\mathcal{M}^\dagger$ a la propriété* (**NL**\*\*) *dans chaque classe singulière on a les égalités:*

$$\phi(0) = \sum_{x \in \mathrm{sing}(\mathcal{M}^\dagger)} \dim_K \mathrm{Ext}^1_{\mathcal{D}_{\mathcal{X}^\dagger/K}}((j_*\mathcal{M})^\dagger, \mathcal{A}_{Kx}(1))$$

*et si $\mathcal{M}^\dagger$ est soluble dans les classes singulières on a la formule d'Euler-Poincaré $p$-adique*

$$\chi(U_k, \mathrm{DR}(\mathcal{M}^\dagger)) = m\chi(U_k) - \sum_{x \in \mathrm{sing}(\mathcal{M}^\dagger)} \mathrm{Irr}_x(\mathcal{M}, p).$$

En particulier, le théorème 8.4-1 s'applique pour tous les fibrés *rationnels* sur la droite projective à connexion muni d'une structure de Frobenius dans les classes singulières. Le théorème de semi-continuité est une conséquence du théorème de dualité ([C-M$_1$], 4.4.1) et du fait que l'algèbre †-adique de $U_k$ est un espace de type $\mathcal{LF}$ où l'on peut appliquer le théorème des homomorphismes ([G$_1$, Chap. IV]). Remarquons que le caractère *algébrique* du fibré $\mathcal{M}$ est essentiel dans cette démonstration.

*Remarque* 8.4-2. 1) La cohomologie $p$-adique garde un sens pour un fibré $p$-adique $\mathcal{M}^\dagger$ ainsi que les termes locaux $\mathrm{Irr}_x(\mathcal{M}^\dagger, p)$. Le théorème d'algébrisation [C-M$_4$] ramène la formule d'Euler-Poincaré pour les fibrés $p$-adiques au cas algébrique précédent. En particulier la cohomologie $p$-adique d'un fibré à connnexion muni d'une structure de Frobenius est de *dimension finie* ce qui constitue un résultat essentiel pour la théorie $p$-adique de la fonction zêta sur les corps finis.

2) Le théorème d'existence de réseaux qui entraîne le théorème d'algébrisation [C-M$_4$] est aussi un intermédiaire essentiel pour les propriétés de finitude locales et globales pour les modules holonomes, sur le faisceau $\mathcal{D}^{\dagger}_{\mathcal{X}^{\dagger}/K}$ des opérateurs différentiels d'ordre infini sur le schéma †-adique $\mathcal{X}^{\dagger}$ ([M-N$_1$, 4.4.5], [M-N$_2$, 4.2.1]), définis dans [C-M$_4$] grâce au théorème de l'indice.

Université de Paris 6, 4 place Jussieu, F-75252 Paris
*E-mail address:* christol@math.jussieu.fr

Université de Paris 7, 2 place Jussieu, F-75251 Paris
*E-mail address:* mebkhout@math.jussieu.fr

## SIGLES